%This file starts with my \define commands,
%most of which are probably not used 
%in any given file.
%Hopefully they do not cause any trouble.

%The file uses AMSTEX and AMSppt.
%Depending on how your system is set up,
%you may need one of the following commands:

\input amstex
\documentstyle{amsppt}

\NoBlackBoxes

\magnification\magstep1

\define\p{\Bbb P}

\define\a{\Bbb A}

\redefine\c{\Bbb C}

\redefine\o{\Cal O}

\define\q{\Bbb Q}

\define\r{\Bbb R}

\define\z{\Bbb Z}

\define\n{\Bbb N}

\define\map{\dasharrow}

\define\qtq#1{\quad\text{#1}\quad}

\define\section#1{

\bigpagebreak{\smc#1}\bigpagebreak

}

\define\rup#1{\ulcorner #1\urcorner}

\define\rdown#1{\llcorner #1\lrcorner}

\define\demop{\demo{Proof}}

\define\aut{\operatorname{Aut}}

\define\pic{\operatorname{Pic}}

\define\spec{\operatorname{Spec}}

\define\depth{\operatorname{depth}}

\define\sing{\operatorname{Sing}}

\define\supp{\operatorname{Supp}}

\define\discrep{\operatorname{discrep}}

\define\tdiscrep{\operatorname{totaldiscrep}}

\define\ext{\operatorname{Ext}}

\redefine\hom{\operatorname{Hom}}

\define\codim{\operatorname{codim}}

\define\mult{\operatorname{mult}}

\define\cent{\operatorname{Center}}

\define\bs{\operatorname{Bs}}

\define\deq{:=}

\define\glafg{generically large algebraic fundamental group}

\define\norm#1{\vert\vert#1\vert\vert}

\define\nklt{\operatorname{Nklt}}

\define\broot#1{[@,@,@,\root{n}\of{#1}@,@,@,]}

\define\dual{\overset{\text{dual}}\to{\sim}}

\topmatter
\title Singularities of Pairs
\endtitle
\author J\'anos Koll\'ar
\endauthor
\address Department of Mathematics, University of Utah, Salt Lake City
UT 84112, USA  
\endaddress 
%\date  \enddate 
%\keywords{} 
%\subjclass{}
\email kollar\@{}math.utah.edu \endemail 
%FAX (801)-581-4148}}
\endtopmatter 

%Date: 10.30.1996.

%Please send all remarks, corrections, suggestions, misprints to me:

%kollar\@{}math.utah.edu

%Changes from the Jan.26,1996 version: 

%\qquad Major:  7.1--2; 

%\qquad Medium: 3.7, 3.11, 5.3.3, 7.9.2, 9.8.

\head Contents
\endhead

1. Introduction

2. Vanishing Theorems

3. Singularities of Pairs

4. Bertini Theorems

5. Effective Base Point Freeness

6. Construction of Singular Divisors

7. The $L^2$ Extension Theorem and Inversion of Adjunction

8. The Log Canonical Threshold

9. The Log Canonical Threshold and the Complex Singular Index

10. The Log Canonical Threshold and the Bernstein-Sato Polynomial

11. Rational and Canonical Singularities

\head 1. Introduction
\endhead

Higher dimensional algebraic geometry has been one of
the most rapidly developing research areas in the past twenty
years.  The first decade of its development 
 centered around
the formulation of the minimal model program and finding
techniques to carry this program through.  The proof of the existence
of flips, given in  \cite{Mori88},   completed the program in dimension three. 
These results, especially the progress leading up to \cite{Mori88}, are
reviewed in several surveys.  A very general  overview is given in
\cite{Koll\'ar87b}; many of the methods are explained in the series of lectures
\cite{CKM88}; a technically complete review for experts is found in
\cite{KaMaMa87}.

The methods of \cite{Mori88} allow us to understand three dimensional flips,
but the question of how to proceed to higher dimensions remains a baffling one.
Therefore the focus of the field moved in one of  two major directions.

\demo{Internal Developments} 

There has been a considerable internal development, as we have understood the
implications of the minimal model program to the structure of threefolds.

Two major achievements in this direction are the study of log flips by
\cite{Shokurov92}, and its subsequent application to the proof of the
abundance theorem for threefolds in \cite{Kawamata92}. Both of these results
have been simplified and explained in detail in \cite{Koll\'ar et al.92}.
Unfortunately, many of the necessary methods are rather technical
and require considerable preparatory work. A further significant advance along
these lines is the proof of  the log abundance theorem by
\cite{KeMaM$^{\text{c}}$95}. 

A     study of rational curves on varieties was undertaken in
\cite{KoMiMo92,a,b,c} and in \cite{Keel-M$^{\text{c}}$Kernan95}. Many of these
results are described in  \cite{Koll\'ar95b}. 

\cite{Alexeev93,94} studies some questions about surfaces which were inspired
by 3-dimensio\-nal problems. These results lead to a geometrically meaningful
  compactification of the moduli of surfaces of general type. 

\cite{Corti94,96} has  been developing a method to use the  
minimal model program in order to understand birational transformations between
varieties which are close to being rational.

A short overview of the above four directions  is given in \cite{Corti95}. 
\enddemo

\demo{Applications of higher dimensional geometry}

Another major theme of the last decade has been the spreading of the ideas of
the minimal model program to other fields of algebraic geometry and beyond.

One of the most dramatic changes is that people have been discovering flips
in many places. An early example is  \cite{Thaddeus94}. Later 
\cite{Dolgachev-Hu94; Thaddeus96} showed that geometric invariant theoretic
quotients are frequently related to each other by series of flips. A
similar phenomenon was discovered in \cite{Kapovich-Millson95}. The cone of
curves appeared in the study of symplectic manifolds \cite{Ruan93}. These
articles use relatively few of the results of higher dimensional geometry.
One of the reasons is that in the development of the minimal model program,
the study of singular varieties attracted the greatest attention. For the
minimal model program this is an indispensable part, but in the above
applications most varieties are smooth.

Another collection of concepts developed in higher dimensional geometry 
is a new way of looking at singularities of pairs $(X,D)$ where $X$ is a
variety and $D$  a $\q$-linear combination of divisors. 
Traditional approaches studied either the singularities of a variety $X$, or
the singularities of a divisor $D$ in a smooth variety, but  did not
concentrate on problems that occur when both $X$ and $D$ are singular.

The class of all pairs $(X,D)$ is
usually referred to as the ``log category". (Everybody is rather vague about
what the morphisms should be.)  The terminology seems to derive from the
observation that   differential forms on a variety $X$ with logarithmic poles
along a divisor $D$  should be thought of as  analogs of 
holomorphic differential forms. 
Frequently, the adjective ``log" indicates the
analog of a notion or theorem in the log category.  Unfortunately, the notion
``logarithmic pole" is not the log analog of the notion ``pole".

At the beginning, the log category was
viewed by many as  a purely technical construct, but during the last decade the
importance of this concept  gradually became indisputable. A large part of the
evidence is provided by the numerous applications of these ideas and results
in diverse questions of algebraic geometry.
\enddemo

\demo{The aim of these notes}

My intention  is to explain the basic concepts and results
of the log category, with a strong emphasis on  applications. I am
convinced that in the coming years these techniques will become an
essential tool for algebraic geometers.

These notes are written for two very different kinds of reader.  
First, they are intended to serve as a first introduction 
for algebraic geometers not familiar with higher dimensional methods.
Second, they also contain  new results and simpler proofs of
old results of interest to experts in higher dimensional geometry.
Each section starts with the introductory parts, which in themselves
constitute a coherent treatment and can be read without reference
to the more advanced parts intended for experts.  The technical
parts in each section are separated from the introductory ones by the symbol
$*****$.
\enddemo

\demo{Description of the sections}

Section 2 is a survey of the various generalizations of the Kodaira vanishing
theorem which have been developed in connection with log varieties.
Paragraphs (2.8--14) explain the simplest known proof of the basic versions. At
the end, (2.16--17) provide a summary, without proofs,  of the most general
forms of the vanishing results. The  general results can be reduced to the 
 basic versions by some technical arguments which I do not find
too illuminating. The interested reader should consult
\cite{KaMaMa87} or \cite{Koll\'ar95a}. In many instances the
general variants are easier to apply, so at least the statements should be
widely known. 

Section 3 gives the basic definitions concerning the log category. The most
important notion is   the discrepancy (3.3). This provides a measure of
how singular  a pair $(X,D)$ is. The most significant classes are defined in
(3.5). (3.6--14) give examples and various methods of computing  
discrepancies. Finally (3.18--20) relates our notion to singular metrics on
line bundles.

The first major application is in section 4. Inspired by \cite{Xu94}, we study
Bertini type theorems for linear systems with base points. One of the nicest
applications is (4.5). Its statement has nothing to do with the log category,
but its proof uses log techniques in an essential way. (4.8) gives a slew of
more technical  Bertini--type theorems in the log category.

The best developed application of log techniques is presented in section 5.
This concerns the study of linear systems of the form $K_X+L$,
where $L$ is ample. Numerous people have contributed to this direction
\cite{Fujita87; Demailly93; Ein-Lazarsfeld93; Koll\'ar93a,b; Fujita94;
Demailly94; Angehrn-Siu95; Helmke96; Kawamata96; Smith96}.  The lectures of
\cite{Lazarsfeld96} provide a very readable introduction. My aim is to
explain   a  version which works in all dimensions, but other versions give
better results in  low dimensions. 

Section 6 contains the hardest part of the proof of the results in the
previous section.  The question is quite interesting in itself: Let $L$ be
an ample divisor on a variety $X$ and $x\in X$ a point. We would like to
construct a divisor
$B$, such that $B$ is numerically equivalent to $L$ and $B$ is rather singular
at
$x$ but not too singular near $x$.  It turns out that in order to get a
reasonable answer we need to allow $B$ to be a $\q$-divisor. Also, the
traditional measures of  singularities, like the multiplicity, are not suitable
for this problem. The precise result is given in (6.4--5). For the analogous
theorems concerning singular metrics which blow up at a single point, see
\cite{Angehrn-Siu95}. 

In section 7  we compare the singularities of a pair $(X,D)$ with the pair
$(H,D|H)$ where $H\subset X$ is a  hypersurface. This problem is closely
related to the 
$L^2$ extension theorem (7.2)  of \cite{Ohsawa-Takegoshi87}.
The precise conjecture (7.3), called  ``inversion of
adjunction", was proposed by
 \cite{Shokurov92}. This trick frequently allows us to reduce an
$n$-dimensional problem to an $(n-1)$-dimensional question. For the
applications the most important variant is (7.5). This implies that the
notions {\it klt} and {\it lc}  behave well under deformations (7.6--8). A
significant application of inversion of adjunction is in the study of log
flips \cite{Shokurov92; Koll\'at et al.92, 17--18}; these results are not
discussed here.

The notion of log canonical threshold is introduced in section 8. This
concept provides a new way of measuring the singularities of pairs $(X,D)$
which do not fit in the previous framework. The most striking aspect of this
approach is a conjecture of \cite{Shokurov92} (8.8). The section is devoted
to  some computations that tend to support the conjecture.

Sections 9--10 compare the log canonical threshold to previously known
invariants of a hypersurface singularity, namely the complex singular index,
the quasiadjunction constants of \cite{Libgober83}
and the Bernstein--Sato polynomial. This raises the possibility that
conjecture (8.8) can be approached through the theory of variations of Hodge
structures or through the study of D-modules.  My hope is that experts of
these fields will get interested in such questions.

Finally, section 11 contains a simplified proof of an old result of
\cite{Elkik81; Flenner81}, asserting  that canonical singularities are
rational.  The proof is simpler in that it does not use Grothendieck's general
duality theory, but it is still not very short.
\enddemo

\demo{Terminology}

The terminology follows  \cite{Hartshorne77} for algebraic geometry. Some
other  notions, which are in general use in higher dimensional algebraic
geometry, are defined below.

\demo{1.1.1} I use the words {\it line bundle} and {\it invertible
sheaf} interchangeably. If $D$ is a {\it Cartier divisor} on a variety $X$ then
$\o_X(D)$ denotes the corresponding line bundle. {\it Linear equivalence} of
line bundles (resp. Cartier divisors) is denoted by $\cong$ (resp. 
$\sim$).  {\it Numerical equivalence} is denoted by $\equiv$.
\enddemo

\demo{1.1.2} Let $X$ be a normal variety. A {\it $\q$-divisor}  is a
$\q$-linear combination of Weil  divisors $\sum a_iD_i$.  A  $\q$-divisor
 is called {\it $\q$-Cartier} if 
it is a
$\q$-linear combination of Cartier  divisors $\sum e_iE_i$. Thus, a
$\q$-Cartier Weil divisor is a Weil divisor which is $\q$-Cartier. 
(If $X$ is smooth then any Weil divisor is also Cartier, but not in general.)
The notion of $\r$-divisor etc. can be defined analogously. 
\enddemo

\demo{1.1.3}  Let $X$ be a scheme and $D\subset X$ a divisor. A {\it log
resolution} of $(X,D)$ is a proper and birational morphism $f:Y\to X$
such that $f^{-1}(D)\cup (\text{exceptional set of $f$})$ is a divisor with
global normal crossings. Log resolutions exist if $X$ is an excellent scheme
over a field of characteristic zero. 
\enddemo

\demo{1.1.4} The {\it canonical line bundle} of a smooth variety $X$ is denoted
by $K_X$. By definition, $K_X=(\det T_X)^{-1}$ where $T_X$ is the holomorphic
tangent bundle. Thus $c_1(K_X)=-c_1(X)$.

If $X$ is a normal variety, there is a unique divisor class
$K_X$ on $X$ such that
$$
\o_{X-\sing(X)}(K_X|X-\sing(X))\cong K_{X-\sing(X)}.
$$
$K_X$ is called the {\it canonical class of $X$}. The switching between the
divisor and line bundle versions should not cause any problems.
\enddemo

\demo{1.1.5} A {\it morphism} between algebraic varieties is assumed to be
everywhere defined. It is denoted by a solid arrow $\to$. 
A {\it map}  is defined only on a dense open set. It is sometimes called a
rational or meromorphic map to emphasize this fact. It is denoted by a broken
arrow $\map$.
\enddemo

\demo{1.1.6} Let $g:U\map V$ be a map which is a morphism over the open set
$U^0\subset U$. Let
$Z\subset U$ be a subscheme such that every generic point of $Z$ is in $U^0$.
The closure of $g(U^0\cap Z)$ is called the {\it birational transform} of
$Z$. It is denoted by $g_*(Z)$. (This notion is frequently called the proper
or strict transform.)

If $f:V\map U$ is a birational map and $g=f^{-1}$ then  we get the slightly
strange looking notation $f^{-1}_*(Z)$.
\enddemo

\demo{1.1.7}  As usual, $\rdown{x}$ (round down) denotes the integral part of
a real number $x$ and
$\{x\}\deq x-\rdown{x}$ the fractional part. We also use the notation
 $\rup{x}\deq -\rdown{-x}$ (round up). 

If $D=\sum a_iD_i$ is a $\q$-divisor, then
$\rdown{D}\deq \sum \rdown{a_i}D_i$, similarly for $\{D\}$ and  $\rup{D}$. 
When using this notation we always assume that the $D_i$ are prime
divisors and $D_i\neq D_j$ for
$i\neq j$, since otherwise these operations are not well defined.
\enddemo

\demo{1.1.8} MMP stands for minimal model program, and log MMP for the 
log minimal model program. These notions are used only occasionally  and
familiarity with them is not necessary. See 
\cite{Koll\'ar87b; CKM88; KaMaMa87, Koll\'ar et al.92} for details. 
\enddemo

\enddemo

\demo{Acknowledgements} 
I would like to thank D. Abramovich, J.-P. Demailly, F. Dong, R. Lazarsfeld,
Y. Lee, L. Lempert, A. Libgober and K. Oh for their comments and corrections.  
 Partial financial support was provided by the
NSF under grant numbers  DMS-9102866 and DMS 9622394.
These notes were typeset by \AmSTeX, the \TeX\ macro system of the American
Mathematical Society.
\enddemo

\head 2. Vanishing Theorems
\endhead

Some of the most useful results in higher dimensional geometry are the
various generalizations of Kodaira's vanishing theorem. These
results are not new, and they have been surveyed ten years
ago at Bowdoin \cite{Koll\'ar87a}.  Since then we have understood
their proofs much better, and a whole new range of applications was also
discovered. The aim of this section is to explain the main ideas behind the
proof and to present a typical application (2.5--7).  Several other
applications are contained in subsequent sections. 

Throughout this section, the characteristic is zero.

For other treatments of these and related vanishing theorems see
\cite{KaMaMa87; Esnault-Viehweg92; Koll\'ar95a}. 

Let us first recall Kodaira's vanishing theorem:

\proclaim{2.1 Theorem} \cite{Kodaira53}
 Let $X$ be a smooth projective variety and  $L$  an ample line
bundle on $X$.  Then
$$
H^i(X,K_X\otimes L)=0\qtq{for $i>0$.}\qed
$$
\endproclaim

This result will be generalized in two directions.

(2.1.1) The first step 
 is to weaken the condition ``$L$ ample" while keeping all the
vanishing.    The guiding principle is that if $L$ is sufficiently close
to being ample, vanishing should still hold. 

(2.1.2) In order to understand the second step, we need to look at a
typical application of a vanishing result. Let
$$
A@>{\partial_A}>> B@>>> C@>{\partial_C}>> D
$$
be an exact sequence (in applications this is part of a long exact sequence
of cohomology groups).  If
$C=0$, then
$\partial_A$ is surjective. More generally, if $\partial_C$ is injective, then
$\partial_A$ is still surjective. Thus the injectivity of a map between
  cohomology groups can be viewed as a generalization of a vanishing
theorem.

The formulation of the first step requires some definitions.

\demo{2.2 Definition}  Let $X$ be a proper variety and $L$  a line bundle
or a Cartier divisor on
$X$.

(2.2.1) $L$ is {\it nef} iff $\deg_C(L|C)\geq 0$ for every irreducible curve
$C\subset X$.

(2.2.2) $L$ is {\it big} iff $H^0(X,L^m)$ gives a birational map $X\map
X'\subset
\p$ to some projective space for $m\gg 1$. Equivalently, $L$ is big iff
$h^0(X,L^m)>const\cdot m^{\dim X}$ for $m\gg 1$. 

(2.2.3) Both of these notions extend by linearity to $\q$-Cartier divisors.
\enddemo

We are ready to formulate the simplest form of the general vanishing theorem 
about perturbations of ample line bundles:

\proclaim{2.3 Theorem} \cite{Kawamata82; Viehweg82}
Let $X$ be a smooth projective variety and  $L$ 
a line bundle on $X$. Assume that we can write $L\equiv M+\sum d_iD_i$
where $M$ is a nef and big $\q$-divisor, $\sum D_i$ is a normal crossing
divisor and
$0\leq d_i<1$ for every $i$.
 Then
$$
H^i(X,K_X\otimes L)=0\qtq{for $i>0$.}
$$
\endproclaim

\demo{2.4 Questions} The two main questions raised by this result are the
following:

(2.4.1.1) Where do we find line bundles $L$ which can be decomposed as
$L\equiv M+\sum d_iD_i$?

(2.4.1.2) How to use a result like this?
\enddemo

There are two basic situations where line bundles with a decomposition
$L\equiv M+\sum d_iD_i$ naturally arise.

(2.4.2.1) Let $S$ be a normal surface   and $f:S'\to S$ a resolution of
singularities. Assume for simplicity that $S$ has rational singularities. 
Then $f^*\pic(S)$ and the exceptional curves $D_i$ of
$f$   generate a finite index subgroup of  $\pic(S')$. Thus any
line bundle $L$ on $S'$ can be decomposed as
$L\equiv f^*M+\sum d_iD_i$ where $M$ is a $\q$-Cartier divisor on $S$ and the
$D_i$ are the $f$-exceptional curves. Rational coefficients are usually
inevitable.

(2.4.2.2) Let $L$ be a line bundle on a variety $X$ such that $L^n$ has a
section with zero set $D=\sum a_iD_i$. Then $L\equiv \sum (a_i/n)D_i$.
In general $\sum D_i$ is not a normal crossing divisor. By Hironaka,
there is a proper birational morphism $f:X'\to X$
such that $f^*L\equiv \sum b_iB_i$, where $\sum B_i$ is a normal crossing
divisor. 

This suggests that the best hope of using (2.3) is in questions which are
birational in nature. 
The following example shows a rather typical application of this kind.

\demo{2.5 Example} 
Let  $X$ be a smooth proper variety of general type. Our aim is to 
 express $H^0(X,sK_X)$ as an Euler characteristic, at least for $s\geq 2$. 

(2.5.1) First approach. 

Assume that $X$ is a surface and let $X'$ be its
minimal model. Then $K_{X'}$ is nef and big, so
$$
H^0(X,sK_X)=H^0(X',sK_{X'})=\chi(X',sK_{X'})\qtq{for $s\geq 2$.}
$$
In higher dimensions this still works if  $X$ has a minimal model $X'$.
Unfortunately,  the existence of a minimal model is unknown in general.

(2.5.2) Second approach.

\noindent  We try to follow (2.4.2.2) with some modifications. 
Choose an auxiliary number $m\gg 1$ and $f:X'\to X$ such that:

$X'$ is smooth, and

$|mK_{X'}|=|M|+F$, where 
$M$ is free and big, and
$F=\sum a_iF_i$ is a normal crossing divisor. 

\noindent This is always  possible by Hironaka. Further conditions on $m$ will
be imposed later. Thus
$$
sK_{X'}=K_{X'}+(s-1)K_{X'}\equiv K_{X'}+\frac{s-1}{m}M+ \frac{s-1}{m}F.
$$
This is not quite what we want, since   $\frac{s-1}{m}F$ can have 
coefficients that are bigger than one. To remedy this problem we just get rid
of the excess coefficients in $\frac{s-1}{m}F$. We want to get a Cartier
divisor, so we can subtract only integral multiples of Cartier divisors, and we
also want to end up with coefficients between 0 and 1. These two conditions
uniquely determine the choice of
$\sum c_iF_i\deq \sum
\rdown{(s-1)a_i/m}F_i$.
Set
$$
\aligned
L\deq& (s-1)K_{X'}-\sum c_iF_i\\
\equiv& \frac{s-1}{m}M+ \sum \left\{\frac{(s-1)a_i}{m}\right\}F_i.
\endaligned
\tag 2.5.3
$$
The choice of $L$ in (2.5.3) is set up so that vanishing applies to
$K_{X'}+ L$. We still need to check that we have not lost any sections
of $\o(sK_{X'})$ by subtracting $\sum c_iF_i$.

Let $D=D'+\sum b_iF_i$ be any divisor in the linear system $|sK_{X'}|$.
Assume that $m=rs$. Then
$rD'+\sum rb_iF_i\in |mK_{X'}|$, thus $rb_i\geq a_i$ for every $i$. Therefore
$$
b_i\geq a_i/r=sa_i/m\geq\rdown{(s-1)a_i/m}=c_i.
\tag 2.5.4
$$
Thus  $D-\sum c_iF_i=D'+\sum (b_i-c_i)F_i$ is effective. Therefore
$$
H^0(X,sK_X)=H^0(X',sK_{X'})=H^0(X',K_{X'}+L).
$$
 Thus by (2.3)
$$
H^0(X,sK_X)=\chi(X',K_{X'}+L).
$$
Note. The choice of $m$ has been left rather free. Different choices
do lead to different models $X'$. Also,  the estimate (2.5.4) is far from
being sharp. In delicate situations it is worthwhile to check how much room
it gives us.
\enddemo

The following theorems  use this construction to compare plurigenera of
\'etale covers:

\proclaim{2.6 Theorem}  \cite{Koll\'ar95a, 15.4} Let $p:Y\to X$ be a finite
\'etale morphism between smooth and proper varieties of general type.  Then 
$$
h^0(X,sK_X)=\frac1{\deg p}h^0(Y,sK_Y ) \qtq{for $s\geq 2$.}
$$
\endproclaim

\demop Fix $s$ and choose  $f: X'\to X$ as in (2.5.2). Set $Y'=Y\times_XX'$. 
Then $p':Y'\to X'$ is \'etale. As in (2.5.2) we construct $L$ such that,
$$
h^0(X,sK_X )=\chi(X', K_{X'}+ L),
\qtq{and} 
h^0(Y,sK_Y )=\chi(Y', K_{Y'}+ {p'}^*L).
$$
(This requires a little extra care; see 
\cite{Koll\'ar95a, 15.5} for details.)
An Euler characteristic is multiplicative in \'etale covers, thus
$$
\chi(X', K_{X'}+ L)=\frac1{\deg p}\chi(Y', K_{Y'}+ {p'}^*L).
$$
The two formulas together imply (2.6). 
\qed\enddemo

This is just the baby version of the following result which compares
plurigenera in possibly infinite covers. See \cite{Koll\'ar95a, 15.5}
for the necessary definitions and the proof.

\proclaim{2.7 Theorem} 
Let $X$ be a  smooth, proper variety of general type and 
$p:Y\to X$  a (possibly infinite)  \'etale Galois
cover with Galois group $\Gamma$.
 Then 
$$
h^0(X,K_X^m)=\dim_{\Gamma}H^0_{(2)}(Y,K_Y^m) \qtq{for $m\geq 2$.}
$$
(Here $H^0_{(2)}$ is the Hilbert space of holomorphic $L^2$ sections with
respect to a metric pulled back from $X$ and $\dim_{\Gamma}$ is the usual
dimension in the theory of von Neumann algebras.)
\endproclaim

\demo{2.8 Idea of the proof of the vanishing theorems}

The two main steps of the proofs, as outlined in \cite{Koll\'ar86b,Sec.5}, are
the following:

(2.8.1) Step 1.

Find several examples where  the coherent cohomology of a
sheaf comes from topological cohomology. The simplest example of this
situation is given by  Hodge theory. For the proof see
any of the standard textbooks on K\"ahler geometry (e.g.
\cite{Wells73, V.4.1;  Griffiths-Harris78, p.116}).

\proclaim{2.8.1.1 Theorem} 
Let $X$ be a smooth proper variety 
(or compact K\"ahler manifold) with structure sheaf $\o_X$. Let
$\c_X\subset \o_X$ denote the constant sheaf. Then the natural map
$$
H^i(X,\c_X)\to H^i(X,\o_X)\qtq{is surjective for every $i$.}\qed
$$
\endproclaim

We need that this also holds if $X$ has quotient
singularities. This is due to \cite{Steenbrink77;
Danilov78}. The more up-to-date ``orbifold approach" is to 
 notice that the usual
proof for manifolds works with essentially no changes. We should still view
$X$ as being patched together from smooth coordinate charts, but instead of
allowing patching data between different charts only, we admit patching data
between a chart and itself, corresponding to the local group action. Once the
conceptual difficulties are behind, the proof is really the same.
\enddemo

\demo{2.8.1.2 Remarks} 

The analog of (2.8.1.1) also holds if $X$ has rational
singularities, but I do not know any simple proof, cf. 
\cite{Koll\'ar95a, Chap.12}.

More generally, any variation of Hodge structures over $X$ gives rise to a
similar situation, see \cite{Koll\'ar86b,Sec.5; Saito91}. 
\enddemo

(2.8.2) Step 2.

By an auxiliary construction, which in this case is the study of cyclic
covers, we find many related situations of a topological sheaf $\bold F$ and a
coherent sheaf $\Cal F$ together with natural surjections
$$
H^i(X,\bold F)\to H^i(X,\Cal F).
$$
(By ``topological sheaf" I mean a sheaf of abelian groups which is defined in
terms of the classical topology of $X$.)  Moreover, we try to achieve that
$\bold F$ is supported over an open set
$X^0\subset X$. If $X^0$ is sufficiently small, for instance affine, then
many cohomologies vanish over $X^0$, giving the vanishing of certain  coherent 
cohomology groups.

(In the treatment of \cite{Esnault-Viehweg92}, a De Rham complex takes the
place of  the topological sheaf $\bold F$.)

We start the proof of (2.3) by
constructing cyclic covers and   studying their basic properties. This is the
third time that I give a somewhat different treatment of cyclic covers
(cf. \cite{Koll\'ar95a,Chap.9; Koll\'ar95b,II.6}) but I am still unable to find
one which  I consider optimal.

\demo{2.9 Local construction of cyclic covers}  

Let $U$ be a smooth variety,
$f\in
\o_U$ a regular function and $n\geq 1$ a natural number. 
Let $D=(f=0)$ be the zero divisor of $f$. 
We want to construct the
cyclic cover corresponding to  $\root{n}\of{f}$. It is denoted by
$U\broot{f}$. 

Let $y$ be a new variable.  $U\broot{f}\subset U\times \a^1$ is defined
by the equation $y^n=f$.  Thus $U\broot{f}$ is smooth at a point $(u,\ast)$ iff
either $f(u)\neq 0$ or  $u$ is a smooth point of the divisor $D$.

 Let
$p:U\broot{f}\to U$ be the coordinate projection.   $p$ is finite over $U$
and \'etale over $U-D$. The cyclic
group $\z_n$ acts on $U\broot{f}$ and the $\z_n$-action gives an eigensheaf
decomposition
$$
p_*\o_{U\broot{f}}=\o_U+y\o_U+\dots+y^{n-1}\o_U.
$$
Let $\bar U\broot{f}\to U\broot{f}$ be the normalization and $\bar p:\bar
U\broot{f}\to U$ the projection. The $\z_n$-action lifts to a $\z_n$-action
on $\bar U\broot{f}$, thus we get an eigensheaf
decomposition
$$
\bar p_*\o_{\bar U\broot{f}}=\sum_{k=0}^{n-1}F_k,\qtq{where} F_k\supset
y^k\o_U.
$$
Each $F_i$ is a rank one reflexive sheaf, hence invertible since $U$ is
smooth. Thus there are divisors $D^k$ such that  $F_k=y^k\o_U(D^k)$ and
$\supp D^k\subset
\supp D$. 

In order to figure out the coefficients in $D^k$, we may localize at a smooth
point
$u\in\supp D$. Thus we can assume that $f=vx_1^d$ where $v$ is a unit at
$u$ and $x_1$ is a local coordinate at $u$. 

A rational section $y^kx_1^{-j}$ of $y^k\o_U$ is integral over $\o_U$ iff
$$
(y^kx_1^{-j})^n=y^{nk}x_1^{-nj}=f^kx_1^{-nj}=v^kx_1^{kd-nj}
$$
 is a regular function. Thus $j\leq dk/n$. This shows that
$$
\bar p_*\o_{\bar U\broot{f}}=\sum_{k=0}^{n-1}y^k\o_U(\rdown{kD/n}).
$$
\enddemo

\demo{2.10 Local structure of $\bar U\broot{f}$} 

Assume that $D$ is a normal crossing divisor. Pick $u\in U$ and let 
$D_1,\dots,D_s$ be the irreducible components of $D$ passing through $u$. 
Choose local  coordinates $x_i$ at $u$ such that
$D_i=(x_i=0)$. Let  $B\subset
U$ be a polydisc around $u$ defined by $|x_i|<1$ for every $i$. 

$\pi_1(B-D)\cong \z^s$ is generated by the loops around the
divisors $D_i$. 

Let $\bar B\subset \bar p^{-1}(B)$ be an irreducible component. 
$\bar B\to B$ is \'etale over $B-D$, thus  it corresponds to
a  quotient  $\z^s\to \z_n$. 
 By Galois theory,  $\bar B$ is a quotient of the
cover corresponding to the subgroup $(n\z)^s\subset \z^s$.

Let $\Delta\subset \c$ be the unit disc. The cover
corresponding to  $(n\z)^s\subset \z^s$ is
$$
\Delta^m\to B,\qtq{given by} (z_1,\dots,z_m)\mapsto
(z_1^n,\dots,z_s^n,z_{s+1},\dots,z_m).
$$
 This cover is
smooth, hence $\bar U\broot{f}$ has only quotient singularities.
\enddemo

\demo{2.11 Global construction of cyclic covers} 

Let $X$ be a smooth variety, $L$  a line bundle on $X$ and $s\in H^0(X,L^n)$
a section.

Let $U_i\subset X$ be an affine cover such that  $L|U_i$ has a nowhere zero
global section $h_i$. $L$ is given by transition functions $h_i=\phi_{ij}h_j$.
Let $s=f_ih_i^n$. $s$ is a well defined section, thus
$f_i=\phi_{ij}^{-n}f_j$.

The local cyclic covers are given by equations $y_i^n=f_i$. These are
compatible if we set the transformation rules $y_i=\phi_{ij}^{-1}y_j$.
This gives the global cyclic cover
$$
 X\broot{f}=\cup_i  U_i\broot{f_i}.
$$
The invertible sheaves $y_i^k\o_{U_i}$  patch together to the line
bundle $L^{-k}$, and so
$$
p_*\o_{X\broot{f}}=\o_X+L^{-1}+\dots+L^{-(n-1)}.
$$
Let $\bar X\broot{f}\to X\broot{f}$ be the normalization and $\bar p:\bar
X\broot{f}\to X$ the projection. The $\z_n$-action gives the eigensheaf
decomposition
$$
\bar p_*\o_{\bar X\broot{f}}=\sum_{k=0}^{n-1}L^{-k}(\rdown{kD/n}).
$$
\enddemo

\demo{2.12 Decomposing $\bar p_*\c$}

Until now everything worked in the Zariski as well as in the classical
topology. From now on we have to use the classical topology.

In order to simplify notation set $Z=\bar X\broot{f}$ and let 
 $\c_Z\subset \o_Z$ denote the  sheaf of locally constant functions. We have an
eigensheaf decomposition
$$
\bar p_*\c_Z=\sum G_k, \qtq{such that} G_k\subset L^{-k}(\rdown{kD/n}).
$$
It is not  hard to write down the sheaves $G_k$ explicitly (cf.
\cite{Koll\'ar95a,9.16}), but the  arguments are clearer  and simpler if
we do not attempt to do this. Their basic cohomological  properties are  easy
to establish:
\enddemo

\proclaim{2.12.1 Proposition} Notation as above. Write  $D=\sum d_iD_i$.

(2.12.1.1) $G_0\cong \c_X$.

(2.12.1.2) For every $x\in X$ 
there is an open neighborhood
$x\in U_x\subset X$ such that $H^i(U_x,G_j|U_x)=0$ for
$i>0$.

(2.12.1.3) If $U\subset X$ is connected,
 $U\cap D\neq \emptyset$ 
and $n$ does not divide $d_i$ for every $i$, 
then $H^0(U,G_1|U)=0$.
\endproclaim
 
\demop $G_0$ is the invariant part of $\bar p_*\c_Z$, which is  $\c_X$.

Choose $x\in U_x\subset X$  such that $V_x=(\bar p)^{-1}U_x\subset Z$ 
retracts to $(\bar p)^{-1}(x)$. Then $H^i(V_x,\c_{V_x})=0$ for $i>0$. Since
$\bar p$ is finite,
$$
H^i(W,\bar p_*\c_Z|W)=H^i((\bar p)^{-1}W,\c_Z|(\bar p)^{-1}W),
$$
for every $W\subset X$. In particular for $W=U_x$ we obtain that
$H^i(U_x,\bar p_*\c_Z|U_x)=0$. This implies (2.12.1.2) since
$H^i(U_x, G_j|U_x)$ is a direct summand of
$H^i(U_x,\bar p_*\c_Z|U_x)$.

If $U$ intersects $D$ then we can find a point
$u\in U\cap D$ which has a neighborhood where  $D$
is defined by a  function $f=vx_1^d$ where $v$ is a unit and $x_1$ a local
coordinate.

In the local description of (2.9),  any rational section of $y\o_U$ can
be written as $yg$ where $g$ is a rational function. $yg$ is a locally constant
section iff $(yg)^n=fg^n$ is a locally constant
section of $\o_U$. 
That is, when
$$
g=cf^{-1/n}=cv^{-1/n}x^{-d/n}\qtq{for some $c\in \c$.}
$$
Since $n$ does not divide $d$, this gives a rational function only for $c=0$.
\qed\enddemo

 The following result  is a very general theorem about the injectivity of
certain maps between cohomology groups.  In (2.14) we see that it
implies (2.3), at least when $M$ is ample. The general case of
(2.3) requires a little more work. 

In some applications the injectivity part is important \cite{Koll\'ar86a,b;
Esnault-Viehweg87}, though so far the vanishing theorem has found many more
uses.

The theorem is a culmination of the work of several authors
\cite{Tankeev71; Ramanujam72; Miyaoka80; 
Kawamata82; Viehweg82; Koll\'ar86a,b,87a;
Esnault-Vieh\-weg86,87}.

\proclaim{2.13 Theorem} Let $X$ be a smooth proper variety and 
 $L$  a line bundle on
$X$. Let $L^n\cong \o_X(D)$ where $D=\sum d_iD_i$ is an effective divisor.
Assume that $\sum D_i$ is a normal crossing divisor and $0< d_i<n$ for every
$i$. Let
$Z$ be the normalization of
$X\broot{s}$. Then:

(2.13.1) \
$H^j(Z,\c_Z)\to H^j(Z,\o_Z)$ is surjective for every $j$.

(2.13.2) \
$H^j(X,G_1)\to H^j(X,L^{-1})$ is surjective for every $j$.

(2.13.3) For every $j$ and $b_i\geq 0$ the natural map
$$
H^j(X,L^{-1}(-\sum b_iD_i))\to H^j(X,L^{-1})\qtq{is surjective.}
$$

(2.13.4) For every $j$ and $b_i\geq 0$ the natural map
$$
H^j(X,K_X\otimes L)\to H^j(X,K_X\otimes L(\sum b_iD_i))\qtq{is
injective.}
$$
\endproclaim

\demop  By (2.10) $Z$ has only quotient singularities, thus (2.13.1)
follows from (2.8.1.1). 

 The assumption $0<d_i<n$ implies that
$\rdown{D/n}=\emptyset$, thus
$F_1=L^{-1}$.  Therefore
(2.13.2) is just (2.13.1) restricted to one
$\z_n$-eigenspace.

(2.13.3) and (2.13.4) are equivalent by Serre duality, thus it is
sufficient to prove (2.13.3).

The main step is the following:

\proclaim{2.13.5 Claim} 
$G_1$ is a subsheaf of $L^{-1}(-\sum b_iD_i)$.
\endproclaim

\demop Both of these are subsheaves of $L^{-1}$, so this is  a local question.
We need to show that if $U\subset X$ connected, then
$$
H^0(U,G_1|U)\subset H^0(U,L^{-1}(-\sum b_iD_i)|U).\tag 2.13.6
$$
 $L^{-1}(-\sum b_iD_i)$ and $L^{-1}$   are equal over $X-D$, thus (2.13.6)
holds if $U\subset X-D$.  If $U\subset X$ is connected and it intersects $D$
then  by (2.12.1.3) $H^0(U,G_1|U)=\{0\}$, thus  
$$
H^0(U,G_1|U)=\{0\}\subset H^0(U,L^{-1}(-\sum b_iD_i)|U)\qtq{trivially.}
\qed
$$
\enddemo

This gives a factorization
$$
H^j(X,G_1)\to H^j(X,L^{-1}(-\sum b_iD_i))\to H^j(X,L^{-1}).
$$
The composition is surjective by (2.13.2) hence the second arrow is also
surjective.
\qed\enddemo

As a corollary, let us prove (2.3) in a special case:

\demo{2.14 Proof of (2.3) for  $M$ ample} 

We have $L\equiv M+\sum_{i\geq 1}d_iD_i$. 
Choose $n$ such that   $nd_i$ is an integer for every $i$ and
$$
nM\equiv nL-\sum_{i\geq 1}(nd_i)D_i\sim D_0, 
$$
 where $D_0$ is a smooth divisor
which intersects $\sum D_i$ transversally. 
(This is possible since $M$ is ample.)
Let $d_0=1/n$. 
Then $nL\sim \sum_{i\geq 0}(nd_i)D_i$, thus 
by (2.13.4),
$$
H^i(X,K_X+L)\to H^i(X,K_X+L+ b_0D_0)\qtq{is injective.}
$$
By Serre vanishing, the right hand side is zero for $b_0\gg 1$ and $i>0$.

Thus $H^i(X,K_X+L)=0$ for $i>0$.\qed\enddemo

As an exercise in using (2.3), derive the following relative version of
it.

\demo{2.15 Exercise}  Let $X$ be a smooth projective variety and  $L$ 
a line bundle on $X$. Assume that we can write $L\equiv M+\sum d_iD_i$
where $M$ is a nef and big $\q$-divisor, $\sum D_i$ is a normal crossing
divisor and
$0\leq d_i<1$ for every $i$.

Let $f:X\to Y$ be a proper and birational morphism. 
 Then
$$
R^jf_*(K_X\otimes L)=0\qtq{for $j>0$.}
$$
\enddemo

For the applications it is frequently useful that we have a vanishing even if
$\sum d_iD_i$ is not a normal crossing divisor. This approach was first used 
extensively by
\cite{Nadel90} in the analytic setting which we discuss in (3.18--20).  (See
(3.5) for the definition of klt.)

\proclaim{2.16 Theorem}  
Let $X$ be a 
normal and proper variety and $N$ a  line bundle  on $X$.
 Assume that $N\equiv K_X+\Delta+M$ where $M$ is a  nef and big $\q$-Cartier
$\q$-divisor on
$X$ and
$\Delta$ is an effective $\q$-divisor.
Then there is an ideal sheaf $J\subset \o_X$ such that

(2.16.1) $\supp (\o_X/J)=\{x\in X| (X,\Delta) \text{ is not klt at $x$}\}$.

(2.16.2) 
$H^j(X,N\otimes J)=0$ for $j>0$.
\endproclaim

\demop  Let $f:Y\to X$ be a log resolution of $(X,\Delta)$. Write
$$
K_Y\equiv f^*(K_X+\Delta)+\sum a_iE_i.
$$
This can be rewritten as
$$
\align
f^*N+\sum a_iE_i&\equiv K_Y+f^*M,\qtq{or as}\\
f^*N+\sum \rup{a_i}E_i&\equiv K_Y+f^*M+\sum (\rup{a_i}-a_i)E_i.
\endalign
$$
By (2.3) and (2.15) we know that
$$
H^j(Y,f^*N(\sum \rup{a_i}E_i))=0\qtq{and} 
R^jf_*(f^*N(\sum \rup{a_i}E_i))=0\qtq{for $j>0$.}
$$
Thus by the Leray spectral sequence,
$$
H^j(f_*(f^*N(\sum \rup{a_i}E_i)))=0\qtq{for $j>0$.}
$$
By the projection formula,
$$
f_*(f^*N(\sum \rup{a_i}E_i))=N\otimes f_*(\o_Y(\sum \rup{a_i}E_i)).
$$
If $E_i$ is not an exceptional divisor, then $a_i\leq 0$. Thus in 
$\sum \rup{a_i}E_i$ only $f$-exceptional divisors appear with positive
coefficient.  Therefore
$$
f_*(\o_Y(\sum \rup{a_i}E_i))=f_*(\o_Y(\sum_{a_i\leq -1} \rup{a_i}E_i)).
$$
The latter is an ideal sheaf in $\o_X$ whose cosupport is exactly
the set of points over which there is a divisor with  
$a_i\leq -1$.\qed\enddemo

The following is a summary of the most general versions of vanishing theorems.
Proofs can be found in \cite{KaMaMa87; Esnault-Viehweg92; Koll\'ar95a}.
For the latest results in the holomorphic category, see \cite{Takegoshi95}. 

\proclaim{2.17 Theorem} (General Kodaira Vanishing)  

Let $X$ be a 
normal and proper variety and $N$ a  rank one, reflexive 
sheaf  on $X$.
 Assume that $N\equiv K_X+\Delta+M$ where $M$ is a   $\q$-Cartier
$\q$-divisor on
$X$ and
$(X,\Delta)$ is klt. 

(2.17.1) (Global vanishing)

Assume that $M$ is nef and big.
Then
$$
H^i(X,N)=0\qtq{for $i>0$.}
$$

(2.17.2)  (Injectivity theorem)

Assume  that $M$ is nef   and  let $D, E$ 
be   effective   Weil divisors on $X$ such that $D+E\equiv mM$ for some $m>0$.
Then 
$$
H^j(X,N)\to H^j(X,N(D))\qtq{is injective for $j\geq 0$.}
$$

(2.17.3) (Relative vanishing)

Let $f:X\to Y$ be a surjective morphism with generic fiber $X_{gen}$.
Assume that $M$ is $f$-nef and $M|X_{gen}$ is big. Then
$$
R^jf_*N=0 \qtq{  for $j> 0$.}
$$
(Note that if 
 $f$ is generically finite then
$M|X_{gen}$ is always big.)

(2.17.4)  (Torsion freeness)

Let $f:X\to Y$ be a surjective morphism.  Assume that $M\equiv f^*M_Y$ where
$M_Y$ is a  $\q$-Cartier
$\q$-divisor on
$Y$. Then

$R^jf_*N$ is torsion free for $j\geq 0$.
\endproclaim

\demo{2.17.5 Remark} One has to be careful
with the localization of (2.17.4). Namely, (2.17.4) is not true if $f$ is
a projective morphism but $Y$ is not proper. An example is given by the
Poincar\'e  sheaves on Abelian varieties. The local version holds if the
numerical equivalence is everywhere replaced by $\q$-linear equivalence.
\enddemo

\demo{2.17.6 Remark}  The  Grauert-Riemenschneider vanishing theorem is  the
special case of  (2.17.3)  when $X$ is smooth and $N=K_X$.
\enddemo

\head 3. Singularities of Pairs
\endhead

There are many ways to measure how singular a variety is. In higher dimensional
geometry  a new notion, called the discrepancy, emerged. This concept  was
developed to deal with the following two situations:

(3.1.1) Let $X$ be a singularity and $f:Y\to X$ a resolution. We would like to
measure  how singular $X$ is by comparing 
 $K_Y$   with $f^*K_X$, as $Y$
runs through all resolutions.

(3.1.2) Let $D\subset \c^{n+1}$ be a hypersurface with defining equation
$h=0$. Let $f:Y\to \c^{n+1}$ be a birational morphism, $Y$ smooth. We
would like to
measure  how singular $D$ is by comparing the order of vanishing of the
Jacobian of $f$ with the order of vanishing of $f^*h$ along exceptional
divisors, as $Y$ runs through all birational morphisms.

In order to illustrate the final definitions, it is very useful to look at the
following examples. For simplicity we consider the first of the above cases.

\demo{3.2 Example}  Let $X$ be  normal and assume that   $mK_X$ is Cartier. Let
$f:Y\to X$ be a   birational morphism, $Y$ smooth. We can write 
$$
mK_Y= f^*(mK_X)+\sum (ma_i)E_i,
$$
where the
$E_i$ are  exceptional divisors and the $a_i$ are rational.  We frequently
write this in the form
$$
K_Y\equiv f^*K_X+\sum a_iE_i.
$$
Our aim is to get a single invariant out of the numbers $a_i$, preferably one
that is independent of the resolution. The straightforward candidates are
$\min\{a_i\}$ and  $\max\{a_i\}$. One can easily see that the latter depends on
$f:Y\to X$, and its supremum as $Y$ changes is always $+\infty$.

$\min\{a_i\}$ is somewhat better behaved, but it still depends on the choice of
$f:Y\to X$ in most cases. To make things  better, assume that 
$\sum E_i$ is a normal crossing divisor.

 Let $Z\subset Y$ be a smooth
subvariety and $p:B_ZY\to Y$ the blow up with  exceptional divisor $F\subset
B_ZY$. Set
$g\deq f\circ p:B_ZY\to X$. Then
$$
K_{B_ZY}\equiv g^*K_X+cF+\sum a_iE'_i,
$$
where $E'_i$ is the birational transform of $E_i$ on $B_ZY$.

If $Z$ is not contained in $\sum E_i$, then $c\geq 1$. Otherwise
it is not hard to check that
$$
c\geq \min\{a_i\} \qtq{if} \min\{a_i\}\geq -1,
$$
but not in general.
In fact, if  $\min\{a_i\}<-1$ then there is a  sequence of  resolutions such
that $\min\{a_i\}\to -\infty$ (3.4.1.4). In this case we say that  $X$ is not
log canonical. For these singularities our invariant does not
give anything. From the point of view of general singularity theory this is
rather unfortunate, since most singularities are not log canonical. 
(In section 8 we introduce another invariant which is meaningful in the non
log canonical case.)

Our point of view is, however, quite different. Our main interest is in smooth
varieties, and we want to deal with singularities only to the extent they
inevitably appear in the course of the minimal model program. In many situations
it  is precisely
$\min\{a_i\}$ which tells us which varieties need to be considered.

If $\min\{a_i\}\geq -1$, then by (3.13) the minimum is independent of
the choice of $f:Y\to X$ (assuming that $Y$ is
smooth and $\sum E_i$ is a normal crossing divisor). 
\enddemo

More generally, one can put the two aspects mentioned in (3.1.1--2)  together,
and consider pairs $(X,D)$ where $X$ is a normal scheme and $D$ a formal linear
combination of Weil divisors $D=\sum d_iD_i$, $d_i\in \r$. It took people about
10 years to understand that this is not simply a technical  generalization but
a very fruitful ---  even  basic --- concept. For now just believe
that this  makes sense.

Since one cannot pull back arbitrary Weil divisors, we always have to assume
that $K_X+D$ is $\r$-Cartier, that is, it is an $\r$-linear combination of
Cartier divisors. In the applications we almost always use $\q$-coefficients,
but for the basic definitions the coefficients  do not matter.

The resulting notion  can be related to more traditional ways of measuring
singularities (for instance, multiplicity or arithmetic genus), but it
is a truly novel way of approaching the study of singularities.

\demo{3.3 Definition}  Let $X$ be a normal, integral scheme and $D=\sum
d_iD_i$ an
$\r$-divisor (not necessarily effective) such that $K_X+D$ is $\r$-Cartier. Let
$f:Y\to X$ be a   birational morphism, $Y$ normal. We can write 
$$
K_Y\equiv f^*(K_X+D)+\sum a(E,X,D)E,
\tag 3.3.1
$$
where $E\subset Y$ are  
distinct prime divisors and $a(E,X,D)\in \r$. The right hand side is not
unique because we allow nonexceptional divisors in the summation. In order to
make it unique we adopt the following:

\demo{3.3.2 Convention}  A nonexceptional divisor
$E$ appears in the sum  (3.3.1) iff $E=f^{-1}_*D_i$ for some $i$, and then with
the coefficient
$a(E,X,D)=-d_i$. (Note the negative sign!)

Similarly, if we write $K_Y+D'\equiv f^*(K_X+D)$, then $D'=-\sum a(E,X,D)E$.
\enddemo

$a(E,X,D)$ is called the {\it discrepancy} of $E$ with respect to $(X,D)$.  
We   frequently write $a(E,D)$ or $a(E)$ if no confusion is likely.

If $f':Y'\to X$
is another  birational morphism and $E'\subset Y'$ is the birational
transform  of $E$ on $Y'$ then $a(E,X,D)=a(E',X,D)$.  In this sense $a(E,X,D)$
   depends only
on the divisor $E$ but not on $Y$.  This is the reason why $Y$ is suppressed in
the notation. 

A more invariant description is obtained by considering a rank
one discrete valuation $\nu$ of the function field $K(X)$. $\nu$  corresponds
to a divisor $E\subset Y$ for some $f:Y\to X$.  The closure of $f(E)$ in $X$ is
called the {\it center}  of $\nu$ (or of $E$) on $X$. It is denoted by 
$\cent_X(\nu)$ or $\cent_X(E)$.

Thus we obtain
a function 
$$
 a(\ ,X,D):\{\text{divisors of $K(X)$ with nonempty center on $X$}\}\to \r.
$$
(If $X$ is proper over a field $k$ then every divisor of $K(X)$ over $k$ has a
nonempty center.) 
\enddemo

\demo{3.4 Definition}
In order to get a global measure of the singularities of the pair $(X,D)$ we 
define  
$$
\align
\discrep(X,D)&\deq\inf_E\{a(E,X,D)|\text{$E$ is exceptional with nonempty center
on $X$}\},\\
\tdiscrep(X,D)&\deq\inf_E\{a(E,X,D)|\text{$E$ has nonempty
center on $X$}\}.
\endalign
$$
\enddemo

\demo{3.4.1 Examples} (3.4.1.1) Let $E\subset X$ be a divisor different from
any of the
$D_i$. Then $a(E,X,D)=0$, thus $\tdiscrep(X,D)\leq 0$.

(3.4.1.2) Let $E$ be a divisor obtained by blowing up a
codimension 2 smooth point $x\in X$ which is not contained in  any
of the
$D_i$. Then $a(E,X,D)=1$, thus $\discrep(X,D)\leq 1$.

(3.4.1.3) If $X$ is smooth then  $K_Y=f^*K_X+E$ where $E$ is
effective and its support is the whole exceptional divisor.
Thus $\discrep(X,0)=1$. 

(3.4.1.4) (cf. \cite{CKM88, 6.3})  Show that
$$
\alignat4
&\qtq{ either}&\discrep(X,D)&=-\infty,&&\qtq{or} 
&&-1\leq \discrep(X,D)\leq 1,\qtq{and}\\
&\qtq{ either}
&\tdiscrep(X,D)&=-\infty,&&\qtq{or} &&-1\leq \tdiscrep(X,D)\leq 0.
\endalignat
$$

(3.4.1.5) Assume that the $D_i$ are $\r$-Cartier. Let $D=\sum d_iD_i$ and
$D'=\sum d'_iD_i$. If $d'_i\geq d_i$ for every $i$, then
$\discrep(X,D')\leq \discrep(X,D)$. \qed
\enddemo

Every restriction on $\discrep(X,D)$ defines a class of pairs $(X,D)$. The
following cases  emerged as the most important ones:

\demo{3.5 Definition} Let $X$ be a normal scheme and $D=\sum d_iD_i$  
a (not necessarily effective) $\r$-divisor such that
 $K_X+D$ is $\r$-Cartier. We say that $(X,D)$ or $K_X+D$ is
$$
\cases
\text{terminal}\\
\text{canonical}\\
\text{klt (or Kawamata log terminal)}\\
\text{plt (or purely log terminal)}\\
\text{lc (or log canonical)}\\
\endcases
\qtq{iff} \discrep(X) 
 \quad
\cases
>0,\\
\geq 0,\\
>-1\qtq{and $\rdown{D}\leq 0$,}\\
>-1,\\
\geq -1.\\
\endcases
$$
Equivalently, one can  define klt by the condition $\tdiscrep(X,D)>-1$. 
\enddemo

In order to get a feeling for these concepts, let us give some examples. In
dimension two these notions correspond to well-known classes of singularities.
The proof of the first two parts is an easy exercise using the minimal
resolution.  The last two cases are trickier.  See, for instance, 
\cite{Koll\'ar et al.92,3}. 

\proclaim{3.6 Theorem} Let $0\in X$ be a (germ of a)
normal surface singularity over $\c$. Then $X$
is
 $$
\align
\text{terminal}&\quad\Leftrightarrow\quad\text{ smooth;}\\
\text{canonical}&\quad\Leftrightarrow\quad\c^2/\text{(finite subgroup of 
$SL(2,\c)$);}\\
\text{klt}&\quad\Leftrightarrow\quad\c^2/\text{(finite subgroup of
$GL(2,\c)$);}\\ 
\text{lc}&\quad\Leftrightarrow\quad
\text{simple elliptic, cusp, smooth, or a}\\
&\hphantom{\quad\Leftrightarrow\quad}\text{ quotient of these
 by a finite group.}\qed\\ 
\endalign
$$
\endproclaim

Log canonical pairs appear naturaly in many contexts:

\proclaim{3.7 Proposition} Let $X$ be a normal toric variety with open orbit
$T\subset X$ and set $D=X-T$. Then $(X,D)$ is lc. If $K_X$ is $\q$-Cartier,
then $(X,0)$ is klt.
\endproclaim

\demop  Let $f: Y \to X$ be a toric resolution and set   $E=Y-f^{-1}(T)$.
By \cite{Fulton93,4.3}, $K_X\sim -D$ and $K_Y\sim -E$. Thus
$K_Y+E\sim f^*(K_X+D)$, and so $(X,D)$ is lc.  The rest is easy (cf.
(3.4.1.5)).
\qed\enddemo

As was observed by \cite{Alexeev96, Sec.3}, this easily implies that
 Baily--Borel
compactifications are also log canonical:

\proclaim{3.7.1 Corollary}  Let $D$ be a bounded symmetric domain
and
$\Gamma$ an arithmetic subgroup of $\aut(D)$. Let $(D/\Gamma)^*$ denote the
Baily--Borel compactification of $D/\Gamma$. There is a natural choice for a
$\q$-divisor $\Delta$, supported on the boundary, such that
$((D/\Gamma)^*,\Delta)$ is lc.\qed
\endproclaim

\demo{3.8 Example: Cones}  

Let $Y$ be a smooth variety and $E\subset Y$  a smooth divisor with normal
bundle $L^{-1}$. If $L$ is ample then $E$ is contractible to a point, at
least as an analytic space. Let
$f:Y\to X$ be this contraction. 

(3.8.1) $K_X$ is $\q$-Cartier iff $K_E$ and $L|E$ are linearly dependent in
$\pic E$. Assume that this is the case and write 
$$
K_E\equiv -(1+a)(L|E), \qtq{thus} K_Y\equiv f^*K_X+aE.
$$

(3.8.2) If $K_E$ is ample then $a<-1$; $X$ is not log canonical.

(3.8.3) If $K_E=0$ then $a=-1$ and $X$ is log canonical.

(3.8.4) If $-K_E$ is ample, that is $E$ is Fano, we have 3 cases.

(3.8.4.1) $0<L<-K_E$. Then $a>0$ and $X$ is terminal. Notice that there are
very few Fano varieties for which this can happen, so this is  a very rare case.

(3.8.4.2) $L=-K_E$. Then $a=0$ and $X$ is canonical. For every Fano variety
we get one example, still only a few cases.

(3.8.4.3) $L>-K_E$. Then $-1<a<0$ and $X$ is klt. For every Fano 
variety we get infinitely many cases.
\enddemo

\demo{3.9 Example}  Let $X$ be a smooth variety and $D\subset X$ a
 divisor. The discrepancies are all integers. 

(3.9.1) Show that $(X,D)$ is terminal iff $D=\emptyset$. 

(3.9.2) We see in (7.9) that $(X,D)$ is canonical iff $D$ is reduced, normal
and has rational singularities only.

(3.9.3) The case when $(X,D)$ is log canonical does not seem to have a
traditional name. Show that if $\dim X=2$, then  $(X,D)$ is log canonical iff $D$
has normal crossings only. In dimension three the list is much longer. $D$ can
have pinch points, rational double points, simple elliptic points (like
$(x^2+y^3+z^6=0)$) and cusps (like $(xyz+z^p+y^q+z^r=0)$). 

(3.9.4) Let $X$ be a smooth variety of dimension $n$ and $B$ a $\q$-divisor.
Assume that $(X,B)$ is terminal (resp. canonical, plt, lc). Study the blow
up of  $x\in X$ to show that
$\mult_x(B)<n-1$ (resp. $\leq n-1$, $<n$, $\leq n$).
(The converse is not true, see (3.14).)
\enddemo

The definition (3.5) requires some understanding of all exceptional divisors of
all
  birational modifications of $X$. The following  lemmas reduce the
computation of $\discrep(X,D)$ to a finite computation  in principle.

\proclaim{3.10 Lemma}  Notation as above.
 Let $f:X'\to X$ be a proper,
birational morphism  and write
$K_{X'}+D'\equiv f^*(K_X+D)$ (using   (3.3.2)). 

(3.10.1) 
$a(E,X,D)=a(E,X',D')$ for every divisor $E$ of $K(X)$.

(3.10.2)  $(X,D)$ is klt (resp. lc) iff  $(X',D')$ is klt (resp. lc).

(3.10.3)  $(X,D)$ is plt iff  $(X',D')$ is plt and $a(E,X,D)>-1$  for  every
exceptional divisor $E\subset X'$ of $f$. 

(3.10.4)  $(X,D)$ is terminal (resp. canonical) iff  $(X',D')$ is terminal
(resp. canonical) and $a(E,X,D)>0$ (resp. $a(E,X,D)\geq 0$) for 
every exceptional divisor $E\subset X'$ of $f$.\qed
\endproclaim

\proclaim{3.11 Lemma}  Let $X$ be a smooth scheme and 
$\sum D_i$ a normal crossing divisor.  Set $D=\sum d_iD_i$ and assume that
$d_i\leq 1$ for every $i$.  Let
$x\in X$ be a (not necessarily closed) point and
 $E$  a divisor of $K(X)$ such
that $\cent_X(E)=x$. Then

(3.11.1) 
$a(E,X,D)\geq \codim(x,X)-1-\sum_{j:x\in D_j}d_j$.

(3.11.2) $\tdiscrep(X,D)=\min\{0, -d_i\}$, and

(3.11.3) $\discrep(X,D)=\min\{1,1-d_i, 1-d_i-d_j\: D_i\cap D_j\neq
\emptyset\}$.
\endproclaim

\demop 
Pick a birational morphism $f:Y\to X$ such that $ E\subset Y$  is an exceptional
divisor with   general point
$y\in E$. By localizing at $x=f(y)$ we may assume that $x$ is a closed point. 
 Pick a local coordinate system $\{y_j\}$ such that $E=(y_1=0)$.
Possibly after reindexing, let $D_1,\dots, D_k$ be those divisors which pass
through
$f(y)$.
Let $x_i$ be a local coordinate system  at $x$ such that
$D_i=(x_i=0) $ for $i=1, \dots, k$.  
Set $c_i=d_i$ for $i\leq k$ and $c_i=0$ for $i>k$. 
We can write $f^*x_i=y_1^{a_i}u_i$ where $a_i>0$ for $ i\leq k$ and  
$u_i$ is a unit at $y$. Then
$$
f^*\frac{dx_i}{x_i^{c_i}}=a_iy_1^{(1-c_i)a_i-1}u_i^{1-c_i}dy_1+
y_1^{(1-c_i)a_i}\omega_i,\qtq{where
$\omega_i$ is regular at $y$.}
$$
Therefore in 
$$
f^*\frac{dx_1\wedge\dots\wedge dx_n}{x_1^{c_1}\cdots x_n^{c_n}}
$$
the only terms which could have a pole at $y$ are of the form
$$
\align
&y_1^{A_i}dy_1\wedge \omega_1\wedge \dots \wedge \widehat{\omega_i}\wedge \dots
\wedge\omega_n\qtq{where} \\
&A_i=-1+\sum_{j=1}^n (1-c_j)a_j= -1+\sum_{j=1}^k (1-d_j)a_j
\geq k-1-\sum_{j=1}^k d_j.
\endalign
$$
The rest is a simple computation.\qed\enddemo

\proclaim{3.12 Corollary}  Let $X$ be a smooth scheme and 
$\sum D_i$ a normal crossing divisor. Then $(X,\sum d_iD_i)$ is 
$$
\cases
\text{terminal}\\
\text{canonical}\\
\text{klt}\\
\text{plt}\\
\text{lc}\\
\endcases
\qquad
\text{   iff  }
 \quad
\cases
d_i<1 \qtq{and $d_i+d_j<1$ if $D_i\cap D_j\neq\emptyset,$}\\
d_i\leq 1 \qtq{and $d_i+d_j\leq 1$ if $D_i\cap D_j\neq\emptyset,$}\\
d_i<1,\\
d_i\leq 1 \qtq{and $d_i+d_j<2$ if $D_i\cap D_j\neq\emptyset,$}\\
d_i\leq 1.\hfill  \qed\\
\endcases
$$
\endproclaim

\proclaim{3.13 Corollary}  Let $(X,D=\sum d_iD_i)$ be a pair and $f:Y\to X$ a
log resolution of singularities.  Let $E_j\subset Y$ be the exceptional divisors
of
$f$. Then $(X,D)$ is lc iff $\min\{a(E_j,X,D), -d_i\}\geq -1$. 

If $(X,D)$ is lc then the following hold:

(3.13.1) $\tdiscrep(X,D)=\min\{a(E_j,X,D), -d_i\}$. 

(3.13.2) $\discrep(X,\emptyset)=\min\{ a(E_j,X,D), 1\}$.

(3.13.3)  If $f^{-1}_*(\supp D)$ is smooth, then 
$\discrep(X,D)=\min\{ a(E_j,X,D), 1-d_i\}$. \qed
\endproclaim

In the canonical case, the conditions (3.12) imposed by codimension
one and two points of $X$ impose stronger restrictions than those imposed by
higher codimension points.
One might expect that this is also the case for arbitrary divisors. 
The following exercise shows that this is not the case:

\demo{3.14 Exercise} 
Let $X$ be a smooth variety and $D$ an effective $\q$-divisor. Show that

(3.14.1) If $\mult_xD\leq 1$ for every $x\in X$ then $(X,D)$  is canonical.

(3.14.2) If $\dim X=2$ then the converse also holds.

(3.14.3) Give an example of a pair $(\c^2,D)$ such that $\mult_xD\leq 1$ outside
the origin, $\mult_0D=1+\epsilon$ and $(\c^2,D)$ is not even log canonical.

(3.14.4) If $\dim X\geq 3$ then the converse of (3.14.1) does not hold.

(3.14.5) For every $\epsilon>0$ give an example of a pair $(\c^n,D)$ such that
$\mult_xD\leq 1$ for $x$ outside the origin, $\mult_0D=1+\epsilon$ and $(\c^n,D)$
is not even log canonical.
\enddemo

Unfortunately, (3.4.1.3)  does not characterize smooth points,
except in dimension 2. The  problem is that the discrepancy 1 is caused by
codimenson 2 effects, and it gives very little information about $X$ in higher
codimension.  The following question corrects this:

\proclaim{3.15 Conjecture} \cite{Shokurov88} Let $X$ be a  normal scheme, 
$D$ an effective $\r$-divisor and $x\in X$ a closed point. Assume that $K_X+D$ is
$\r$-Cartier. Then 
$$
\inf_E\{a(E,X,D)|\cent_X(E)=x\}\leq \dim X-1,
$$
and equality holds iff $X$ is smooth and $0\notin \supp D$.
(This is proved for $\dim X\leq 3$ by  \cite{Markushevich96}.)
\endproclaim

It is also possible to compare the discrepancies for pull-backs by finite
morphisms, though the relationship is not as close as in (3.10).

\proclaim{3.16 Proposition} Let $p:X\to Y$ be a  finite and
dominant morphism between normal varieties. Let $D_Y$ be a $\q$-divisor on $Y$
and define $D_X$ by the formula
$$
K_X+D_X=p^*(K_Y+D_Y),\qtq{that is,} D_X=p^*D_Y-K_{X/Y}.
$$
Then 

(3.16.1) $1+\discrep(Y,D_Y)\leq 1+\discrep(X,D_X)\leq
\deg(X/Y)(1+\discrep(Y,D_Y))$.

(3.16.2) $(X,D_X)$ is lc (resp. klt) iff $(Y,D_Y)$ is lc (resp. klt).
\endproclaim 

\demop  (3.16.2) is a special case of (3.16.1). 

In order to prove (3.16.1), 
let $f_Y:Y'\to Y$ be a proper birational morphism and
$X'\to Y'\times_YX$ the normalization of the (dominant component of the) fiber
product with projection maps $f_X:X'\to X$ and $q:X'\to Y'$. Write
$$
K_{X'}+D_{X'}=f_X^*(K_X+D_X)\qtq{and} K_{Y'}+D_{Y'}=f_Y^*(K_Y+D_Y).
$$
$D_{X'}$ and $D_{Y'}$ are related by the formula
$$
K_{X'}+D_{X'}=q^*(K_{Y'}+D_{Y'}).
$$
 In order to compare the
coefficients in $ D_{X'}$ and in $D_{Y'}$, we may localize at the generic
point of a component  $E_Y\subset \supp D_{Y'}$.  
Let $E_X\subset \supp D_{X'}$ be a component which dominates $E_Y$.
Thus we are reduced to the
case when
$y\in Y'$ and $x\in X'$ are smooth pointed curves, $D'_Y=d_y[y],\ D'_X=d_x[x]$
and
$q: X'\to Y'$ has ramification index $r$ at $x$. 
Here $d_x=-a(X,D_X,E_X)$, $d_y=-a(Y,D_Y,E_Y)$  and
$r\leq
\deg(X/Y)$. Then
$d_x=rd_y-(r-1)$, or equivalently,
$$
a(X,D_X,E_X)+1=r(a(Y,D_Y,E_Y)+1).
$$
As $Y'\to Y$  runs through all proper birational morphisms, the
corresponding morphisms $X'\to X$ do  not give all possible proper birational
morphisms. Nonetheless, every algebraic valuation of $K(X)$ with center on $X$
appears on some $X'\to X$ by (3.17). This shows (3.16.1).\qed\enddemo

The proof used the next result which is essentially due to \cite{Zariski39}.
See also \cite{Artin86, 5.1; Koll\'ar95b,VI.1.3}. In many instances it can be
used instead of the resolution of singularities.

\proclaim{3.17 Theorem}   Let $X,Y$ be  integral
schemes of finite type  over a field or over $\z$, and   $f:Y\to X$  a 
dominant  morphism. Let
$D\subset Y$ be an irreducible  divisor  and  $y\in D$  the generic point. 
Assume that $Y$ is normal at $y$.
 We define a
sequence of schemes and maps as follows:

$X_0=X, f_0=f$. 

If  $f_i:Y\map X_i$ is already defined, then let
$Z_i\subset X_i$ be the closure of $f_i(y)$.  Let $X_{i+1}=B_{Z_i}X_i$ and
$f_{i+1}:Y\map X_{i+1}$   the induced map.

Then   $\dim Z_n\geq \dim X-1$ 
and $X_n$ is regular at the generic point of $Z_n$ for some $n\geq 0$.\qed
\endproclaim

The notion of a klt pair $(X,\Delta)$ also emerged  naturally in the
theory of singular metrics on line bundles. I just give the basic definition
and prove its equivalence with the algebro-geometric one. For a more detailed
exposition of the theory, see, for instance, \cite{Demailly92,94}.

\demo{3.18 Singular Metrics on Line Bundles} Let $L$ be a line bundle on a
complex manifold $M$. A {\it singular Hermitian metric} 
$\norm{\ }$ on $L$ is a 
Hermitian metric on 
$L|M-Z$  (where $Z$ is a measure zero set) such that if $U\subset M$
is any open subset, $u:L|U\cong U\times
\c$  a local trivialization  and $f$   a local generating section over $U$
then
$$
\norm{f}=|u( f)|\cdot e^{-\phi}
$$
where $|\ |$ is the usual absolute value on $\c$ and 
$\phi\in L^1_{loc}(U)$. (The latter assumption  assures that
$\partial\bar\partial \phi$ exists as a current on $M$. We do not use it.)

We say that the metric is $L^p$ on $M$ if $e^{-\phi}$ is locally $L^p$ for
every point. (This is clearly independent of the local trivializations.)
\enddemo

\demo{3.19 Examples} (3.19.1) Let $D$ be a  divisor and  $L=\o_X(D)$.
$L$ has a natural section $f$ coming from the constant  section  1 of $\o_X$.
A natural choice of the metric on $L$ is to set
$\norm{f}=1$ everywhere. This metric is singular along $D$. If $h$ is  a local
equation of $D$ at a point $x\in D$  
then $h^{-1}f$ is
a local generating section of $L$ at $x$ and
$$
\norm{h^{-1}f}=e^{-\log |h|}.
$$
 
(3.19.2) Let $L$ be a line bundle on $X$. Assume that $L^n\cong M(D)$ for some
line bundle $M$ and effective divisor $D$.
Let $\norm{\ }_M$ be  a continuous Hermitian metric on $M$. We construct  a
singular metric 
$\norm{\ }_L$ on $L$ as follows. 

Let $f$ be a local section of $L$ at a point $x\in D$ and  $h$  a local
equation for $D$ at $x$.   $hf^n$ is a local section of 
$M$, thus we can set
$$
\norm{f}_L\deq (\norm{hf^n}_M)^{1/n}e^{-log|h|/n}.
$$
The first factor on the right is continuous
and positive if $f$ is a local generator. Thus $\norm{\ }_L$ is $L^p$
iff the exponential factor
$$
e^{-log|h|/n}=|h|^{-1/n}\qtq{is $L^p$.}
$$

(3.19.3) Assume that  $D=n\sum d_i\Delta_i$
where $\Delta_i=(x_i=0)$ for a local coordinate system. Then
$$
\text{$\norm{\ }_L$ is $L^2$ at $x$} \quad\Leftrightarrow\quad
\text{$\prod |x_i|^{-d_i}$ is $L^2$ near $x$} \quad\Leftrightarrow\quad
\text{$d_i<1$ for every $i$.}
$$
\enddemo

More generally 
we have:

\proclaim{3.20 Proposition}  Let $X$ be a smooth manifold
and  $D$  a divisor on $X$. Let $L$ be a line bundle on $X$ and assume
that $L^n=M(D)$ for some line bundle $M$.
Set $\Delta=D/n$. 
 Let $\norm{\ }_L$ be the singular
metric constructed on $L$ as in (3.19.2). Then
$$
\norm{\ }_L \qtq{is $L^2$}\Leftrightarrow\quad (X,\Delta) \qtq{ is klt.}
$$
\endproclaim

\demop Let $f:Y\to X$ be a log resolution of $(X,D)$ and $E\subset Y$ the
exceptional divisor.  Both properties are local in $X$, so pick a point $x\in
D$ and fix a local coordinate system $\{x_i\}$. 
Let $h$ be a local equation for $D$.
Set $\omega_x=dx_1\wedge \dots\wedge dx_k$.   $\norm{\ }_L$ is $L^2$ iff
$$
\int |h|^{-2/n}\omega_x\wedge \bar \omega_x<\infty.\tag 3.20.1
$$
This is equivalent to saying that $|h|^{-1/n}\omega_x$ is $L^2$.
The advantage of putting $\omega_x$ in is that in this form the condition is
invariant under pull backs. Thus (3.20.1) is equivalent to
$f^*(|h|^{-1/n}\omega_x)$ being $L^2$ on $Y$.

This is a local condition on $Y$, so pick a point and  a local coordinate system
$\{y_i\}$ such that every irreducible component $F_i\subset E+f^{-1}_*(D)$ is
defined
 locally  as
$F_i=(y_i=0)$.
 Set $\omega_y=dy_1\wedge \dots\wedge dy_k$.  We can write
$$
f^*(|h|^{-1/n}\omega_x)=\omega_y\cdot u\cdot \prod |y_i|^{a_i},
$$
where $u$ is  invertible and 
$a_i=a(F_i,X,\Delta)$. Thus 
$f^*(|h|^{-1/n}\omega_x)$ is 
$L^2$ iff $a(F_i,X,\Delta)>-1$ for every $i$.
This happens precisely when $(X,\Delta)$ 
is klt.
\qed\enddemo

\head 4. Bertini Theorems
\endhead

The classical Bertini theorem says that  on a smooth variety a general
member of a base point free linear system is again smooth.  Actually,  the
original version of the Bertini theorem applies in the case of linear systems
with base points and it says the following:

\proclaim{4.1 Theorem} (Bertini) Let $X$ be a smooth variety over a field of
characteristic zero and $|B_1,\dots,B_k|$ the linear system spanned by the
(effective) divisors $B_i$. Let $B\in |B_1,\dots,B_k|$ be a general member.
Then  
$$
\mult_xB\leq 1+\inf_i\{\mult_xB_i\}\qtq{for every $x\in X$.}\qed
$$
\endproclaim

If $x\not\in B_i$ then $\mult_xB_i=0$, and so this theorem includes the usual
form as a special case. 

Recently
\cite{Xu94} studied the case of linear systems with base points.
\cite{Xu94} proved a variant of (4.1)  which also  implies similar results about
infinitely near singularities of $B$, though it is not clear to me whether his
results can be interpreted in terms of multiplicities alone. 
The results of this section grew out of trying to understand his results in
terms of discrepancies.

In order to get a better idea of what is possible, let us consider the simplest
case  of linear systems with base points: at each point there is a smooth
member. The general member will not, in general, be smooth:

\demo{4.2 Examples} (4.2.1) Let $X=\c^n$ and $f\in \c[x_3,\dots,x_n]$
such that    $(f=0)$ has an isolated singularity at
the origin. Consider the linear system $|B|=(\lambda x_1+\mu x_1x_2+
\nu f=0)$. At each point there is a smooth
member and the general member is
singular at 
$(0,-\lambda/\mu,0,\dots,0)$. All general members are isomorphic to 
$(x_1x_2+f=0)$. 
This way we can get any isolated $cA$-type singularity (4.3) for suitable choice
of
 $f$.

(4.2.2) As above, let $n=3$ and $f=x_3^{m+1}$, then general members have an
$A_m$ singularity, which is canonical but not terminal.

(4.2.3)  Consider the linear system $\lambda(x^2+zy^2)+\mu y^2$. At any point
$x\in \c^3$ its general member has a $cA$-type singularity, but the general
member has a moving pinch point. 
\enddemo

\demo{4.3 Definition} Let $0\in H\subset
X$ (where $X$ is smooth at $0$) be a hypersurface singularity. 
In local coordinates $H=(g=0)$; let $g_2$ denote the quadratic part of $g$.
We say that $H$ has
type  $cA$ if  either $H$ is smooth or $g_2$  has rank at least 2 (as a
quadratic form).
Equivalently, there are
 suitable local analytic (or formal) coordinates 
$x_1,\dots,x_n$ such that $H=(x_1x_2+f(x_3,\dots,x_n)=0)$
\cite{AGV85,I.11.1}.

 Let $0\in Y\subset X$ be a smooth hypersurface. If
$H\cap Y$ has a $cA$-type singularity at $0$ then so does $H$.
By following the Hironaka resolution process, it is easy to see that a normal
$cA$-type singularity is canonical. \enddemo

In the smooth case one can give a very precise description of the possible
singularities occurring on general members of linear systems.

\proclaim{4.4 Theorem}  Let $X$ be a smooth variety over a field of
characteristic zero and  $|B|$  
a linear system of Cartier divisors.
Assume that for every $x\in X$ there is a $B(x)\in |B|$ such that $B(x)$ is
 smooth at $x$ (or $x\not\in B(x)$). 

Then a   general  member 
$B^g\in |B|$ has only type $cA$ singularities.
\endproclaim

\demop The result is clear if $\dim X=1$, thus assume that it holds for
smaller dimensional schemes. By Noetherian induction it is sufficient to prove
the following:

(4.4.1) For every irreducible subvariety $Z\subset X$ there is an open subset
$Z^0\subset Z$ such that a   general  member 
$B^g\in |B|$ has only type $cA$ singularities at points of $Z^0$.

If $Z\not\subset \bs|B|$ then let $Z^0=Z-\bs|B|$. The base point free Bertini
then implies (4.4.1).

Next assume that $Z\subset\bs|B|$ and $\codim(Z,X)=1$. The assumptions imply
that $Z$ is smooth  and $|B|-Z$ induces a base point free linear system on $Z$.
Again by Bertini, the general member $B^g-Z$ intersects $Z$ transversally. Thus
at every point of
$Z$ the divisor $B^g$ is either smooth or its local equation is
$x_1x_2=0$.

Finally assume that $Z\subset\bs|B|$ and $\codim(Z,X)>1$. Pick a smooth point 
$z\in Z$. Let $Y$ be a hypersurface in $X$ such that:

(4.4.2.1) $Z\subset Y$ and $Y$ is smooth at $z$, and

(4.4.2.2) $Y$ is transversal to $B(z)$ at $z$.

Let $Y^0\subset Y$ be an open set containing $z$ such that $Y^0$  and $B(z)\cap
Y^0$ are smooth.  Let $|B_{Y^0}|$ be the restriction of  the linear system $|B|$
to $Y^0$.  By induction $B^g\cap Y^0$ has only type $cA$ singularities, hence by
(4.3) $B^g$ has 
 only type $cA$ singularities at points of $Z\cap Y^0$.\qed\enddemo

The above results  say that if $|B|$ is a linear system and for each $x\in
X$ there is a $B(x)\in |B|$ which is not too singular at $x$, then the general
$B^g\in |B|$ has only somewhat worse singularities. This raises the question:
is there a class of singularities for which the general member does not get
any worse?

Assume that $\bold S$ is such a class. By (4.2.1)  $\bold S$ contains all
$cA$-type singularities. Thus by (4.2.1) it also has  to contain pinch
points and maybe many more singularities. It is not at all clear that this
process ever terminates with a reasonably small class $\bold S$.  I
do not know what is the smallest class $\bold S$ (it is clear that it
exists). 

 The
following result provides one such example for $\bold S$, under a mild
assumption on the linear series. More examples  are contained in (4.8).

\proclaim{4.5 Theorem}  Let $X$ be a smooth variety over a field of
characteristic zero and  $|L|$  
a linear system of Cartier divisors such that $\bs|L|$ has codimension at
least two.  Assume that for every $x\in X$ there is a $B_x\in |L|$ such that
$B_x$ has a rational singularity at $x$. 

Then a   general  member 
$B\in |L|$ has only rational singularities.
\endproclaim

The above result  has nothing to do with discrepancies or with
canonical pairs. Still, I have no idea how to prove it without   the
machinery of adjunction and  canonical pairs. After translating the problem
to our language, it becomes easy.

\demop  $B$ is a Cartier divisor on a smooth variety, hence $\omega_B$ is
locally free. By (11.1.1) $B$ has rational singularities iff it has canonical
singularities.
By (7.9)  the latter holds
iff the pair $(X,B)$ is canonical. Thus (4.5) is equivalent to the following
version:

\proclaim{4.5.1 Theorem}  Let $X$ be a smooth variety over a field of
characteristic zero and  $|L|$  
a linear system of Cartier divisors such that $\bs|L|$ has codimension at
least two. 
Assume that for every $x\in X$ there is a $B_x\in |L|$ such that $(X,B_x)$
is canonical at $x$. 

Let  $B\in |L|$ be a   general  member. Then $(X,B)$ is canonical.
\endproclaim

\demop 
Let $f:Y\to X$ be a proper birational morphism
such that $Y$ is smooth, $f^*|L|=F+|M|$ where $|M|$ is base point free
 and $F+(f\text{-exceptional divisors})$ has only normal
crossings.   Let $K_Y=f^*K_X+E$.  For a given divisor $B\in |L|$, let 
$B^Y\deq f^*(B)-F\in |M|$ denote the corresponding member. We can write
$B^Y=f^{-1}_*B+N$ where $N$ is effective and empty for general $B\in |L|$. 
(If $\bs|L|$ contains a divisor, then $N$ is not effective and the rest of
the proof does not work.) Then
$$
K_Y+f^{-1}_*B_x=K_Y+B^Y_x-N_x=f^*(K_X+B_x)+(E-F-N_x).
$$
By assumption $(X,B_x)$ is canonical at $x$, thus $E-F$ is effective over a
neighborhood of $x$.  This holds for any $x$, thus $E-F$ is effective.

Choose $B\in |L|$ such that the corresponding $B^Y$ is irreducible and
intersects $E-F$ transversally. These are both nonempty and open conditions.
Then
$$
K_Y+f^{-1}_*B=K_Y+B^Y=f^*(K_X+B)+(E-F).
$$
This shows that $(X,B)$ is canonical.\qed\enddemo
\enddemo
$$
*****
$$
In the rest of the section we study Bertini-type theorems
that compare the properties   klt, lc and so on for the general members and for
generators of linear systems. There are two ways of approaching this problem.
One can look at singularities of divisors $B$ or singularities of pairs $(X,B)$.
The second variant is better
suited for the present purposes. In some cases the two versions are
equivalent (4.9).

\demo{4.6 Definition}  Let $X$ be a normal, integral scheme, $D=\sum
d_iD_i$ a
$\q$-divisor (not necessarily effective) and $|B_j|$ (not necessarily
complete)  linear systems of Weil divisors.   Let $0\leq b_j\leq 1$ be
rational numbers such that 
$K_X+D+\sum b_jB_j$ is
$\q$-Cartier. Let
$E$ be a divisor of the function field $\c(X)$ and define
$$
a(E,X,D+\sum b_j|B_j|)\deq \sup\{a(E,X,D+\sum b_jB'_j)\vert B'_j\in |B_j|\}.
$$
In the above formula it is sufficient to let $B'_j$ run through a finite set of
divisors spanning $|B_j|$. In  particular the supremum is a maximum and
if the $b_j$ are rational then so is $a(E,X,D+\sum b_j|B_j|)$.
 We  define as in (3.4)
$$
\align
&\discrep(X,D+\sum b_j|B_j|)=\\
&\quad =\inf_E\{a(E,X,D+\sum b_j|B_j|)|\text{$E$ is
exceptional with nonempty center on
$X$}\}.
\endalign
$$

As in (3.5),  we say that $(X,D+\sum b_j|B_j|)$ or $K_X+D+\sum b_j|B_j|$ is
{\it terminal, canonical, klt, plt} resp. {\it lc} if 
$\discrep(X,D+\sum b_j|B_j|)>0, \geq 0, >-1$ and $d_i,b_j<1\ \forall j, >-1$
resp.
$\geq -1$. 
\enddemo

The following properties are straightforward from the definitions. 

\proclaim{4.7 Lemma}  Notation as above.

(4.7.1) If $|B_1|$ is base point free then
$$
a(E,X,D+\sum b_j|B_j|)=a(E,X,D+\sum_{j\geq 2} b_j|B_j|).
$$

(4.7.2) If $F_j\subset \bs |B_j|$ is a divisor, then
$$
a(E,X,D+\sum b_j|B_j|)=a(E,X,(D+\sum b_jF_j)+\sum b_j|B_j-F_j|).
$$

(4.7.3) Assume that the $B_j$ are $\q$-Cartier. Let $f:X'\to X$ be a proper,
birational morphism  and write
$K_{X'}+D'=f^*(K_X+D)$. Then
$$
a(E,X,D+\sum b_j|B_j|)=a(E,X',D'+\sum b_jf^*|B_j|).\qed
$$
\endproclaim

The following is a summary of the  Bertini-type  theorems for linear
systems. It should be clear from the proof that there are several other variants
involving different values of the discrepancy. Also, one can be more precise
concerning the interplay of the allowable coefficients $d_i$ and $b_j$.

\proclaim{4.8 Theorem}  Let $X$ be a normal, integral, excellent scheme over a
field of characteristic zero, 
$D=\sum d_iD_i$ a
$\q$-divisor (not necessarily effective) and $\sum b_j|B_j|$ a formal sum of
(not necessarily complete)  linear systems of Weil divisors, $0\leq b_j\leq 1$. 
Assume that
$K_X+D$ and the $B_j$ are
$\q$-Cartier. Let $B_j^g$ be a general member of $|B_j|$. Then:

(4.8.1) $(X,D+\sum b_j|B_j|)$ is lc $\Leftrightarrow$  $(X,D+\sum b_jB_j^g)$
 is lc.

(4.8.2)  $(X,D+\sum b_j|B_j|)$ is klt $\Leftrightarrow$ 
$(X,D+\sum b_jB_j^g)$ is klt, for $0\leq b_j<1$. 

\noindent Assume   that   general members of $|B_1|$ are irreducible
and $(X,D)$ is klt. Then: 

(4.8.2')  $(X,D+ b_1|B_1|)$ is
plt  $\Leftrightarrow$ 
$(X,D+b_1B_1^g)$ is plt. 

\noindent Assume   that   $d_i\leq 1/2$ for every $i$. Then:

(4.8.3)  $(X,D+\sum b_j|B_j|)$ is
canonical  $\Leftrightarrow$ 
$(X,D+\sum b_jB_j^g)$ is canonical, for $0\leq b_j\leq 1/2$.

(4.8.4)  $(X,D+\sum b_j|B_j|)$ is
terminal  $\Leftrightarrow$ 
$(X,D+\sum b_jB_j^g)$ is terminal, for $0\leq b_j< 1/2$. 

\noindent Assume   that   general members of $|B_1|$ are irreducible
and $D=\emptyset$. Then: 

(4.8.3')  $(X, b_1|B_1|)$ is
canonical  $\Leftrightarrow$ 
$(X,b_1B_1^g)$ is canonical. 

(4.8.4')  $(X, b_1|B_1|)$ is
terminal  $\Leftrightarrow$ 
$(X, b_1B_1^g)$ is terminal, for $0\leq b_1< 1$. 
\endproclaim

\demop Let $f:X'\to X$ be a proper birational morphism
such that $X'$ is smooth, $f^*|B_j|=F_j+|M_j|$ where $|M_j|$ is base point free
for every $j$ and $f^*D+\sum F_j+(f\text{-exceptional divisors})$ has only normal
crossings.  By (4.7.1--3) we see that 
$(X,D+\sum b_j|B_j|)$ is lc (resp. klt) iff 
$(X',D'+\sum b_jF_j)$ is lc (resp. klt).  Let $D'+\sum b_jF_j=\sum e_kE_k$.
By (3.12), $(X',\sum e_kE_k)$ is lc (resp. klt) iff $e_k\leq 1$ 
(resp. $e_k< 1$) for every $k$. 

Let $M_j^g\deq f^*B_j^g-F_j\in |M_j|$ be a general member. Then
$$
K_{X'}+\sum e_kE_k+\sum b_jM_j^g=f^*(K_X+D+\sum b_jB_j^g),
$$
and $\sum E_k+\sum M_j^g$ is still a normal crossing divisor by the usual
Bertini theorem.  Therefore 
$$
\align
&K_X+D+\sum b_jB_j^g\qtq{is lc (resp. klt) $\Leftrightarrow$ }\\
&K_{X'}+\sum e_kE_k+\sum b_jM_j^g\qtq{is lc (resp. klt), $\Leftrightarrow$ }\\
&e_k,b_j\leq 1\text{ (resp. $e_k,b_j< 1$) for every $j,k$.}
\endalign
$$
This shows
(4.8.1--2).

In order to obtain (4.8.2') we need the additional remark that if 
$\sum E_k+ M_1^g$ is  a normal crossing divisor, $M_1^g$ is irreducible and
$e_k<1$ for every
$k$ then
$(X',\sum e_kE_k+b_1M_1^g)$ is plt for every $b_1\leq 1$.

The proofs of the remaining assertions are   similar. 
$(X,D+\sum b_j|B_j|)$ is canonical  (resp. terminal) iff
$e_k\leq 0$ 
(resp. $e_k< 0$) for every $k$ such that $E_k$ is $f$-exceptional. Therefore 
$(X',\sum e_kE_k+\sum b_jM_j^g)$ is canonical  (resp. terminal)
iff $(X',f^{-1}_*D+\sum b_jM_j^g)$ is canonical  (resp. terminal). 

In (4.8.3--4) the
coefficient $1/2$ comes in because two divisors $M_i^g$ and $M_j^g$ may
intersect. In   (4.8.3'--4') we know that $M_1^g$ is irreducible, thus we
only have to make sure that its coefficient is at most 1 (resp. less than 1). 
\qed\enddemo

It is frequently more convenient to have Bertini-type theorems which
give information about the singularities of the general member of  a linear
system directly.  This is
rather straightforward  for base point free linear systems (7.7). 
For linear systems with base points the situation is more difficult to analyze.
The known cases of the inversion of adjunction conjecture (7.3) can be used to
transform the
$b_j=1$ cases of (4.8)  to Bertini-type results that concern only the
singularities of divisors. 
I formulate it only for Cartier
divisors. Using the notion of different (cf. \cite{Koll\'ar et al.92,16.6})  it
can be stated in case
$(X,\Delta)$ is klt and $|B|$ is a linear system of $\q$-Cartier Weil divisors.

\proclaim{4.9 Corollary} Let $X$ be a scheme over a field of characteristic
zero and $|B|$ 
a linear system of Cartier divisors such that the general member
 of $|B|$ is irreducible.

(4.9.1) 
Assume that   $X$ is  klt and for every $x\in X$ there is a $B(x)\in |B|$
such that $B(x)$ is
 klt  at $x$ (or $x\not\in B(x)$). 

Then $B^g$ is klt for   general 
$B^g\in |B|$.

(4.9.2) 
Assume that   $X$ is  canonical and of index 1 (for instance smooth).  Assume
also that for every
$x\in X$ there is a
$B(x)\in |B|$ such that $B(x)$ is
 canonical  at $x$ (or $x\not\in B(x)$). 

Then $B^g$ is canonical for   general 
$B^g\in |B|$.
\endproclaim

\demop  Let $B'\in |B|$ be any divisor. By (7.5)
we see that $K_X+B'$ is  plt at a point $x\in X$ iff
$B'$ is  klt at  $x$. Thus (4.9.1) follows from 
(4.8.2').

If $X$ has index 1 and $|B|$ is
a linear system of Cartier divisors, then every member of $|B|$ has index 1.
Hence klt is the same as canonical. Thus (4.9.1) implies (4.9.2).
\qed\enddemo

\demo{4.9.3 Remark}   It is expected that (4.9.2) remains true
even if $X$ has higher index canonical singularities.
\enddemo

\head 5. Effective Base Point Freeness
\endhead

In its simplest form the problem is the following:

\demo{5.1 Problem}  Let $X$ be a   projective variety and $L$ an ample
line bundle on $X$. Try to construct a very ample (or maybe just globally
generated) line bundle, using as little information about $X$ and $L$ as
possible.
\enddemo

The first major result of this type is ``Matsusaka's big theorem"
 which asserts the following:

\proclaim{5.2 Theorem} \cite{Matsusaka72}
There is a function
$\Phi(x,y,z)$ with the following property:

 If $X$ is an $n$-dimensional 
smooth projective variety over a field of
characteristic zero and
$L$   an ample line bundle on $X$, then 
$$
L^m\qtq{is very ample for} m\geq 
\Phi(n,(L^n),(K_X\cdot L^{n-1})). 
$$
\endproclaim

\cite{Matsusaka86}
generalized this to the case when $X$ has   at worst rational singularities.
The methods of  \cite{Matsusaka72,86} do not give
any information about $\phi$, beyond its existence. 

Mukai pointed out that a  reasonable bound can be expected if one tries to
find a very ample line bundle of the form $K_X\otimes L^m$. 
 The precise conjecture was formulated by \cite{Fujita87}:

\proclaim{5.3 Conjecture}  Let $Y$ be a  smooth projective variety,
  and $L$ an ample line bundle on $Y$.
Then:

(5.3.1)  $K_Y\otimes L^m$ is globally generated for $m\geq \dim Y+1$, and

(5.3.2)  $K_Y\otimes L^m$ is very ample for $m\geq \dim Y+2$.
\endproclaim

Both of the bounds are sharp for $Y=\p^n$ and $L=\o(1)$. The following example
gives many more  such cases:

\demo{5.3.3 Example} (Lazarsfeld)  Let $X\subset \p^{n+1}$ be a hypersurface of
degree $d$ and $L\subset \p^{n+1}$   a line intersecting $X$ in distinct
points
$P_1,\dots,P_d$. Let $p:Y\to X$ be the blow up of $P_1,\dots,P_{d-1}$
with exceptional divisors $E_i$. Set $L=p^*\o_X(1)(-\sum_{i=1}^{d-1}E_i)$. 

Show that $L$ is nef and big, and in fact it is generated by global sections
outside $P_d$.  $L$ is not always ample (for instance, if $X$ contains a line
through $P_1$) but $L$ is ample for general $X$ and $L$ for  $d$ sufficiently
large. 

$K_Y\otimes L^n\cong p^*\o_X(d-2)(-\sum_{i=1}^{d-1}E_i)$, and this has $P_d$ as its
base point. 

Another series of examples is in \cite{Kawachi96}. 
\enddemo

\demo{5.3.4 Remarks} This conjecture is true in low dimensions. The case $\dim Y=1$  
is very easy. The surface case follows from  \cite{Reider88}. (5.3.1) is quite
hard in dimension three 
\cite{Ein-Lazarsfeld93}. 
A very readable introduction is provided by the lectures \cite{Lazarsfeld96}.
The first step in all dimensions  was proved by \cite{Demailly93} who showed
that under the above assumptions $K_Y^2\otimes L^m$ is very ample for $m\geq
12n^n$ where
$n=\dim Y$.

 These results seem to furnish rather strong evidence in
favor of (5.3), but looking at the proofs gives a less optimistic picture.
The method of \cite{Ein-Lazarsfeld93}  gives a  base point freeness result
assuming that $(L^3)$ is large. In some vague sense only finitely many types of
cases remain to be analyzed. The study of these cases  
 requires considerable care and several ad hoc arguments. This is especially
the case for the proof of variants of (5.4) in dimension 3 given by 
\cite{Ein-Lazarsfeld93} and improved by 
\cite{Fujita94}.

Recently \cite{Kawamata96} proved (5.3) in dimension 4, and \cite{Smith95}
showed that (5.3) holds in positive characteristic for $L$ very ample. 
\cite{Helmke96} considerably improved on the earlier methods. With his
approach the low dimensional cases are now quite satisfactory, but for large
dimensions (5.3) is still out of reach.
\enddemo

The above conjecture was  given a more local form in \cite{Ein-Lazarsfeld93}.
As before, the actual bounds are inspired by the worst known example
$(\p^n,\o(n))$. 

\proclaim{5.4 Conjecture}  Let $Y$ be a  smooth projective variety,
$y\in Y$ a closed point 
  and $L$ a nef and big line bundle on $Y$.
Assume that if $y\in Z\subset Y$ is an irreducible (positive
dimensional) subvariety then
$$
 (L^{\dim Z}\cdot Z)\geq 
(\dim Y)^{\dim Z},\qtq{and} 
 (L^{\dim Y})>
(\dim Y)^{\dim Y}.
$$
Then $K_Y\otimes L$ has a section which is nonzero at $y$.
\endproclaim

The higher dimensional situation was recently greatly clarified by
\cite{Demailly94; Angehrn-Siu95; Tsuji95}.  Their results are weaker than
(5.4), but the proofs are very natural and the bounds  quite good.

\proclaim{5.5 Theorem}   Let $Y$ be a  smooth projective
variety  and $L$ an ample line bundle on $Y$.
Then

(5.5.1) $K_Y\otimes L^m$ is generated by global sections for $m>\binom{\dim
Y+1}{2}$.

(5.5.2) Global sections of $K_Y\otimes L^m$ separate points  for
$m\geq\binom{\dim Y+2}{2}$.
\endproclaim

\demo{5.6 Idea of the proof}
 Assume for simplicity that the linear system
$|L|$ is base point free. Let $D\in |L|$ be a general smooth member. For
$m>0$  we have an exact sequence
$$
H^0(X,K_X\otimes L^{m+1})\to H^0(D,K_D\otimes L^{m})\to
H^1(X,K_X\otimes L^m)=0.
$$
 Thus an induction on the dimension produces sections, giving even
the original  conjecture (5.3). 

In general  we cannot assume that $|L|$ is nonempty, let alone
that it has a smooth member. Assume instead that 
$$
L\equiv M+D+\Delta,
\tag 5.6.1
$$
 where $M$ is a nef and big $\q$-divisor, $D$
 a smooth (integral) divisor and $\Delta$ a $\q$-divisor such that
$\rdown{\Delta}=\emptyset$ and everything is in normal crossing. As before we
  get an exact sequence
$$
H^0(X,K_X\otimes L^{m+1})\to H^0(D,K_X\otimes L^{m+1}|D)\to
H^1(X,K_X\otimes L^{m+1}(-D)).
$$
Observe that 
$$
\align
K_X\otimes L^{m+1}(-D)&\equiv K_X+\Delta+mL+M,\qtq{and}\\
 K_X\otimes L^{m+1}|D&\equiv K_D+\Delta|D+(mL+M)|D.
\endalign
$$
Thus by (2.3) 
$H^1(X,K_X\otimes L^{m+1}(-D))=0$, and we can run induction on the dimension as
before, assuming that we can make the whole process work in the log version.

The assumption (5.6.1) seems strong, but it is easy to
achieve.  Pick a point $x\in X$ and assume that $(L^n)>1$. By (6.1) we can
find a
$\q$-divisor
$L\equiv B$ such that $c\deq \mult_xB> 1$.
 Let $\pi:X'\to X$ be the blow up of
$x$ and $E\subset Y$ the exceptional divisor. Then
$\pi^*B\equiv cE+B'$ and $c> 1$.  After further blowing ups the normal
crossing assumption can be satisfied, and  we obtain a proper birational
morphism  $p:Y\to X$ such that $p^*B=\sum c_iE_i$ and $\max\{c_i\}>1$. 
For suitable indexing the maximum is achieved for $c_0$. 
Assume for simplicity that $c_i<c_0$ for $i>0$. (Paragraph (6.3.5)  shows
what to do otherwise.) Then
$$
p^*L\equiv (1-c_0^{-1})p^*L+E_0+\sum_{i>0} \frac{c_i}{c_0}E_i,
$$
exactly as required  for (5.6.1).

The problem is that the pull back of $L$ is no longer ample, only nef. In the
worst case $p^*L|D\equiv 0$ and  $p^*(K_X\otimes L^{m+1})|D$ may have negative
degree.  The induction breaks down completely. This happens   already for
surfaces.

One way to get around this problem is to find 
$L\equiv B$ such that $c=\mult_xB\geq n$. At the level of the first blow
up $p:Y\to X$ we get that
$$
p^*(K_X\otimes L^{m+1})\equiv K_Y+B'+(c-(n-1))E+p^*L^m.\tag 5.6.2
$$
The advantage of this situation is that the divisor in (5.6.2) is a pull back,
so it has sections over the fibers of $p$. 

The inductive assumptions become rather messy and there are further technical
problems. Still, this idea can be made to work to get   some results, see
\cite{Koll\'ar95a, Ch.14}.

The idea of \cite{Angehrn-Siu95; Tsuji95} is to try to get  a section right
away. This is possible if we can  prove that
$H^1(X,K_X\otimes L^m\otimes m_x)=0$ where $m_x$ is the ideal sheaf of  $x\in
X$. This seems very hard to do. Fortunately, it is sufficient to  produce an
ideal sheaf
$J\subset\o_X$ such that 

(5.6.3.1) $x$ is an isolated point of $\spec (\o_X/J)$, and

(5.6.3.2) $H^1(X,K_X\otimes L^m\otimes J)=0$.

\noindent The  variant (2.16) of the vanishing theorem suggests such an
approach:

Try to find a $\q$-divisor $B$ such that

(5.6.4.1) $B\equiv (m-\epsilon)L$, and

(5.6.4.2) $B$ is not  klt at $x$ but is klt in a neighborhood of $x$.

The construction of such a divisor is not easy but turns out to be
feasible once a few technical points are settled. Thus the essential part of
the proof is postponed until section 6. \qed\enddemo

Properties as in (5.6.4.2) will appear frequently, so we introduce a notation
for it.

\demo{5.7 Definition}  Let $(X,D)$ be a pair. The set of points where $(X,D)$
is klt is open, it is called the {\it klt locus} of $(X,D)$.  The complement of
the klt locus is denoted by $\nklt(X,D)$; it is called the {\it non-klt
locus}.

Some  authors call this the ``locus of log canonical singularities". In
my view this may be  misleading. 
\enddemo

Here I give an algebraic version of the proof of \cite{Angehrn-Siu95}. I
state a  more general form which also applies to singular varieties.

\proclaim{5.8 Theorem} Let $(X,\Delta)$ be a proper klt pair and  $M$  a
line bundle. Assume that $M\equiv K_X+\Delta+N$, where $N$ is  a nef and big
$\q$-Cartier
$\q$-divisor on
$X$. Let
$x\in X$ be a closed point and assume that    there are positive numbers
$c(k)$ with the following properties:

(5.8.1) If
$x\in Z\subset X$  is an irreducible (positive dimensional) subvariety then 
$$
(N^{\dim Z}\cdot Z)>c(\dim Z)^{\dim Z}.
$$ 

(5.8.2) The numbers $c(k)$ satisfy the inequality
$$
\sum_{k=1}^{\dim X} \frac{k}{c(k)}\leq 1.
$$
Then $M$ has a global section not vanishing at $x$.
\endproclaim

Analogous results hold for separating points:

\proclaim{5.9 Theorem} Let $(X,\Delta)$ be a proper klt pair and  $M$  a
line bundle. Assume that $M\equiv K_X+\Delta+N$, where $N$ is  a nef and big
$\q$-Cartier
$\q$-divisor on
$X$. Let
$x_1,x_2\in X$ be  closed points and assume that    there are positive numbers
$c(k)$ with the following properties:

(5.9.1) If
$Z\subset X$  is an irreducible subvariety which contains $x_1$ or $x_2$ then 
$$
(N^{\dim Z}\cdot Z)>c(\dim Z)^{\dim Z}.
$$ 

(5.9.2) The numbers $c(k)$ satisfy the inequality
$$
\sum_{k=1}^{\dim X} \root{k}\of{2}\frac{k}{c(k)}\leq 1.
$$
\noindent Then global sections of $M$ separate $x_1$ and $x_2$.
\endproclaim

\demop Let us prove first (5.8). 

First I claim that the inequalities (5.8.1) are satisfied if we replace $N$ by
$(1-\epsilon)N$   for $0<\epsilon\ll 1$.  This is a minor technical step which
could easily have been built into the assumptions. A proof is given in
 (6.6.2).

Thus, by (6.4), there is a $\q$-divisor
$B\equiv (1-\epsilon)N$ such that
$x$ is an isolated point of $\nklt(X,\Delta+B)$. We can write
$$
M\equiv K_X+\Delta+B+\epsilon N.
$$
By (2.16) there is an
ideal sheaf
$J\subset\o_X$ such that $\supp (\o_X/J)=\nklt(X,\Delta+B)$ and
$H^i(X,M\otimes J)=0$ for $i>0$. In particular, the $i=1$ case
implies that
$$
H^0(X,M)\twoheadrightarrow H^0(X,M\otimes (\o_X/J))
\twoheadrightarrow H^0(X,M\otimes (\o_X/m_x))
\qtq{is surjective.}
$$
Thus $M$ has a section which does not vanish at $x$.  

The proof of (5.9) is similar.  We already know that $M$ has
sections which do not vanish at $x_1$ and at $x_2$. Thus global sections
separate $x_1$ and $x_2$ iff there is an $i\in \{1,2\}$  and a global section
$s\in H^0(X,M)$ such that $s(x_i)\neq 0$ and $s(x_{3-i})=0$. 
By (6.5) there is an $i\in \{1,2\}$ and a  $\q$-divisor $B\equiv
(1-\epsilon)N$ such that
$x_i$ is an isolated point of $\nklt(X,\Delta+B)$ and $x_{3-i}\in 
\nklt(X,\Delta+B)$. As before, this implies the existence of the required
section $s$.\qed\enddemo

\demo{5.10 Proof of (5.5)} Apply (5.8) and (5.9) with
$X=Y$, 
$\Delta=\emptyset$ and $N=L^m$.  Set $n=\dim Y$.  In the first case set
$c(k)=\binom{n+1}{2}$. This gives (5.5.1). 

 (5.9) implies (5.5.2)  by setting $c(k)=\binom{n+2}{2}$ and using the
inequality
$$
\sum_{k=1}^{n} \root{k}\of{2}k< \sum_{k=1}^{n}
\bigl(1+\frac1{k}\bigr)k=\binom{n+2}{2}-1.
\qed
$$
\enddemo
$$
*****
$$
The following is an   application of (5.9) to 
varieties with \glafg. See \cite{Koll\'ar95a} for the relevant definitions
and results.

\proclaim{5.11 Theorem} Let $X$ be a smooth proper variety with \glafg.
 Then:

(5.11.1) If $M$ is a big Cartier divisor on $X$ the $K_X+M$ is also big.

(5.11.2) $K_X$ is the limit of effective $\q$-divisors.
\endproclaim

\demop One can choose a suitable birational model $p:Y\to X$ such that 
$$
p^*M\equiv N+\Delta+ R
$$
where $N$ is an ample $\q$-divsior, $\Delta$  is a fractional normal crossing
divisor and
$R$ is effective. It is sufficient to prove that  
$$
K_Y+p^*M-R\equiv K_Y+N+\Delta
$$
is big.  If $N$ has sufficiently large intersection number with any
subvariety through a given point $y\in Y$, then by (5.9) global sections of
$K_Y+p^*M-R$  separate points, thus it is big.

 $X$ has \glafg, thus  there is a finite
\'etale cover
$q:Y'\to Y$ such that $q^*N$ has 
 large intersection number with any
subvariety through a given point $y'\in q^{-1}(y)$. (In fact these
two properties are equivalent.)  Thus
$$
K_{Y'}+q^*N+q^*\Delta\equiv  q^*(K_Y+N+\Delta)
$$
is big, hence so is $K_Y+N+\Delta$. 

By induction, $mK_Y+M$ is big for any $m\geq 1$, thus
$$
K_Y\equiv \lim_{m\to \infty} \frac1{m}(mK_Y+M)
$$
is  the limit of big  $\q$-divisors. A big $\q$-divisor is also effective,
proving (5.11.2).\qed\enddemo

\head 6. Construction of Singular Divisors
\endhead

 The aim of this section is to construct divisors
which are ``very singular" at a given point and not too singular elsewhere. The
precise measure of what we mean by ``very singular" is given by the
notion of discrepancy. Actually, the construction is even weaker in the
sense that we are   able to guarantee only that the resulting divisor is 
not too singular in a   neighborhood of our point.

The first step is to construct divisors which are as singular at a point as
possible. The optimal result is achieved by  an easy dimension count:

\proclaim{6.1 Lemma} Let $X$ be a proper and irreducible variety of
dimension $m$,
$M$ a nef and big $\q$-Cartier $\q$-divisor on $X$ and $x\in X$ a smooth point.
For every
$\epsilon>0$ there is an effective  $\q$-Cartier  $\q$-divisor divisor
$D=D(x,\epsilon)$ such that
$M\equiv D$ and
$$
\mult_xD\geq \root{m}\of{(M^m)}-\epsilon.
$$
\endproclaim

\demop Fix $s,t>0$ such that $tM$ is Cartier and let $m\subset \o_X$ be the
ideal sheaf of
$x$.  From the sequence
$$
 0\to m^s\otimes \o_X(tM)\to \o_X(tM)\to (\o_X/m^s)\otimes \o_X(tM)\cong
\o_X/m^s\to 0
$$ 
we see that 
$$
 h^0(X,m^s\otimes \o_X(tM))>0\qtq{if} H^0(X,\o_X(tM))>H^0(X,\o_X/m^s).
$$
 Since $x\in X$ is a smooth point, 
$$ H^0(X,\o_X/m^s)=\dim_kk[x_1,\dots,x_m]/(x_1,\dots,x_m)^s=\binom{m+s-1}{m}=
\frac{s^m}{m!}+O(s^{m-1}).
$$
By Riemann-Roch,
$$
 H^0(X,\o_X(tM))=\frac{(M^m)}{m!}t^m+O(t^{m-1}).
$$
 Choose $t\gg 1$ and $s$
such that $\root{m}\of{(M^m)}>s/t>\root{m}\of{(M^m)}-\epsilon$. Let
$D(s,t,x)$ be the zero set of a nonzero section of $m^s\otimes \o_X(tM)$ and
$D(x,\epsilon)=D(s,t,x)/t$. By construction $\mult_xD(x,\epsilon)\geq s/t$, as
required.\qed\enddemo

The above divisor $D(x,\epsilon)$ has high multiplicity at $x$, but we cannot
guarantee that it has low multiplicity elsewhere. The following example shows
that, even for surfaces, forcing high multiplicity at one point can cause high
multiplicities to appear at other points.

\demo{6.2 Exercises} (6.2.1) Let $S=\p^1\times \p^1$ and $M=\o(1,m)$. Then
$(M^2)=2m$, so by (6.1) for  any point $x=(a,b)\in S$ and $d<\sqrt{2m}$ there
is a
$D_x\equiv M$ such that $\mult_xD_x>d$. Show that $D_x$ contains the curve 
$\p^1\times
\{b\}$ with multiplicity at least $d-1$.

(6.2.2) For any $m,n, d>0$ construct a smooth surface $S$  and
an irreducible curve
$C\subset S$  such that $(C^2)=d$ and there are two points $p,q\in C$ such
that $\mult_pC=m$ and $\mult_qC=n$.

Assume that $m^2>d$. 
Let $D=aC+C'$ ($C\not\subset \supp C'$) be a $\q$-divisor such that $D\equiv C$
and
$\mult_pD\geq r$.  $D=C$ is such a divisor, but we want to understand all such
divisors as well.
 Show that $a\geq
(mr-d)/(m^2-d)$.  Thus if $mr-d>0$ then any such divisor has $C$ as an
irreducible component.

Assume now that $m^2>d$ and
$n>2(m^2-d)/(m-d)$. Then any   $\q$-divisor $D$ which has multiplicity at
least 1 at $p$ has multiplicity at least 2 at $q$. Thus there is no way to
make $D$ not klt at $p$ without making it much worse at $q$.
\enddemo

In order to illustrate the techniques involved, I first prove the surface
version of (6.4). The proof is set up to emphasize the general methods, and 
it does not give the optimal  result for surfaces.

\proclaim{6.3 Example-Theorem} Let $S$ be a normal projective surface, $L$ an
ample
 $\q$-divisor on $S$ and $x\in S$ a smooth point.  Set
$a^2=(L^2)$ and $b=\min_{x\in C} (L\cdot C)$.
Assume that
$$
1> \frac2{a}+\frac1{b}.
$$
Then there is an effective $\q$-divisor $D\equiv L$ such that

(6.3.1) $D$ is not log canonical at $x$, and

(6.3.2) $D$ is klt in a punctured neighborhood of $x$.
\endproclaim

\demop   The proof is in three steps, corresponding to the 3 main steps
(6.7--9) of the higher dimensional  argument.

(6.3.3) Step 1.  Construction of a divisor which is singular at $x$.

Choose $c_1>2/a$. Then $(c_1L)^2>4$, thus by (6.1) there is a divisor
$D_1\equiv c_1L$ such that $\mult_xD_1>2$. Thus $D_1$ is not log
canonical at
$x$ (3.9.4).
If $D_1$ is klt (or even lc) in a punctured neighborhood of $x$, then  go to
step 3.

We are left with the case when $D_1$ is not lc in any punctured neighborhood
of $x$.  Write $D_1=\sum e_iE_i$  where the $E_i$ are irreducible curves and
set $e=\max\{e_i\vert x\in E_i\}$. By assumption $e> 1$.  Let 
$$
D_2=\frac1{e}D_1=\sum \frac{e_i}{e}E_i\qtq{and} C=\sum_{i\:e_i=e, x\in E_i}
E_i.
$$
Then $D_2\equiv (c_1/e)L$ and $D_2-C$ is an effective $\q$-divisor where each
irreducible component containing $x$ has coefficient less than 1. 

If $x$ is a singular point of $C$, then $D_2$ is not klt at $x$, again go to
step 3. Thus we
are left with the case when $x$ is a smooth point of $C$, in particular, $C$ is
irreducible.

(6.3.4) Step 2. Induction on the dimension.

Let $c_2>1/b$ and choose $n\gg 1$ such that $nc_2L$ is Cartier.
Then 
$$
\deg_C(nc_2L)=nc_2\deg_C(L)\geq nc_2 b >n+2g(C)-1\qtq{for $n\gg 1$,}
$$
 thus $\o_C(nc_2L|C)$  has a section 
$s_C$ which vanishes at $x$ to order $n$. That is, $(1/n)(s_C=0)$ is a
divisor on $C$ which is numerically equivalent to $c_2L|C$ and which is not
klt at $x$. 

We may assume that
$H^1(S,\o_S(nc_2L)(-C))=0$, thus $s_C$ can be lifted to a section $s_S$ of
$\o_S(nc_2L)$.  By generic choice of $s_S$ we may assume that $s_S$ does not
vanish along any irreducible component of $D_1$.  Let $D'_1=(1/n)(s_S=0)$.  
Then
$$
D'_1+D_2\equiv ((c_1/e)+c_2)L,\qtq{and it is not klt at $x$ by (7.3.2).}
$$
We can choose $c_1, c_2$ such that $(c_1/e)+c_2<1$.

(6.3.5) Step 3. Tie breaking.

The previous steps frequently yield a divisor 
 $D_1\equiv cL$ for some $c<1$ such that
$D_1$  is not klt at
$x$ and it is
lc in a punctured neighborhood of $x$.  We show that a small perturbation of
$D_1$ gives the required $D$.

Choose $n\gg 1$ such that $n(1-c)L$
is Cartier and very ample. Let 
$D'_1\in |n(1-c)L|$ be a general divisor  passing through $x$ and set
$D_2=D_1+(1/n)D'_1$. Then $D_2$  is not lc at $x$ but lc 
in a punctured neighborhood of $x$.  Choose $m\gg 1$ such that $mL$
is Cartier and very ample. Let 
$D'_2\in |mL|$ be a general divisor. Then for $0<\delta\ll 1$,
$$
D=(1-\delta)D_2+\frac{\delta}{m}D'_2
$$
satisfies the
requirements of (6.3).
\qed\enddemo

\demo{6.3.7 Remarks}   

(6.3.7.1) Observe that the multiplicity of the  divisor $D$ at $x$ does not
necessarily predict that it is not lc at $x$. In step 2, the best lower bound
for the multiplicity is $\mult_xD\geq 1+(1/m)>1$.

(6.3.7.2) The proof of (6.4) proceeds along the same lines.
First we find a very singular divisor, and then we try to correct it, improving
things one dimension at a time. There are some technical problems. In the
surface case, if a divisor is klt at a point, it is smooth. In higher
dimensions this is not true, and the main technical innovation of
\cite{Angehrn-Siu95} is to figure out how to deal with the resulting
singularities.
\enddemo

The main  result of this section is the following:

\proclaim{6.4 Theorem} Let $(X,\Delta)$ be a projective klt pair and  $N$ 
a nef and big $\q$-Cartier $\q$-divisor on $X$. Let $x\in X$ be  a closed
point and $c(k)$ positive numbers such that   if
$x\in  Z\subset X$  is an irreducible (positive dimensional) subvariety  then 
$$
(N^{\dim Z}\cdot Z)>c(\dim Z)^{\dim Z}.
$$ 
 Assume that 
$$
\sum_{k=1}^n \frac{k}{c(k)}\leq 1.
$$
Then there is an effective $\q$-divisor $D\equiv N$ and an open
neighborhood $x\in X^0\subset X$ such that

(6.4.1) $(X^0,\Delta+D)$ is lc;

(6.4.2) $(X^0,\Delta+D)$ is klt on $X^0-x$;

(6.4.3) $(X,\Delta+D)$ is not klt at $x$.
\endproclaim

In order to separate points  by global sections we need a version of the
above result with two points. One might try to find a divisor
which is lc at two given points and klt in a neighborhood of them. 
This is impossible in general (6.2.2). The following
proof gives  a  weaker result which, however, is sufficient for our
purposes.

\proclaim{6.5 Theorem} Notation as above.  Let $x,x'\in X$ be  closed
points  such that   if
$ Z\subset X$  is an irreducible (positive dimensional) subvariety such that
$x\in Z$  or $x'\in Z$ then 
$$
(N^{\dim Z}\cdot Z)>c(\dim Z)^{\dim Z}.
$$ 
 Assume also  that 
$$
\sum_{k=1}^n \root{k}\of{2}\frac{k}{c(k)}\leq 1.
$$
Then, possibly after switching $x$ and $x'$,  one can choose
$D$ as above such that in addition to (6.4.1--3) it also satisfies:

(6.5.1)  $(X,\Delta+D)$ is not lc at $x'$.
\endproclaim
$$
*****
$$
\demop  The proof has several steps. Many of the intermediate results
are of interest in their own right.
\demop

\demo{6.6 Step 0}  Reduction to  $N$  ample. 

This step relies on two lemmas:

\proclaim{6.6.1 Lemma} Let $X$ be a proper scheme, $N$ a nef  and big divisor
on
$X$. Let
$x\in X$ be a point and assume that there are numbers $c(k)>0$ such that if
$x\in Z\subset X$  is an irreducible subvariety then $(N^{\dim Z}\cdot
Z)>c(\dim Z)$. 

Then we can write $N\equiv M+F$ where $M$ is ample, $F$ is effective and very
small and if
$x\in Z\subset X$  is an irreducible subvariety then $(M^{\dim Z}\cdot
Z)>c(\dim Z)$. 
\endproclaim

\demop  Write $N=A+E$ where $A$ is ample and $E$ effective.
Set 
$A_{\epsilon}=(1-{\epsilon})N+{\epsilon}A$. Then  $A_{\epsilon}$ is  ample
and $N=A_{\epsilon}+\epsilon E$. 
Furthermore, if $Z\subset X$ is a $k$-dimensional irreducible subvariety then
$$
(A_{\epsilon}^k\cdot Z)=(1-{\epsilon})^k(N^k\cdot Z)+\sum_{i=0}^{k-1}
\epsilon(1-{\epsilon})^i(A\cdot
N^i\cdot A_{\epsilon}^{k-1-i}\cdot Z) \geq (1-{\epsilon})^k(N^k\cdot Z).
$$

This says that a nef and big  divisor can be approximated by ample
ones with uniform control over intersection numbers.
We are done if we can exclude the possibility that $\inf_Z (N^k\cdot Z)=c(k)$.
This is implied by (6.6.2).
\qed\enddemo

\proclaim{6.6.2 Lemma}  Let $X$ be a proper scheme, $N$ a nef  and big divisor
on
$X$. Let
$x\in X$ be a point and assume that  if
$x\in Z\subset X$  is an irreducible subvariety then $(N^{\dim Z}\cdot Z)>0$. 

For every constant $C>0$ there are only finitely many families of 
irreducible subvarieties $x\in Z\subset X$  such that 
$(N^{\dim Z}\cdot Z)<C$.
\endproclaim

\demop Write $N=A+E$ where $A$ is ample and $E$ effective.
If $Z\not\subset E$ and $\dim Z=k$ then 
$$
(N^k\cdot Z)=(A^k\cdot Z)+\sum_{i=0}^{k-1}(E\cdot A^i\cdot
(A+E)^{k-1-i}\cdot Z) \geq (A^k\cdot Z),
$$
 and there are only  finitely many
families of $k$-dimensional irreducible subvarieties of $X$ such that 
$(A^k\cdot Z)<C$.  By induction on the dimension, there are 
only 
finitely many families of $k$-dimensional
irreducible subvarieties of $E$ containing $x$ such that 
$(N^k\cdot Z)<C$.
\qed\enddemo

By (6.6.1) we
can write
$N=M+F$ where $M$ is ample, satisfies the assumptions of (6.4) and we
can choose $F$  small enough such that $(X,\Delta+F)$ is still klt. Find
$D'\equiv M$ as required and then set $D=D'+F$.\qed\enddemo

\demo{6.7 Step 1} Finding a singular divisor at  $x$.

\proclaim{6.7.1. Theorem} Let  $(X,\Delta)$ be  klt,  projective of
dimension $n$ and 
$x\in X$ a closed point.  Let
$H$ be an ample
$\q$-divisor on
$X$ such that $(H^n)>n^n$. Then  there is an
effective
$\q$-divisor
$B_x\equiv H$  such that 
 $(X,\Delta+B_x)$ is not lc at  
 $x$.
\endproclaim

\demop If $x$ is smooth, this follows directly from (6.1).  Moreover, from the
proof we see that there is an $m>0$ 
(depending on $(X,\Delta)$ and $H$ but not on $x$)
such that we can choose
$B_x=(1/m)D_x$ where $D_x\in |mH|$. 

In general, let $0\in C$ be a smooth affine  curve and $g:C\to X$ a
morphism such that $x=g(0)$ and $g(c)\in X$ is smooth for  $0\neq c\in C$. For
 $c\neq 0$ pick a $\q$-divisor 
$B_{g(c)}$  as above such that $B_{g(c)}$ is not lc at $g(c)$.   It is natural
to take  
$B_x\deq \lim_{c\to 0} B_{g(c)}$. Limits of $\q$-divisors do not make too
much sense in general (except as currents), but in our case one can attach a
clear and precise meaning.  By our construction, $B_{g(c)}=(1/m)D_{g(c)}$
where the
$D_{g(c)}$ are Cartier divisors from the same linear system $|mH|$. After
 passing to a finite cover of $C$, we may assume that $g$ lifts to a morphism
$\tilde g:C\setminus\{0\}\to |mH|$ such that 
$D_{g(c)}=\tilde g(c)$. Thus we can take $B_x\deq (1/m)\tilde g(0)$.

By (7.8), $B_x$  is not lc  at $x$.\qed\enddemo
\enddemo

\demo{6.8 Step 2} Inductive step.

The main part of the proof is the following:

\proclaim{6.8.1. Theorem} Let  $(X,\Delta)$ be  klt,  projective and 
$x\in X$ a closed point. Let
$D$ be an effective $\q$-Cartier $\q$-divisor on $X$ such that  $(X,\Delta+D)$
is lc in a neighborhood of $x$. Assume that 
$\nklt(X,\Delta+D)= Z\cup Z'$ where $Z$ is irreducible, $x\in Z$ and $x\not\in
Z'$. Set 
$k=\dim Z$.  Let
$H$ be an ample
$\q$-divisor on
$X$ such that $(H^k\cdot Z)>k^k$. Then  there is an
effective
$\q$-divisor
$B\equiv H$ and rational numbers $1\gg \delta>0$ and $0<c<1$ such that 

(6.8.1.1) $(X,\Delta+(1-\delta)D+ cB)$ is lc 
in a neighborhood of $x$, and

(6.8.1.2) 
$ \nklt(X,\Delta+(1-\delta)D+ cB)=Z_1\cup Z'_1$, $x\in Z_1, x\not\in Z'_1$ and 
$\dim Z_1<\dim Z$.
\endproclaim

\proclaim{6.8.1.3 Complement} 
 Assume  in addition that $(X,\Delta+D)$
is not lc at $x'\in X$. Then we can choose $B,\delta$ and $c$  such that
$(X,\Delta+(1-\delta)D+ cB)$ is not lc at $x'$.
\endproclaim

\demop By assumption there is a proper birational morphism $f:Y\to X$ and a
 divisor $E\subset Y$ such that $a(X,\Delta+D,E)=-1$ and $f(E)=Z$. 
Write $K_Y\equiv f^*(K_X+\Delta+D)+\sum e_iE_i$ where $E=E_1$ and so $e_1=-1$. 
Let $Z^0\subset Z$ be an open subset such that

(6.8.2.1)  $f|E:E\to Z$ is smooth over $Z^0$, 
and 

(6.8.2.2)  if $z\in Z^0$ then  $(f|E)^{-1}(z)\not\subset E_i$  for $i\neq 1$.

The  following claim essentially proves (6.8.1) in case $x\in Z^0$.

\proclaim{6.8.3 Claim}  Notation as above. 
Choose  $m\gg 1$ such that $mH$ is Cartier. Then for every $z\in Z^0$ the
following assertions hold:

(6.8.3.1) There is a divisor $ F_z\sim mH|Z$ such that
$\mult_zF_z>mk$.

(6.8.3.2) $\o_X(mH)\otimes I_Z$ is generated by global sections and 
$H^1(X,\o_X(mH)\otimes I_Z)=0$. In particular
$H^0(X,\o_X(mH))\to H^0(Z,\o_Z(mH|Z))$ is surjective. 

(6.8.3.3) For any $F\sim mH|Z$ there is $F^X\sim mH$ such that
$F^X|Z=F$ and $(X,\Delta+D+(1/m)F^X)$ is klt  on $X-(Z\cup Z')$.

(6.8.3.4) Let $F_z^X\sim mH$ be such that $F_z^X|Z=F_z$. Then 
$(X,\Delta+D+(1/m)F_z^X)$ is not lc at $z$.
\endproclaim

\demop (6.8.3.1) is the usual multiplicity estimate (6.1),  and (6.8.3.2) is
satisfied once
$m\gg 1$. 

Let $|B|\subset |mH|$ be the linear subsystem consisting of those divisors 
$B'$ such that either $Z\subset B'$ or $B'|Z=F$. By (6.8.3.2) $|B|$ is base
point free on $X-Z$. Thus (4.8.2) implies (6.8.3.3).

Finally consider (6.8.3.4). Let  $y$ be the generic point of $(f|E)^{-1}(z)$.
Write
$$
\align
&K_Y\equiv f^*(K_X+\Delta+D)+\sum e_iE_i,\qtq{where $ E=E_1$, and}\\
& f^*F_z^X=F_z^Y+\sum mf_iE_i,\qtq{where $F_z^Y=f^{-1}_*F_z^X$. Thus}\\
&K_Y+(1/m)F_z^Y+\sum (f_i-e_i)E_i\equiv f^*(K_X+\Delta+D+(1/m)F_z^X).
\endalign
$$
 $(X,\Delta+D+(1/m)F_z^X)$ is not lc at $z$ if
$(Y,(1/m)F_z^Y+\sum (f_i-e_i)E_i)$ is not lc at $y$. $Z\not\subset F_z^X$,
and therefore $f_1=0$. Thus $\sum (f_i-e_i)E_i=E+\sum_{i\neq 1} (f_i-e_i)E_i$
and by assumption none of the $E_i$ for $i\neq 1$ contains $y$. Thus 
$(Y,(1/m)F_z^Y+\sum (f_i-e_i)E_i)$ is not lc at $y$ iff 
$(Y,(1/m)f^*F_z^X+E)$ is not lc at $y$. By (7.5.2) the
latter holds iff
 $(E, (1/m)f^*F_z^X|E=(1/m)(f|E)^*(F_z))$  is not lc at $y$. 
$E$ is smooth at $y$, $y$ has codimension $k$ in $E$ and 
$(1/m)(f|E)^*(F_z)$ has multiplicity $>k$. Thus  
$(E, (1/m)(f|E)^*(F_z))$  is not lc at $y$.\qed\enddemo

Next we intend to show that  by continuity, there are divisors $F_z^X$
as in (6.8.3) even if $z\in Z-Z^0$.

 Pick $z_0\in Z$
arbitrary. Let $0\in C$ be a smooth affine  curve and $g:C\to Z$ a
morphism such that $z_0=g(0)$ and $g(c)\in Z^0$ for general $c\in C$. For
general $c\in C$ pick
$F_c\deq F_{g(c)}$ as in (6.8.3.1). Let
$F_0=\lim_{c\to 0} F_c$. (The limit is defined as at the end of Step 1.)

\proclaim{6.8.4 Claim} Notation as above.  There is a divisor $F_0^X\in
|mH|$ such that

(6.8.4.1) $F_0^X|Z=F_0$,

(6.8.4.2) $(X,\Delta+D+(1/m)F_0^X)$
is klt  on $X-(Z\cup Z')$,

(6.8.4.3) $(X,\Delta+D+(1/m)F_0^X)$ is lc at the generic point of $Z$,

(6.8.4.3) $(X,\Delta+D+(1/m)F_0^X)$ is not lc at $z_0$.
\endproclaim

\demop By (6.8.3.3) we can find   $F_0^X$ such that (6.8.4.1--2) are
satisfied.
$F_0^X$ does not contain $Z$, thus  $(X,\Delta+D+(1/m)F_0^X)$ is lc at the
generic point of $Z$.  

We can lift $F_0^X$ to a family $F_c^X\: c\in C$ such that  $F_c^X|Z=F_c$. If 
$(X,\Delta+D+(1/m)F_0^X)$ is lc at $z_0$ then  by (7.8) 
$(X,\Delta+(D+(1/m)F_c^X)$ is lc in a neighborhood of $z_0$
for general $c\in C$. This, however, contradicts (6.8.3.4).\qed\enddemo

To finish the proof of (6.8.1) set $B=(1/m)F_0^X$. 
$(X,\Delta+(1-\delta)D)$ is klt at the generic point of $Z$ for every 
 $\delta>0$.   Choose
$1\gg \delta>0$ such that   
$(X,\Delta+(1-\delta)D+B)$ is not lc at $z_0$ and then $0<c<1$
such that $(X,\Delta+(1-\delta)D+cB)$ is   lc but not klt at $z_0$.

If $(X,\Delta+D)$ is not lc at $x'$ then the same holds  for
$(X,\Delta+(1-\delta)D)$ for $1\gg\delta>0$ and any choice of $c$
preserves this property. \qed\enddemo\enddemo

(6.8.1) is nearly enough to prove (6.4) by induction. The only problem is
that in (6.8.1) we may end up with $(X,\Delta+(1-\delta)D+cB)$
such that $\nklt(X,\Delta+(1-\delta)D+cB)$ has several irreducible
components passing through $z_0$. This is taken care of by the next step.

\demo{6.9 Step 3} Tie breaking.

\proclaim{6.9.1 Lemma}  Let  $(X,\Delta)$ be  klt,  projective and
$x\in X$ a  point. Let $D$ be an effective $\q$-Cartier $\q$-divisor on $X$
such that 
$(X,\Delta+D)$ is lc in a neighborhood of $x$.
Let
$\nklt(X,\Delta+D)=\cup Z_i$ be the irreducible components; $x\in Z_1$.  Let
$H$ be an ample $\q$-divisor on
$X$. Then for every $1\gg \delta>0$   there is an effective
$\q$-divisor
$B\equiv
 H$   and $0<c<1$ such that 

(6.9.1.1) $(X,\Delta+(1-\delta)D+  cB)$ is
lc in a neighborhood of $x$, and

(6.9.1.2) 
$\nklt(X,\Delta+(1-\delta)D+  cB)=W\cup W'$ where $x\in W, x\not\in W'$ and
$W\subset Z_1$.
\endproclaim

\demop  Choose $m\gg 1$  such that $mH$ is Cartier and so that $\o_X(mH)\otimes
I_{Z_1}$ is generated by global sections. Let $B'\in |\o_X(mH)\otimes I_{Z_1}|$
be a general member. By (4.8.2)
$(X,\Delta+(1-\delta)D+b B')$ is klt outside
$Z_1$ in a neighborhood of $x$ for $b< 1$. It is definitely not lc  along
$Z_1$ for $1>b\gg \delta>0$. 
First choose
$b= 1/m$ and $1/m\gg \delta>0$. Then choose $0<c<1$ such that
$(X,\Delta+(1-\delta)D+(c/m) B')$ is lc but not   klt at $x$. Set
$B=B'/m$. \qed\enddemo
\enddemo

\demo{6.10 Step 4} Proof of (6.4).

We prove by induction the following theorem. The case $j=\dim X$ gives
(6.4):

\proclaim{6.10.1 Theorem} Notation and assumptions as in (6.4). Assume in
addition that $N$ is ample.  Let $j\in \{1,\dots,n\}$. 
Then for every $b\geq \sum_{k=n-j}^nkc(k)^{-1}$ there is an effective
$\q$-divisor
$D_j\equiv bN$ and an open neighborhood $x\in X^0\subset X$ such that

(6.10.1.1) $(X^0,\Delta+D_j)$ is lc;

(6.10.1.2) $\codim(\nklt(X^0,\Delta+D_j), X^0)\geq j$;

(6.10.1.3) $(X,\Delta+D_j)$ is not klt at $x$.
\endproclaim

\demop Set $D_0=\emptyset$.  Assume that we already
found  $D_j$ and we would like to get $D_{j+1}$. 

If $j=0$, then apply (6.7.1).  If $j>0$ then 
by (6.9) for every $\epsilon>0$ there is a divisor $D'_j\equiv
(1-\delta)D_j+\epsilon N$ such that $D'_j$ satisfies (6.10.1.1--3)
and in addition either $Z\deq \nklt(X^0,\Delta+D_j)$ is irreducible of
codimension at least $j$ or it has codimension at least $j+1$. In the latter
case we can take $D_{j+1}=D'_j+\alpha M$ where $M$ is a general member of
$|mN|$ for $m\gg 1$ and $\alpha$ is suitable. (We may assume that
$\epsilon<(j+1)c(j+1)^{-1}$.)

In the first case set 
$H=((j+1)c(j+1)^{-1}-\epsilon)N$. For $0<\epsilon\ll 1$ we have that
$(H^j\cdot Z)>j^j$. Thus we can apply  (6.8.1) to obtain $D_{j+1}$. 
\qed\enddemo

\enddemo

\demo{6.11 Step 5} Proof of (6.5).

The proof is very much similar to the proof of (6.4), I just outline the
necessary changes. As before, we assume that $N$ is ample.

As a first step of the induction we take
$H=(\root{n}\of{2}nc(n)^{-1}-\epsilon)N$. Then $(H^n)>2n^n$ and  as in (6.1) we
can find a divisor $D_1\equiv H$ such that $\mult_xD_1>n$ and
$\mult_{x'}D_1>n$. This shows that $(X,\Delta+D_1)$ is not lc at the points
$x,x'$. Choose $0<c<1$ such that $(X,\Delta+cD_1)$ is not klt at the points
$x,x'$ and is lc at one of them, say at $x$. There are two separate
cases to consider.

(6.11.1) Case 1. If $(X,\Delta+cD_1)$ is not lc at $x'$ then (6.8.1) can be
used to continue exactly as in Step 4 to complete the proof.

(6.11.2) Case 2. What if $(X,\Delta+cD_1)$ is  lc at $x'$? Then first we apply
the tie-breaking method, to reduce to the case when $\nklt(X,\Delta+cD_1)$ is
irreducible in an open set containing $x$ and $x'$. The tie-breaking may put us
in the first case.  Otherwise we are in the situation when
$Z=\nklt(X,\Delta+cD_1)$ is irreducible near $x$ and contains $x'$. 

We then proceed as in (6.8.1) but instead of trying to force  high
multiplicity  at one point only, we do it at two points. Only the notation
has to be changed. Proceeding inductively as in (6.10) we obtain (6.11).

If we always end up in case 2, then we may get a divisor $D_n$ such that both
$x$ and $x'$ are isolated points of $\nklt(X,\Delta+D_n)$ and $(X,\Delta+D_n)$
is lc at both points. For the purposes of (5.9) this is not a problem at all
(and may even be an advantage in general). In this case we can do a last
tie breaking to achieve exactly (6.5).
\qed\enddemo
\enddemo
\enddemo

\head 7. The $L^2$ Extension Theorem and Inversion of Adjunction
\endhead

Let $X$ be a variety and $S\subset X$ a Cartier divisor. If we know
something about the singularities of $S$ then we can usually assert that the
singularities of $X$ near $S$ are not worse. For instance, if $S$ is smooth,
or rational, or CM then the same holds for $X$. In some cases 
the converse implication also holds. This fails for smooth or rational, but
works for CM.

Our aim is to investigate the analogous problem for discrepancies.
It can be formulated in two variants.  This first one  has been
considered in complex analysis. The second version is the natural one from the
algebraic geometry point of view:  When talking about
$(X,B)$, we always compute with $K_X+B$, and $K_X+S+B|S=K_S+B$  by
adjunction. 

\proclaim{7.1 Questions}  Let $0\in S\subset X$ be a Cartier divisor.
 Let $B=\sum b_iB_i$ be a $\q$-divisor such that $S$ is not among the
$B_i$. 

(7.1.1) If $(S,B|S)$ is lc, is then  $(X, B)$ also lc?

(7.1.2) If $(S,B|S)$ is lc, is then  $(X,S+B)$ also lc?
\endproclaim

It is a priori clear that the second form is stronger, but it turns out that
there is no real difference between them:

\proclaim{7.1.3 Lemma} \cite{Manivel93} The above two questions are equivalent.
\endproclaim

\demop  We need to prove that the first version implies the second one.

Let $S=(f=0)$ and let $X_n\subset X\times \a^1$ be given by the
equation
$y^n-f=0$, where $y$ is the coordinate on $\a^1$. 
Let $p_n:X_n\to X$ be the projection and
set $B_n\deq p_n^{-1} B$. $S\cong S_n\deq (y=0)$ appears as a Cartier divisor
on
$X_n$ and $B_n|S_n= B|S$ under this isomorphism.

By (7.1.1),  $(X_n,B_n)$ is lc.  Observe that
$$
K_{X_n}+B_n+(n-1)S_n=p_n^*(K_X+B),\qtq{thus}
K_{X_n}+B_n\equiv p_n^*(K_X+B+(1-1/n)S).
$$
From (3.16.2) we conclude that $(X,B+(1-1/n)S)$ is lc for every $n$,
thus $(X,B+S)$ is lc.\qed\enddemo

The first significant result toward answering (7.1) is the 
$L^2$ extension theorem of \cite{Ohsawa-Takegoshi87}, though the connection
was first realized only later. I state a form of the theorem which is natural
from the point of view of complex analysis.
Instead of  defining the notions ``pseudoconvex" and ``plurisubharmonic",  it
is sufficient to keep  two special cases in mind:

--- every convex subset $\Omega\subset\c^n$
is pseudoconvex,

--- if $g$ is holomorphic then $c\log|g|$ is  plurisubharmonic for $c>0$.

\proclaim{7.2 Theorem} \cite{Ohsawa-Takegoshi87} Let $\Omega\subset\c^n$
be  a bounded pseudoconvex domain, and  $H\subset\c^n$ a hyperplane
intersecting $\Omega$. 
Fix Lebesgue measures $dm_n$ on $\c^n$ and $dm_{n-1}$ on  $H$. 
Then
there is a constant $C_{\Omega}$   with the following
property.

Let  $\phi$  be plurisubharmonic on $\Omega$ and $f$ holomorphic on
$\Omega\cap H$.  Then $f$ can be extended to a holomorphic function $F$ on
$\Omega$  such that
$$
\int_\Omega|F|^2e^{-\phi}dm_n\leq C_{\Omega}
\int_{\Omega\cap H}|f|^2e^{-\phi}dm_{n-1}.\qed
$$
\endproclaim

The following consequence relates this to (7.1):

\proclaim{7.2.1 Corollary}  Let $0\in H\subset\c^n$ be a hyperplane.

(7.2.1.1)  Let  $g$ be
a holomorphic function near $0$ and let $g_H$ denote the restriction of $g$
to $H$. If
$|g_H|^{-c}$ is
$L^2$ near
$0$ then
$|g|^{-c}$ is $L^2$ near $0$.

(7.2.1.2) Let $B=\sum b_iB_i$ be a $\q$-divisor such that $H$ is not among the
$B_i$. If $(H,B|H)$ is klt (resp. lc) then $(\c^n, B)$ is klt (resp. lc).
\endproclaim

\demop  For $\Omega$ choose a small ball around $0$. Pick $\phi=2c\log|g|$ and
$f\equiv 1$.  We do not know what $F$ is, but $|F|\geq 1/2$ in a neighborhood
 $0\in \Omega'\subset \Omega$. Thus
$$
\frac12\int_{\Omega'}|g|^{-2c}dm_n\leq
\int_\Omega|F|^2|g|^{-2c}dm_n\leq C_{\Omega}
\int_{\Omega\cap H}|g|^{-2c}dm_{n-1}.
$$
This shows the first part. In order to see the second part, choose $m$ such
that $mb_i$ are all integers and let $g$ be a function with $(g=0)=\sum
(mb_i)B_i$. Applying (7.2.1.1) with $c=1/m$ gives the klt
part of (7.2.1.2) by (3.20).

Finally $(\c^n,B)$ is lc iff $(\c^n,(1-\epsilon)B)$ is klt for every
$\epsilon>0$, thus the klt case implies the lc version.
\qed\enddemo

\demo{7.2.2 Remark} The application of the $L^2$ extension theorem for the
constant function seems quite silly. After all,  we need only that
if $|g_H|^{-c}$ is  $L^2$ then so is $|g|^{-c}$. A simple manipulation of
 integrals may   give this result, but such a proof is not yet known. 
\enddemo

Unaware of (7.2), \cite{Shokurov92,3.3} proposed 
 a conjecture along the lines of (7.1.1) for algebraic varieties. The
conjecture was  
 subsequently generalized in \cite{Koll\'ar et al.92, 17.3}.

The conjecture (or similar
results and conjectures) is frequently referred to as {\it adjunction}
(if we assume something about
$X$ and obtain conclusions about $S$)  or {\it inversion of adjunction} (if we
assume something about $S$ and obtain conclusions about $X$).

\proclaim{7.3 Conjecture}  Let $X$ be a normal
variety, $S$ a normal Cartier divisor  and   $B=\sum b_iB_i$  a
$\q$-divisor. Assume that $K_X+S+B$ is $\q$-Cartier. Then 
$$
\tdiscrep(S,B|S)= 
\discrep\left(\cent\cap S\neq\emptyset,X,S+B\right), 
$$
where the notation on the right means that we compute the discrepancy
using only those divisors whose center on $X$ intersects $S$. 
\endproclaim

\demo{7.3.1 Remarks} 

(7.3.1.1) The conjecture can also be formulated if $S\subset X$ is only a Weil
divisor. For some applications this is crucial, but for us it is not
necessary. Also, it leads to   additional difficulties. 
The point is that if $S$ is not Cartier, the usual adjunction formula
$K_S=(K_X+S)|S$ fails. A suitable correction term needs to be worked out. Once
this is settled, the proofs require little change.  The interested reader
should consult \cite{Koll\'ar et al.92}, especially  Chapters 16--17.

(7.3.1.2) \cite{Koll\'ar et al.92,17.12} shows that (7.3) is implied by the
logMMP.  Thus it is true if $\dim X\leq 3$.
Various special cases
of this conjecture are very useful in the proof of the logMMP and in many other
contexts, see, for instance, \cite{Koll\'ar et al.92,Ch.18} or
\cite{Corti94}.  Therefore it is rather desirable to find a proof independent
of the logMMP.
\enddemo

One inequality is easy to prove:

\proclaim{7.3.2 Proposition} \cite{Koll\'ar et al.92, 17.2} Let $X$ be a normal
variety, $S$ a  Cartier divisor  and   $B=\sum b_iB_i$  a
$\q$-divisor. Assume that $K_X+S+B$ is $\q$-Cartier.  Then
$$
\tdiscrep(S,B|S)\geq \discrep (X,S+B ). 
$$
\endproclaim

\demop (Strictly speaking, the left hand side is only defined if $S$ itself is
normal. The following proof furnishes a definition of the left hand side, 
which is the correct one for schemes $S$ which are Cartier divisors on a
normal scheme.)

Let  $f:Y\to X$ be a log resolution of $(X,S+B)$ and set
$S'\deq f^{-1}_*S$.  
We may also assume that $S'$  is disjoint from $f^{-1}_*B$.
Write $K_Y+S'\equiv f^*(K_X+S+B)+\sum e_iE_i$. By
the usual adjunction formula,
$$
K_{S'}=K_Y+S'|S',\qtq{and} K_X+S+B|S=K_S+B|S.
$$
This gives that 
$$
K_{S'}\equiv f^*(K_S+B|S)+\sum e_i (E_i\cap S').
$$
By assumption $S'$  is disjoint from $f^{-1}_*B$, thus if $E_i\cap S'\neq
\emptyset$ then $E_i$ is $f$-exceptional. 
This shows that every discrepancy
which occurs in $S'\to S$ also occurs among the
exceptional divisors of $Y\to X$.
It may of course happen that $E_i$ is $f$-exceptional but $E_i\cap S'$ is
not
$f|S'$-exceptional.  This is why we have $\tdiscrep$ on the left hand side
of the inequality.
\qed\enddemo

\demo{7.3.2.1 Remark} In general there are many exceptional divisors  $E_j$  of
$f:Y\to X$ which do not intersect $S'$, and there is no obvious connection
between the discrepancies of such divisors and the discrepancies occurring in
$S'\to S$. This makes the reverse inequality nonobvious.
\enddemo

For most of the applications of (7.3) the crucial case is when one of the
two sides is klt or lc.  The case when $S$ is smooth follows from (7.2.1).
The  singular cases are settled in 
\cite{Koll\'ar et al.92, Chapter 17}. The proof relies on the
following connectedness result which is of interest itself.

\proclaim{7.4  Theorem} \cite{Koll\'ar et al.92, 17.4} Let $X$ be
a normal variety 
(or analytic space)
and   $D=\sum d_iD_i$   an effective  $\q$-divisor on $X$ such that $K_X+D$
is $\q$-Cartier. Let  
$g:Y\to X$ 
be a log resolution of $(X,D)$. Write
$$
\align
&K_Y=g^*(K_X+D)+\sum e_iE_i,\qtq{and  set}\\
&A=\sum_{i:e_i>-1}e_iE_i,\qtq{and} F=-\sum_{i:e_i\leq-1}e_iE_i.
\endalign
$$

Then $\supp F=\supp\rdown{F}$ is connected in a neighborhood of any fiber of
$g$.  \endproclaim

\demop  By definition
$$
\rup{A}-\rdown{F}=K_Y+(-g^*(K_X+D))+\{-A\}+\{F\},
$$
and therefore by (2.17.3)
$$
R^1f_*\o_Y(\rup{A}-\rdown{F})=0.
$$
Applying $g_*$ to  the exact sequence
$$
0\to \o_Y(\rup{A}-\rdown{F})\to \o_Y(\rup{A})\to \o_{\rdown{F}}(\rup{A})\to 0
$$
we obtain that
$$
g_*\o_Y(\rup{A})\to g_*\o_{\rdown{F}}(\rup{A})
\qtq{is surjective.}\tag 7.4.1
$$
 Let $E_i$ be an irreducible component of $A$. Then either
$E_i$ is $g$-exceptional  or $E_i$ is the birational transform of some  $D_j$
whose coefficient in $D$ is less than 1. Thus $\rup{A}$ is $g$-exceptional and
$ g_*\o_Y(\rup{A})=\o_X$.
Assume that $\rdown{F}$ has at least two
connected components  $\rdown{F}=F_1\cup F_2$ in a neighborhood of
$g^{-1}(x)$ for some
$x\in X$. Then 
$$
g_*\o_{\rdown{F}}(\rup{A})_{(x)}\cong
g_*\o_{F_1}(\rup{A})_{(x)}+g_*\o_{F_2}(\rup{A})_{(x)},
 $$
and neither of these summands is zero. Thus $g_*\o_{\rdown{F}}(\rup{A})_{(x)}$
cannot
be the quotient of the  cyclic module 
$\o_{x,X}\cong g_*\o_Y(\rup{A})_{(x)}$.\qed\enddemo

As a corollary we obtain the following results which were proved 
by \cite{Shokurov92} in dimension 3 and by 
\cite{Koll\'ar et al.92, 17.6--7} in general.

\proclaim{7.5 Theorem}  Let $X$ be normal and
 $S\subset X$  an irreducible Cartier divisor. Let $B$ be an effective
$\q$-divisor and assume that $K_X+S+B$ is $\q$-Cartier.  
Then  

(7.5.1) 
$$
(X,S+B)\qtq{is  plt near $S$} \Leftrightarrow\quad (S,B|S)\qtq{is 
klt.}
$$ 

(7.5.2) Assume in addition that $B$ is $\q$-Cartier and $S$ is klt. Then
$$
(X,S+B)\qtq{is  lc near $S$} \Leftrightarrow\quad (S,B|S)\qtq{is 
lc.}
$$ 
\endproclaim

\demop  In both cases the implication $\Rightarrow$ follows from (7.3.2). 

In order to see (7.5.1), 
let
$g:Y@>>> X$
be a resolution of singularities and as in (7.4) let
$$
K_Y\equiv g^*(K_X+S+ B)+A-F.
$$
Let $S'=g^{-1}_*S$ and   $F=S'+ F'$. 
By adjunction
$$
K_{S'}=g^*(K_S+ B|S)+(A-F')|S'.
$$
$(X,S+B)$ is plt iff $F'=\emptyset$ and
$(S,B|S)$ is plt iff $F'\cap S'=\emptyset$. 
By (7.4) $S'\cup F'$ is connected, hence $F'=\emptyset$  iff 
$F'\cap S'=\emptyset$. This shows (7.5.1). 

From the definition 
we see that 
$$
\alignat3
&(X,S+B)\text{ is lc}  && \qtq{iff} && (X,S+cB) \qtq{is plt for
$c<1$, and}\\
 &(S,B|S)\text{ is lc} && \qtq{iff} && (S,cB|S)
\qtq{is klt for
$c<1$.}\\
\endalignat
$$
Thus (7.5.1) implies (7.5.2). 
\qed\enddemo

\demo{7.5.3 Exercise} Notation and assumptions as above. If $(X,S+B)$ is
plt then $\rdown{B}=\emptyset$ in a neighborhood of $S$.
\enddemo

As a corollary we obtain that klt and lc are open conditions in flat families:

\proclaim{7.6 Corollary}  Let $(X,B)$ be a pair such that $K_X+B$ is
$\q$-Cartier and
$g:X\to C$ a flat morphism to a smooth pointed curve $0\in C$.  Let $x\in
X_0=g^{-1}(0)$ be a closed point. 

(7.6.1) Assume that $(X_0,B|X_0)$ is klt at $x$, respectively

(7.6.2) assume that $X_0$
is klt  at $x$ and $(X_0,B|X_0)$ is lc at $x$. 

\noindent Then there is an open
neighborhood
$x\in U\subset X$ such that $(U_c,B|U_c)$ is klt (resp. lc) for every $c\in C$.
\endproclaim

\demop Let $S=X_0$. By (7.5) $(X,S+B)$ is plt (resp. lc) in a
neighborhood $W$ of  $x$. In the first case,  
$\rdown{B}=\emptyset$ by (7.5.3).  Thus $(W-S,(S+B)|W-S)=(W-S,B|W-S)$  is klt
(resp. lc). By (7.7) there
is an open subset
$C^0\subset C$ such that $(W_c,B|W_c)$ is klt (resp. lc) for $c\in C^0$. 
 Set $U=W\cap g^{-1}(\{0\}\cup C^0)$ to conclude.\qed\enddemo

\proclaim{7.7 Proposition} \cite{Reid80, 1.13} Let $X$ be a scheme over a field
of characteristic zero, $D$ a $\q$-divisor and $|B|$  a base point free linear
system of Cartier divisors. Assume that $K_X+D$ is $\q$-Cartier.
Let  $B'\in |B|$ be a  general member. Then
$$
\discrep(B', D|B')\geq  \discrep (X,D).
$$
In particular, if $(X,D)$ is lc (resp. klt, canonical, terminal) then
$(B', D|B')$ is also lc (resp. klt, canonical, terminal).
\endproclaim

\demop Let $f:Y\to X$ be a log resolution of $(X,D)$ and
$C'\deq f^{-1}_*B'=f^*B'$. Then $g\deq f|C':C'\to B'$ is a log resolution of
$(B',D|B')$.  Write  $K_Y\equiv  f^*(K_X+D)+\sum e_iE_i$. Then
$$
\align
K_{C'}=K_Y+C'|C'&\equiv (f^*(K_X+D+B')+\sum e_iE_i)|C'\\
&=g^*(K_{B'}+D|B')+\sum
e_i(E_i|C').\qed
\endalign
$$
\enddemo

\proclaim{7.8 Corollary} Let $Y$ be a klt variety over $\c$ and $B_c:c\in C$ an
algebraic family of $\q$-divisors on $Y$ parametrized by a  smooth pointed 
curve
$0\in C$. Assume that $(Y,B_0)$ is  klt (resp. lc) at $y\in Y$. 

Then there is a Euclidean  open neighborhood $y\in W\subset Y$
such that $(Y,B_c)$ is klt (resp. lc) on $W$  for  $c\in C$ near $0$.
\endproclaim

\demop $B_c$ defines a $\q$-divisor $B$ on $X\deq Y\times C$. Let $g:X
\to C$ be the projection.  Apply (7.6)  to obtain
$(y,0)\in U\subset X$. There are 
Euclidean    neighborhoods $y\in W\subset Y$ and $0\in D\subset C$
such that $W\times D\subset U$.
\qed

In general $W$ cannot be chosen to be Zariski open.\enddemo

The following result, predating (7.3), settles another special case.

\proclaim{7.9 Theorem} \cite{Stevens88} Let $X$ be a normal variety such
that  $\omega_X$ is locally free. Let $S\subset X$ be a Cartier
divisor. Then
$$
S \qtq{is  canonical} \Leftrightarrow \quad (X,S) \qtq{is  canonical near
$S$.}
$$
\endproclaim

\demop  If $(X,S)$ is canonical then $S$ is canonical by (7.3.2). 
In order to see the converse, let  $f:Y\to X$ be a log resolution of
$(X,S)$ and set 
$T\deq f^{-1}_*S$.

We have exact sequences:
$$
\CD
0@>>> \omega_Y@>>> \omega_Y(T)@>>> \omega_T@>>> 0,\qtq{and}\\
0@>>> \omega_X@>>> \omega_X(S)@>r>> \omega_S.@.
%\hphantom{\qtq{and}}
\endCD
$$
(If $X$ is CM, then $r$ is surjective, but we do not need
this.)  Pushing the first sequence forward to $X$, there are natural maps:
$$
\CD
0@>>> f_*\omega_Y@>>>  f_*\omega_Y(T)@>>>  f_*\omega_T@>>> R^1f_*\omega_Y=0\\
@. @VVV @VbVV @VcVV @.\\
0@>>>\omega_X@>>> \omega_X(S)@>r>> \omega_S. @.
\endCD
$$
$b$ is an injection and  $c$ is an isomorphism since $S$ is canonical.
Therefore 
$r\circ b$ is surjective, thus by the Nakayama lemma, $b$ itself is
surjective near $S$. Thus
 $b$ is also an isomorphism near $S$. 
By shrinking $X$, we may assume that $b$ is surjective. 
Since $\omega_X(S)$ is locally free, the
pull back of
$b$ to
$Y$ gives a map
$$
f^*(\omega_X(S))\to \omega_Y(T).
$$
That is, $f^*(\omega_X(S))=\omega_Y(T)(-E)$  for some effective divisor $E$.
Thus
$$
K_Y+T=f^*(K_X+S)+E,
$$
and $(X,S)$ is canonical near $S$.\qed\enddemo

\demo{7.9.1 Remark} Looking at the proof we see that instead of assuming that
$\omega_X$ is locally free and $S$ is Cartier, it is sufficient to assume
that $\omega_X(S)$ is locally free. This is  important in some
applications \cite{Stevens88; Koll\'ar-Mori92,3.1}.
\enddemo

A study of the above commutative diagram also gives the following:

\proclaim{7.9.2 Corollary} \cite{Ein-Lazarsfeld96, 3.1} 
Let
$X$ be a smooth variety (or a variety with index one canonical singularities)
and 
$S\subset X$ a Cartier divisor with desingularization 
$f:\bar S\to S$. 
Then there exists an ideal sheaf $J\subset \o_X$ such that
$$
0\to \omega_X\to \omega_X(S)\otimes J\to f_*\omega_{\bar S}\to 0
\qtq{is exact, }
$$
and $\supp (\o_X/J)$ is precisely the set of points where $S$ 
is not canonical. \qed
\endproclaim

\head 8. The Log Canonical Threshold
\endhead

Let $(X,D)$ be a pair. In section 3 we used the notion of discrepancy to
attach various invariants to $(X,D)$ which provide a way of measuring how
singular $X$ and $D$ are. If $(X,D)$ is not log canonical, the values of these
invariants are
$-\infty$, hence they give very little information. The aim of this section is
to develop another invariant, which becomes nontrivial precisely when the
discrepancy is $-\infty$. 

\demo{8.1 Definition} Let $(X,\Delta)$ be an lc pair, $Z\subset X$ a
closed subscheme and $D$ an effective $\q$-Cartier divisor on $X$.
The {\it log canonical threshold} (or {\it lc-threshold}) of $D$ along $Z$
with respect to
$(X,\Delta)$ is the number
$$
c_Z(X,\Delta,D)\deq \sup\{c\vert (X,\Delta+cD) \text{ is lc in an open
neighborhood of
$Z$}\}.
$$
If $\Delta=0$ then we use $c_Z(X,D)$ instead of $c_Z(X,0,D)$. 
We frequently write $c_Z(D)$ instead of $c_Z(X,\Delta,D)$ if no
confusion is likely.   If $D=(f=0)$ then we also use  the notation
$c_Z(X,\Delta,f)$ and $c_Z(f)$.
\enddemo

(3.20) shows that the lc-threshold has an equivalent analytic definition:

\proclaim{8.2 Proposition} Let $Y$ be a smooth variety over $\c$, $Z\subset
Y$ a closed subscheme and $f$ a
nonzero regular function on $Y$.   Then
$$
c_Z(Y,0,f)= \sup\{c\: |f|^{-c} \text{ is locally $L^2$ near
$Z$}\}.\qed
$$
\endproclaim

The following  properties  are clear from the definition:

\proclaim{8.3 Lemma}  Notation as above. Then:

(8.3.1) $c_Z(X,\Delta,D)\geq 0$ and  $c_Z(X,\Delta,D)=+\infty$ iff $D=0$.

(8.3.2) $c_Z(X,\Delta,D)=\inf_{p\in Z} c_p(X,\Delta,D)$.

(8.3.3) If $D$ is a Weil divisor, then $c_Z(X,\Delta,D)\leq 1$.\qed
\endproclaim

\demo{8.4 Remark} There is a slightly more general situation where the above
definition also makes sense. Instead of assuming that $(X,\Delta)$ is lc, it
is sufficient to assume that $(X,\Delta)$ is lc on $X-\supp D$. In this  case
$c_Z(X,\Delta,D)$ is negative if $(X,\Delta)$ is not lc along $Z$.
\enddemo

Next we turn to various techniques of computing and estimating the
lc-threshold.
 We can rewrite (3.13) to give an effective computational method:

\proclaim{8.5 Proposition}  Let $(X,\Delta)$ be an lc pair, $Z\subset X$ a
closed subscheme and $D$ an effective $\q$-Cartier divisor on $X$.    Let
$p:Y\to X$ be a proper birational   morphism. Using the convention (3.3.2),
write
$$
K_Y\equiv p^*(K_X+\Delta)+\sum a_iE_i,\qtq{and} p^*D=\sum b_iE_i.
$$
 Then:

(8.5.1) 
$$
c_Z(X,\Delta,D)\leq \min_{i\:p(E_i)\cap Z\neq \emptyset}
\left\{\frac{a_i+1}{b_i}\right\}.
$$

(8.5.2) Equality holds if $\sum E_i$ is a divisor with normal crossings
only. In particular, $c_Z(X,\Delta,D)\in \q$. \qed
\endproclaim

(7.6) translates to an upper semicontinuity statement for the lc-threshold:

\proclaim{8.6 Lemma} Let $(X,\Delta)$ be a klt pair, $x\in X$ a closed point
and $D_t\:t\in C$ an algebraic family of effective  $\q$-Cartier divisors.
Pick a point $t_0\in C$. Then
$$
c_0(X,\Delta, D_{t_0})\leq c_0(X,\Delta, D_{t}),\qtq{for $t$ near $t_0$.}\qed
$$
\endproclaim

The log canonical threshold has been investigated earlier in different
contexts. Some of these are explained in sections 9--10. Recent interest
arose following \cite{Shokurov92}  who used various properties of
lc-thresholds in order to establish the existence of log flips in dimension
three.  He proposed a rather striking conjecture (8.8), and proved it for
surfaces. Later this was proved for threefolds in \cite{Alexeev93}. The
general case is still unknown.  Before formulating the conjecture, we need a
definition.

\demo{8.7 Definition} For every $n\in \n$ define a subset of $\r$ by
$$
\Cal T_n\deq
\{c_x(X,D)\vert \text{$x\in X$ is klt, $\dim X=n$ and $D$ is an effective Weil
divisor}\}.
$$
By (8.3.3) and (8.5.2) we see that $\Cal T_n\subset (0,1]\cap \q$.
\enddemo

\proclaim{8.8 Conjecture} \cite{Shokurov92} For every $n$, the set $\Cal T_n$
satisfies the ascending chain condition.
\endproclaim

\demo{8.8.1 Remark}  This conjecture is only one example of a series of
conjectures that assert the ascending or descending chain condition for
various naturally defined invariants coming from algebraic geometry.
See \cite{Shokurov92; Koll\'ar et al.92, Ch.18; Koll\'ar94; Alexeev94;
Ganter95} for further examples and for applications.
\enddemo

The rest of the section is devoted to various methods of computing the
lc-threshold in several examples, and to study (8.8) in those cases.
Most of the computations that we do are for   $X$   smooth. Thus
working analytically we consider only the case $X=\c^n$. Even this special
case of (8.8) is mysterious.

\demo{8.9 Example}   Let $f\in \c[[x,y]]$ be an irreducible power series. By
\cite{Igusa77}, 
$$
c_0(\c^2,f)=\frac1{m}+\frac1{n}\qtq{where $m=\mult_0f$ and $n/m$ is the
first Puiseux exponent of $f$.}
$$
See \cite{Loeser87} for some higher dimensional generalizations.
\enddemo

The following result furnishes the basic estimates for the lc-threshold:

\proclaim{8.10 Lemma}  Let $f$ be a holomorphic function near $0\in \c^n$
and 
$D=(f=0)$.   Set $d=\mult_0f$ and let $f_d$ denote the degree $d$
homogeneous part of the Taylor series of $f$.  Let
$T_0D\deq (f_d=0)\subset \c^n$ be the  tangent cone
of $D$ and 
$\p(T_0D)\deq (f_d=0)\subset \p^{n-1}$  the projectivized tangent cone
of $D$. Then

(8.10.1) 
$\frac1{d}\leq c_0(D)\leq \frac{n}{d}$.

(8.10.2) $c_0(D)= \frac{n}{d}$ iff $(\p^{n-1}, \frac{n}{d}\p(T_0D))$ is lc.

(8.10.3) If $\p(T_0D)$ is smooth (or even lc) then  $c_0(D)=
\min\{1,\frac{n}{d}\}$.

(8.10.4) $ c_0(T_0D)\leq  c_0(D) $.
\endproclaim

\demop  Let  $p:Y\to \c^n$ 
 the blow up of $0\in {\c^n}$ with exceptional divisor $E\subset Y$. 
Then
$$
K_Y=p^*K_{\c^n}+(n-1)E,\qtq{and}  p^*D=p^{-1}_*D+dE.
$$
Thus by (8.5.1), $c_0(D)\leq (n-1+1)/d$.
In particular, 
$c_0(D)= \frac{n}{d}$ iff $(\c^n, \frac{n}{d}D)$ is lc. 
$$
K_Y+E+\frac{n}{d}p^{-1}_*D=p^*(K_{\c^n}+\frac{n}{d}D),
$$
and by (3.10.2) we see that 
$$
c_0(D)= \frac{n}{d}\quad\Leftrightarrow \quad(Y,E+\frac{n}{d}p^{-1}_*D)\qtq{is
lc.}
$$
Observe that  $E\cap p^{-1}_*D=\p(T_0D)$. 
We can apply inversion of adjunction (7.5) to see that 
$$
(Y,E+\frac{n}{d}p^{-1}_*D)\qtq{is lc}\Leftrightarrow \quad
(\p^{n-1}, \frac{n}{d}\p(T_0D))\qtq{is lc.}
$$
This shows (8.10.2) which implies (8.10.3). 

In order to see the lower bound in (8.10.1), choose  
local coordinates  such that  the Taylor series of
$f$ has the form
$x_1^d+\dots$. Consider the deformation
$f_t=t^{-d}f(tx_1,t^2x_2,\dots,t^2x_n)$.  For $t\neq 0$ the singularity of
$(f_t=0)$  is isomorphic to $(f=0)$, and for $t=0$ we get $(x_1^d=0)$.
By (8.6) we see that
$$
\frac{1}{d}=c_0(x_1^d)\leq c_0(D).
$$
To see (8.10.4)  we use  the deformation
$f_t=t^{-d}f(tx_1,tx_2,\dots,tx_n)$.  For $t\neq 0$ the singularity of
$(f_d=0)$  is isomorphic to $(f=0)$, and for $t=0$ we get $(f_k=0)$.
By (8.6) we are done. \qed
\enddemo

\demo{8.11 Remarks}  (8.11.1) It is not true in general that truncation of $f$
yields  a smaller value  for $c_0$. For instance, let $f=x^2+2xy^2+y^4$. Then
$c_0(f)=1/2$, but  $c_0(x^2+2xy^2)=3/4$.

(8.11.2)   (8.10.1) shows that $c_0(D)^{-1}$ behaves roughly as
$\mult_0(D)$. For this reason the number $c_0(D)^{-1}$
is sometimes called the Arnold multiplicity of $D$ or of $f$. 

(8.11.3) The estimate $c_0(D)\leq n/d$ holds for any $\q$-divisor $D$, even
if $D$ contains components with negative coefficients.
\enddemo

Looking at the homogeneous leading term does not give a true indication of the
subtle behaviour of the lc-threshold. In order to get better examples, we
need to look at the weighted homogeneous case. The best way to study it is by
using weighted blow ups, see \cite{Reid80,87}. In many cases weighted blow
ups can be reduced to an ordinary blow up using the following lemma, which is
a direct consequence of (3.16). 

\proclaim{8.12 Lemma} Let $p:X\to Y$ be a  finite and
dominant morphism between normal varieties. Let $\Delta_Y$ be a $\q$-divisor on
$Y$ and define $\Delta_X$ by the formula
$$
K_X+\Delta_X=p^*(K_Y+\Delta_Y),\qtq{that is,} \Delta_X=p^*\Delta_Y-K_{X/Y}.
$$
Let $D_Y$ be an effective  $\q$-divisor on $Y$ and $Z\subset Y$ a closed
subscheme. Then 
$$
c_Z(Y,\Delta_Y, D_Y)=c_{p^{-1}(Z)}(X,\Delta_X, p^*D_Y).\qed
$$
\endproclaim 

\demo{8.12.1 Remark} If $\dim X=2$ and $X$ is klt, then $X$ has quotient
singularities. Thus, by (8.12), $\Cal T_2$ can be determined by computing
$c_0(Y,D)$ where $Y$ is smooth. In general, determining $\Cal T_n$ can be
reduced to the case when $Y$ is canonical, and, assuming MMP, to the case when
$Y$ is terminal (cf. \cite{Koll\'ar94, pp.267-268}). A reduction to the smooth
case is not known.
\enddemo

We can now prove the analog of (8.10) in the weighted case:

\proclaim{8.13 Proposition}  Let $f$ be a holomorphic function near $0\in
\c^n$.  Assign rational weights $w(x_i)$ to the variables and let 
$w(f)$ be the weighted multiplicity of $f$ (= the lowest weight of the
monomials occurring in $f$).  Then
$$
c_0(f)\leq \frac{\sum w(x_i)}{w(f)}.
$$
\endproclaim

\demop   We may assume that the weights $w(x_i)$ are natural numbers. Set
$X\cong
\c^n$ with coordinates
$z_i$, $H_i=(z_i=0)$ and let 
$p:X\to
\c^n$ be given by $z_i\mapsto z_i^{w(x_i)}=x_i$. Let
$F\deq f(z_1^{w(x_1)},\dots,z_n^{w(x_n)})$ and note that $\mult_0F=w(f)$. Then
$$
K_X+\sum (1-w(x_i))H_i +c(F=0)=p^*(K_{\c^n}+c(f=0)).
$$
By (8.11.3) we obtain that if $(X,\sum (1-w(x_i))H_i +c(F=0))$ is lc
then
$$
\sum (1-w(x_i))+cw(f)\leq n,\qtq{or equivalently,} 
c\leq \frac{\sum w(x_i)}{w(f)}.\qed
$$
When do we have equality? Let $b=\sum w(x_i)/w(f)$. As in the proof of (8.10),
let $p:Y\to X$ denote the blow up of the origin with exceptional divisor $E$. 
We obtain that 
$$
\align
&(X,\sum (1-w(x_i))H_i +b(F=0))\qtq{is lc, iff}\\
&(Y,E+\sum (1-w(x_i))p^{-1}_*(H_i)+b\cdot
p^{-1}_*(F=0))\qtq{is lc.}
\endalign
$$
A slight problem is that the coefficients of the
$H_i$ are negative, and inversion of adjunction fails with negative
coefficients.   Thus we can only assert that if the leading term of $F$
defines a smooth (or just lc) hypersurface 
$(F_{w(f)}=0)\subset \p^{n-1}$, then 
$(Y,E+b\cdot p^{-1}_*(F=0))$ is lc. Subtracting the divisors $p^{-1}_*(H_i)$
only helps, thus we obtain the following:
\enddemo

\proclaim{8.14 Proposition}   Let $f$ be a holomorphic function near $0\in
\c^n$ and 
$D=(f=0)$.  Assign integral weights $w(x_i)$ to the variables and let 
$w(f)$ be the weighted multiplicity of $f$. Let
$f_w$ denote the weighted
 homogeneous leading term of the Taylor series of $f$.  
Assume that
$$
(f_w(z_1^{w(x_1)},\dots,z_n^{w(x_n)})=0)\subset \p^{n-1}\qtq{is smooth (or
lc).}
$$
Then
$$
c_0(D)=\frac{\sum w(x_i)}{w(f)}.\qed
$$
\endproclaim

\demo{8.14.1 Remark} Using a weighted blow up, it is not hard to see that
(8.14) also holds if $f$ is  semiquasihomogeneous, that
is, if $f_w$ has an isolated critical point at the origin
(\cite{AGV85,I.12.1}).
\enddemo

The following  examples give some explicit formulas. 

\demo{8.15 Example}  (8.14) shows that
$$
c_0(\sum x_i^{b_i})=\min\{1, \sum\frac{1}{b_i}\}.
$$
 Define  sets of  numbers  by
$$
\Cal F_n\deq \left\{\sum_{i=1}^n\frac{1}{b_i}\vert b_i\in \n\right\}
\cap (0,1].
$$
(8.14) shows that $\Cal F_n\subset \Cal T_n$.  $\Cal F_n$
satisfies the ascending chain condition for every $n$. 
The set of accumulation points  of $\Cal F_n$ is precisely
$\Cal F_{n-1}$. 
This shows that the sets $\Cal T_n$ have plenty of
accumulation points in the interval $[0,1]$.
\enddemo

\demo{8.16 Example} An equivalent formulation of the 
ascending chain condition is that each subset of the set  has a maximal
element. Thus, for instance, $\Cal T_n\cap [0,1)$ has a  maximal element.
Let us denote it by $1-\delta'(n)$ (cf. \cite{Koll\'ar94, 5.3.3}). 
It is known that $\delta'(1)=1/2$, $\delta'(2)=1/6$ and $\delta'(3)=1/42$
\cite{Koll\'ar94, 5.4}.
More generally, define a sequence $a_i$ by the recursive formula
$$
a_1=2, a_{k+1}=a_1a_2\cdots a_k+1.
$$
From (8.14) we obtain that 
$$
c_0(x_1^{a_1}+\dots+x_n^{a_n})=1-\frac1{a_{n+1}-1}.
$$
It is possible that $\delta'(n)=1/(a_{n+1}-1)$ for every $n$.
It is known  that the maximal element of 
$\Cal F_n\cap [0,1)$ is $1-1/(a_{n+1}-1)$ \cite{Soundararajan95}.
\enddemo

\demo{8.17 Example} An analysis of the proof of (8.14) shows that
$$
c_0((\prod x_i^{a_i})(\sum x_i^{b_i}))=
\min\left\{\frac{ \sum_i 1/b_i}{1+\sum_i a_i/b_i},
\frac1{a_1},\dots,\frac1{a_n}\right\}.
$$
 Define  sets of  numbers  by
$$
\Cal G_n\deq \left\{\frac{ \sum_{i=1}^n 1/b_i}{1+\sum_{i=1}^n a_i/b_i}
\vert a_i,b_i\in\n\right\}.
$$
It is quite remarkable that $\Cal G_n$  does not satisfy the
ascending chain condition. For example, fix the numbers $a_i$ and
$b_1,\dots,b_{n-1}$ and let $b_n\to \infty$.  
Then
$$
\lim_{b_n\to \infty} \frac{ \sum_i^n 1/b_i}{1+\sum_i^n a_i/b_i}=
\frac{ \sum_i^{n-1} 1/b_i}{1+\sum_i^{n-1} a_i/b_i},
$$
and the sequence is increasing iff
$$
\frac{ \sum_i^{n-1} 1/b_i}{1+\sum_i^{n-1} a_i/b_i}>\frac1{a_n}.
$$
It is precisely in this case that the lc-threshold is computed
by $\min\{1/a_i\}$ and not by the main part  
$(\sum_i 1/b_i)/(1+\sum_i a_i/b_i)$.  

It is not hard to see that $\Cal G_n\cap \Cal T_n$ satisfies the 
ascending chain condition.
\enddemo

The last question that we   consider in this section is the following. 
Let $f\in \c[[x_1,\dots,x_n]]$ be a power series. How well can one approximate
$c_0(f)$ by computing
$c_0$ of some polynomials? More specifically, let $f_{\leq d}$ be the degree
$\leq d$ part of $f$. What can one say about the difference
$c_0(f)-c_0(f_{\leq d})$?

\demo{8.18 Example} (8.18.1) Assume that $f$ defines an isolated singularity.
Then  $f$ and $f_{\leq d}$ differ only by a coordinate
change for
$d\gg 1$, thus $c_0(f)=c_0(f_{\leq d})$.

(8.18.2)  Let $f=(y+x^2+x^3+\dots)^2\in \c[[x,y]]$. Then $c_0(f)=1/2$.
Furthermore,
$$
f_{\leq d}=(y+x^2+x^3+\dots+x^{d-1})^2-\rdown{(d-1)/2}x^{d+1}-\dots.
$$
Change variables to $z=y+x^2+x^3+\dots+x^{d-1}$. We get that
$$
f_{\leq d}=z^2-\rdown{(d-1)/2}x^{d+1}-\dots.
$$
Setting $w(z)=1/2$ and $w(x)=1/(d+1)$ shows that $c_0(f_{\leq
d})=1/2+1/(d+1)$. 
\enddemo

Thus the best we can hope is that $c_0(f_{\leq d})$ converges to $c_0(f)$ in 
a uniform way. Theorem (8.20) is a much more precise result. Together with
(8.10.1) it implies the following estimate:

\proclaim{8.19 Proposition} Let $f\in \c[[x_1,\dots,x_n]]$ be a power series
and let 
$f_{\leq d}$ denote the degree
$\leq d$ part of $f$.  Then
$$
|c_0(f)-c_0(f_{\leq d})| \leq \frac{n}{d+1}.\qed
$$
\endproclaim

\proclaim{8.20 Theorem} \cite{Demailly-Koll\'ar96} Let $f,g$ be polynomials
(or power series) in $n$ variables. Then
$$
c_0(f+g)\leq c_0(f)+c_0(g)\qtq{and} c_0(fg)\leq \min\{c_0(f),c_0(g)\}.
$$
\endproclaim
$$
*****
$$
\demop  The estimate for $c_0(fg)$ is clear. The only surprising part is that
it cannot be sharpened; for instance $c_0(x^ay^b)=\min\{1/a,1/b\}$. 

The proof of the additive part has two steps. The first is the computation of
the lc-threshold for direct sums of  functions. The formula is proved in
\cite{AGV84,vol.II. sec.13.3.5} for isolated singularities.  Since
the proof does not seem to generalize to the nonisolated case, we give an
alternative argument in a more general setting (8.21).

The second step uses inversion of adjunction to go from  
$f(x_1,\dots,x_n)+g(y_1,\dots,y_n)$ to 
$f(x_1,\dots,x_n)+g(x_1,\dots,x_n)$.  Let $B_X\subset
\c^n$ be  a small ball with coordinates
$x_i$ and
$B_Y$ a small ball with coordinates $y_i$. 
Let
$$
\align
&G=(f(x_1,\dots,x_n)+g(x_1,\dots,x_n)=0)\subset B_X, \qtq{and}\\
&D=(f(x_1,\dots,x_n)+g(y_1,\dots,y_n)=0) \subset B_X\times B_Y.
\endalign
$$
Set $F_i=(x_i=y_i)$.  Applying (7.5) $n$-times to $cD+F_1+\dots+F_n$ we
obtain that
$$
(B_X,cG) \qtq{is lc} 
\Leftrightarrow\quad
(B_X\times B_Y, cD+F_1+\dots+F_n) \qtq{is lc.}
$$
The latter clearly implies that $(B_X\times B_Y, cD)$ is lc. This shows that
$c_0(B_X,G)\leq c_0(B_X\times B_Y,D)$.   By (8.21),
$c_0(B_X\times B_Y,D)\leq c_0(B_X,f)+c_0(B_X,g)$, which completes the
proof.\qed\enddemo

\proclaim{8.21 Proposition}  Let $(X_1,\Delta_1)$ and   $(X_2,\Delta_2)$ be
lc pairs with marked points $x_i\in X_i$. Let 
$$
X=X_1\times X_2,\qtq{and} \Delta=\Delta_1\times X_2+X_1\times \Delta_2
$$
be their product; $x=(x_1,x_2)\in X$. Let $f_i$ be a regular function on
$X_i$, 
$D_i=(f_i=0)$ and $D=(f_1+f_2=0)$.  Then
$$
c_x(X,\Delta, D)=\min\{1, c_{x_1}(X_1,\Delta_1, D_1)+
c_{x_2}(X_2,\Delta_2, D_2)\}.
$$
\endproclaim

\demop  In order to simplify notation, we pretend that $X_i$ is  local with
closed point $x_i$.  Let  $p_i:Y_i\to X_i$ be log resolutions and write
$$
K_{Y_i}=p_i^*(K_{X_i}+\Delta_i)+\sum_j a_{ij}E_{ij},
\qtq{and} p_i^*D_i=\sum_j b_{ij}E_{ij}.
$$
Set $Y=Y_1\times Y_2$ and let $p:Y\to X$ be the product morphism. $p$ is a
resolution of $X$ but it is rarely a log resolution. 
Set $E'_{1j}=E_{1j}\times Y_2$ and $E'_{2j}=Y_1\times E_{2j}$. Then
$$
K_{Y}=p^*(K_X+\Delta)+\sum_{ij} a_{ij}E'_{ij}.
$$
The problem is that $p^*D$ is not a sum of the divisors $E'_{ij}$.

To study the situation, choose indices $j,k$ and points $y_1\in E_{1j}$ and
$y_2\in E_{2k}$ such that  $\sum_jE_{ij}$ is smooth at $y_i$.
Let $v_1$ (resp. $v_2$) be a local defining equation of  $E_{1j}$ (resp. 
$E_{2k}$).  By suitable choice of the $v_i$ we may assume that locally near
$y$
$$
p_1^*f_1=v_1^{b_{1j}},\qtq{and} p_2^*f_2=v_2^{b_{2k}}.
$$
In a neighborhood of $y\in Y$ there are two exceptional divisors
$E'_{1j}=(v_1=0)$ and $E'_{2k}=(v_2=0)$, and 
$$
p^*f=v_1^{b_{1j}}+v_2^{b_{2k}},\qtq{locally near $y$.}
$$
Set $F_{jk}=(v_1^{b_{1j}}+v_2^{b_{2k}}=0)$. 
If $(X,\Delta+cD)$ is lc, then by (3.10) we obtain that
$$
(Y, -a_{1j}E'_{1j}-a_{2k}E'_{2k}+cF_{jk}) \qtq{is lc near $y$.}
\tag 8.21.1
$$
  (8.11.3) shows that (8.21.1) is equivalent to
$$
c-\frac{a_{1j}}{b_{1j}}-\frac{a_{2k}}{b_{2k}}\leq 
\frac{1}{b_{1j}}+\frac{1}{b_{2k}},
\qtq{that is,}
c\leq \frac{a_{1j}+1}{b_{1j}}+\frac{a_{2k}+1}{b_{2k}}.
$$
By (8.5), 
$$
c_{x_1}(X_1,\Delta_1, D_1)=\min_j\frac{a_{1j}+1}{b_{1j}}
\qtq{and} 
c_{x_2}(X_2,\Delta_2, D_2)=\min_k\frac{a_{2k}+1}{b_{2k}}.
$$
Choose $j$ and $k$ such that  they achieve the minima. Then we obtain that
$$
c_x(X,\Delta, D)\leq  c_{x_1}(X_1,\Delta_1, D_1)+
c_{x_2}(X_2,\Delta_2, D_2).
$$
We have not proved equality, since the above procedure controls only those
divisors $F$ of $K(X)$ such that $\cent_Y(F)\supset E'_{1j}\cap E'_{2k}$
for some $j,k$. 
There are two ways to go to the general case.

First, we can do a similar computation at any point of $Y$.  I found this
somewhat cumbersome.

Second, we could try to show that the above computation accounts for all
divisors
$F$ of $K(X)$, if we vary $Y_1$ and $Y_2$. Indeed, for suitable choice of
$Y_i$ we may assume that  the image of 
$F$ on $Y_i$ contains a divisor $E_{1j}$ (resp. $E_{2k}$)
(3.17). 
This means that  $\cent_Y(F)\supset E'_{1j}\cap E'_{2k}$, as
required.\qed\enddemo

\demo{8.21.2 Exercise} Use (8.21) to show that every element of $\Cal
T_n-\{1\}$ is an accumulation point of $\Cal T_{n+1}$. 

It is possible that the set of accumulation points  of $\Cal T_n$ is precisely
$\Cal T_{n-1}$.  In the toric case this was proved by \cite{Borisov95}. 
\enddemo

\head 9. The Log Canonical Threshold and the Complex Singular Index
\endhead

In this section we compare the lc-threshold  and 
the complex singular index of an isolated singularity. The notion of 
complex singular index was introduced by Arnold, using 
the asymptotic behaviour of
certain integrals over vanishing cycles. 
See  \cite{AGV85, II.Chap.13} for the motivation and for basic results.

The classical  case is the following:

\demo{9.1 Definition}  Let $f:(0,\c^{n+1})\to (0,\c)$  be a holomorphic
function in the neighborhood of the origin. Assume  that 
$f$ has an isolated critical point at the origin. Set $D\deq(f=0)$. 
Let $B\subset \c^{n+1}$ be a small ball around the origin and $\Delta\subset
\c$ an even smaller disc around the origin.  Set $X=B\cap f^{-1}(\Delta)$.
From now on, we restrict
$f$ to
$f:X\to\Delta$. 

By \cite{Milnor68}, the only interesting homology of $X_t\deq f^{-1}(t)$ 
for $t\neq 0$ is in dimension $n$. The corresponding cycles are called the {\it
vanishing cycles}.

If $\sigma $ is a section of $\omega_{X/\Delta}$ then $\sigma$ restricts to a
holomorphic $n$-form on each $X_t$. Thus if $\delta(t)$ is an $n$-cycle in
$X_t$, then we can form the integral
$$
\int_{\delta(t)} \sigma,\tag 9.1.1
$$
which depends only on the homology class of $\delta(t)$.

Let $x_1,\dots,x_{n+1}$ be local coordinates on $\c^{n+1}$. A local
generator of $\omega_{X/\Delta}$ can be written down explicitly. Up-to a
sign, it is
$$
\frac{dx_1\wedge\dots\wedge dx_{n+1}}{df}\deq\pm
\frac{dx_1\wedge\dots\wedge \widehat{dx_i}\wedge\dots\wedge dx_{n+1}}
{\partial f/\partial x_i}.\tag 9.1.2
$$
\enddemo

\demo{9.2 General case} More generally, we can consider an arbitrary (normal)
complex space $X$ and  any morphism 
$f:X\to\Delta$ such that $f$ is smooth over $\Delta^*$ and $X-f^{-1}(0)\to
\Delta^*$ is a locally trivial topological fiber bundle over $\Delta^*$. (The
latter can usually be assumed by shrinking $\Delta$.)

Let $\sigma $ be a section of $\omega_{X/\Delta}$.  $\sigma$ restricts to a
section of $\omega_{X_t}$ for each $t$. If $X_t$ is smooth,  then 
$\int_{\delta(t)} \sigma$ makes sense as above.
\enddemo

The following basic theorem describes the asymptotic behaviour of the
integrals $\int_{\delta(t)} \sigma$ for small values of $t$.  It can be
approached from many different points of view. See,  for instance,
\cite{AGV85, II.10.2} for a discussion and several references.

\proclaim{9.3 Theorem}   Notation and assumptions
as above.  Let $t\mapsto \delta(t)\in H_n(X_t,\q)$ be a continuous (multiple
valued) section. 
There is an asymptotic
expansion (as $t\to 0$)
$$
\int_{\delta(t)} \sigma =
\sum_{\alpha\in\q, k\in\n} a(\sigma, \delta, \alpha, k) t^{\alpha} (\log t)^k,
$$
where the $a(\sigma, \delta, \alpha, k)$ are constants.\qed
\endproclaim

\proclaim{9.3.1 Complement} One can get rather precise information about the
possible values of $\alpha$ and $k$. The following are some of these:

(9.3.1.1) There is  a lower bound for the values of $\alpha$, depending only
on
$f:X\to\Delta$. (This also follows from (9.5).)

(9.3.1.2) The values of $\alpha \mod 1$ can be described in terms of the
eigenvalues of the monodromy.

(9.3.1.3)  There is an upper bound for the values of $k$ depending on the size
of the Jordan blocks of the monodromy.

(9.3.1.4) If  $X_0$ is a normal crossing divisor, then 
$a(\sigma, \delta, \alpha, k)=0$ for $\alpha< 0$.
\endproclaim

\demo{9.4 Definition}    Let $X$ be a normal complex space and $f:X\to\Delta$
 a morphism. Set $D=f^{-1}(0)$.  Let $x\in D$ be a closed point and
assume that 
$f$ is smooth on $X-x$. Assume also that $\omega_X$ is locally  free.

The {\it complex singular index}, denoted by  $\beta_{\c}(f)$ or by
$\beta_{\c}(X,D)$ is defined by the formula
$$
\beta_{\c}(X,D)=
\beta_{\c}(f)\deq 1+\inf\{\alpha\vert \exists \sigma, \delta, k \text{ such
that } a(\sigma, \delta, \alpha, k)\neq 0\}.
$$

The definition gives $\beta_{\c}(f)=\infty$ if there are no vanishing cycles
at all. If $X$ is smooth this happens only when $D$ is also smooth. 

In most cases the asymptotic expansion involves negative powers of $t$, thus
the complex singular index measures the maximum rate of divergence of 
the above integrals as $t\to 0$.

The terminology is taken from \cite{Steenbrink85}. In \cite{Varchenko82,
p.477} this is called the complex singular exponent and in
\cite{AGV85,II.13.1.5} the complex oscillation index.
\enddemo

The following theorem relates the  lc-threshold to the complex singular
index. For $X$ smooth it was proved by \cite{Varchenko82, \S 4}. The proof
also works in a more general setting. Further generalizations are pointed
out in (9.7).

\proclaim{9.5 Theorem} \cite{Varchenko82, \S 4}  Notation and assumptions
as in (9.4). 
Then 
$$
c_x(f)= \min\{1,\beta_{\c}(f)\}.
$$
\endproclaim

\demop Let $\pi_n:\Delta_n\to \Delta$ be the morphism $t_n\mapsto t=t_n^n$, and
consider the fiber product diagram
$$
\CD
X_n@>\pi_n>> X\\
@Vf_nVV @VVfV\\
\Delta_n@>\pi_n>> \Delta.
\endCD 
\tag 9.5.1
$$
We can identify the central fiber of $f_n$ with the central fiber $D$ of
$f$. For suitable choice of $n$ we may assume that $f_n:X_n\to \Delta_n$ has
a semistable resolution. That is, there is a proper birational morphism
$g_n:Y_n\to X_n$ such that $Y_n$ is smooth and
$$
 g_n^*D=D_0+\sum_{i>0}D_i
$$
 is a reduced divisor with normal
crossings only, where $D_0$ is the birational transform of $D$.

Using  (3.16) we obtain that
$$
c_x(f_n)=1-n(1-c_x(f)),
$$
where we use the general definition of the lc-threshold (8.4). 
Write
$$
\align
K_{Y_n}&=g_n^*K_{X_n}+\sum_{i>0} a_iD_i,\qtq{hence also}\\
K_{Y_n/\Delta_n}&=g_n^*K_{X_n/\Delta_n}+\sum_{i>0} a_iD_i.
\endalign
$$
 Set
$a_0=0$, 
$c=\min_{i\geq 0}\{a_i\}$  and  note that $c_x(f_n)= c+1$. 
The crucial formula is
$$
K_{Y_n/\Delta_n}=g_n^*(K_{X_n/\Delta_n}+cD)+\sum_{i\geq 0} (a_i-c)D_i,
\tag 9.5.2
$$
where $a_i-c\geq 0$ for every $i\geq 0$. 
Let $\sigma$ be any section of $\o(K_{X/\Delta})$. 
$K_{X_n/\Delta_n}=\pi_n^*K_{X/\Delta}$, thus $\pi_n^*\sigma$ is a  section of
$\o(K_{X_n/\Delta_n})$. Therefore $t_n^{-c}\pi_n^*\sigma$ 
is a  section of
$\o(K_{X_n/\Delta_n}+cD)$, hence by (9.1) it corresponds to a holomorphic
section $\sigma_n'$ of $\o(K_{Y_n/\Delta_n})$.

Up to an $n^{th}$-root of unity, $t_n^{-c}\pi_n^*\sigma=
\pi_n^*(t^{1-c_x(f)}\sigma)$, and therefore
$$
t^{1-c_x(f)}\int_{\delta(t)}\sigma=\int_{\delta(t_n)}\sigma_n',
\qtq{for $t\neq 0$, where $t_n=t^{1/n}$.}
$$
  $\sigma_n'$ is a holomorphic
section of $\o(K_{Y_n/\Delta_n})$, and so, by (9.3.1.4), 
$$
\int_{\delta(t_n)}\sigma_n'\tag 9.5.3
$$
grows at most logarithmically as $t_n\to 0$. This shows that
$c_x(f)\leq \beta_{\c}(f)$, hence also $c_x(f)\leq
\min\{1,\beta_{\c}(f)\}$ and equality holds if $c_x(f)=1$.

In order to see the equality in the remaining cases, 
we have to find $\delta(t_n)$ such that the integral (9.5.3) grows as a
nonzero constant times a power of  $\log t_n$.  Thus assume that 
$c_x(f)<1$, which is equivalent to $c<0$.  Let $E\deq
D_j,\ j>0$ be an irreducible component  such that $a_j=c$. Such a component
exists since $c<0$ and it is proper since $D$ has isolated singularities. 
Set  
$E^0=E-\cup_{i\neq j}D_i$; this is  an open set of $E$.  

$\sigma'_n$ restricts to a holomorphic section of the dualizing
sheaf of $\cup_{i\geq 0}D_i$, thus
 $\sigma'_n|E$ is
a holomorphic $n$-form on
$E$ with at worst simple poles along $E-E^0$ (for top degree forms this
is the same as having logarithmic poles).  By
\cite{Deligne71}, closed forms with logarithmic poles at infinity compute the
cohomology of a smooth  variety,  thus there is an
$n$-cycle
$Z\subset E^0$ such that $\int_Z (\sigma'_n|E)\neq 0$. 

$f_n$ is a locally trivial fibration    near $E^0$, thus $Z$ can be extended
to an $n$-cycle
$\delta(t_n)$ for small values of $t_n$. 
(We even get a monodromy invariant cycle, but this is not important for now.)
Let $\sigma$ be a local generator of
$\o(K_{X/\Delta})$.  By construction,
$$
\lim_{t\to 0}t^{1-c_x(f)}\int_{\delta(t)}\sigma=\lim_{t_n\to
0}\int_{\delta(t_n)}\sigma_n'=\int_Z (\sigma'_n|E)\neq 0,
$$
thus the asymptotic expansion of 
$$
\int_{\delta(t)}\sigma
$$
does contain a nonzero term $const\cdot t^{c_x(f)-1}$. \qed\enddemo
$$
*****
$$
\demo{9.6 Generalizations}  From the point of view of the lc-threshold, the
assumptions that $D$ has an isolated singularity and that $f$ is smooth over
$\Delta^*$ are rather restrictive. Some of these conditions can be weakened.

If $D$ does not have  an isolated singularity, then it is not clear what
exactly happens. For instance, assume that $f(x,y)$ defines an isolated
singularity. Viewed as a map $f_2: \c^2\to \c$ the  1-dimensional
homology of the nearby fibers gives the vanishing cycles. If we view $f$ as a
morphism $f_3:\c^3\to \c$, then the central fiber has a  nonisolated
singularity.  Furthermore, $f_3^{-1}(t)\cong f_2^{-1}(t)\times \c$, thus all
the interesting homology is in 1-dimension and we cannot integrate a 2-form.
I do not know how to overcome this problem, except in some special cases.

If $f$ is not smooth over
$\Delta^*$, we can proceed as follows.

  Let $p:X'\to X$
be a
 resolution of singularities and $f':X'\to
\Delta$ the induced morphism. Assume that 
$f'$ is smooth over $\Delta^*$. Let  $\delta(t)\in H_n(X'_t)$ be a 
continuous (multiple
valued) section. Assume furthermore that $X_t$
has canonical singularities for $t\neq 0$.

Under these assumptions, the integral  $\int_{\delta(t)} p^*\sigma$
makes sense and it behaves like the integral (9.1.1).

 The proof of (9.5) shows that the result also holds more
generally:
\enddemo

\proclaim{9.7 Theorem} Let $X$ be a normal analytic space  and $f:X\to
\Delta$ a morphism. Assume that 

(9.7.1) $\omega_X$ is locally free;

(9.7.2) $f^{-1}(0)$ has rational singularities except at a  single
point $x\in f^{-1}(0)$;

(9.7.3) $f^{-1}(t)$ has rational singularities for $t\neq 0$. 

\noindent Then 
$$
c_x(f)= \min\{1,\beta_{\c}(f)\}.\qed
$$
\endproclaim

The log canonical threshold is also related to the constants of
quasiadjunction intruduced in \cite{Libgober83} and further studied in
\cite{Loeser-Vaqui\'e90}.

\proclaim{9.8 Proposition}  Let $f(x_1,\dots,x_n)$ define a singularity at
the origin. For every $m$ let $\psi(m)$ be the smallest integer such that
$y^{\psi(m)}$ is contained in the adjoint ideal of the hypersurface
$X_m\deq (y^m=f)$. Then $\psi(m)=\rdown{m(c_0(f)+1)}$. 
\endproclaim

\demop  $y^{1-m}dx_1\wedge\dots\wedge dx_n$ is a local generator of
$\omega_{X_m}$, thus $y^{\psi(m)}$ is contained in the adjoint ideal
iff $y^{\psi(m)+1-m}dx_1\wedge\dots\wedge dx_n$ is $L^2$. Pushing down to
$\c^n$, this is equivalent to $|f|^{(\psi(m)+1-m)/m}$ being $L^2$. This
happens precisely when $(\psi(m)+1-m)/m>c_0(f)$, which is eqivalent to
$\psi(m)=\rdown{m(c_0(f)+1)}$. \qed\enddemo

\head 10. The Log Canonical Threshold and the Bernstein-Sato Polynomial
\endhead

The aim of this section is to compare the log canonical threshold of a
function $f$ to the  Bernstein-Sato polynomial of $f$. The basic definitions
are given bellow.

\proclaim{10.1 Theorem} \cite{Bernstein71; Bj\"ork79} Let $f=f(z_1,\dots,z_n)$
be a polynomial (resp. a convergent power series)  and   $s$   a
variable. There is a nonzero polynomial $b(s)\in \c[s]$ and  a linear
differential operator 
$$
P=\sum_{I,j}f_{I,j}s^j\frac{\partial^I}{\partial z^I},
$$
whose coefficients $f_{I,j}$ are  polynomials (resp.  convergent power
series) such that
$$
b(s)f^s=Pf^{s+1}.\qed
\tag 10.1.1
$$
\endproclaim

\demo{10.1.2 Remark}  It is easiest   to interpret (10.1.1)  as a formal
equality, where we do not assign any meaning to the powers $f^s$, just handle
them as symbols with the usual roles of differentiation assumed. If the
powers have  a well defined meaning as functions (for instance,  $f$ is
everywhere nonnegative on $\r^n$) then the formal equality becomes an actual
equality of functions.
\enddemo

\demo{10.2 Definition}  All the polynomials satisfying (10.1) form an ideal in 
$\c[s]$. The unique generator of this ideal with leading coefficient 1 is
called the {\it Bernstein-Sato polynomial} of $f$. It is denoted by $b_f(s)$.

In singularity theory, many people  use the defining equation
$b(s)f^{s-1}=Pf^{s}$; this corresponds to
the substitution $s\deq s+1$ in the polynomial $b_f$.
\enddemo

\demo{10.3 Remark} The polynomial $b_f$ is a very interesting invariant of the
singularity $(f=0)$. It can be connected with with other types of invariants
in many different ways, see, for instance, \cite{Malgrange75; Loeser87} and the
references there.
\enddemo

\demo{10.4 Definition} 
Setting $s=-1$, (10.1.1) becomes $b(-1)f^{-1}=\sum_j f_{0,j}$, which implies
that  $b(-1)=0$. Thus $b_f(s)=(s+1)\tilde b_f(s)$.  $\tilde b_f(s)$ is called
the {\it reduced Bernstein-Sato polynomial} of $f$. 
\enddemo

\demo{10.5 Examples}  It is not easy to compute $b_f$ and $P$ in concrete
examples. 

(10.5.1) For quadratic forms the answer is rather obvious. Set $Q(z)=\sum
z_i^2$, then
$$
(s+1)(s+\frac{n+1}{2}) Q(z)^s=
\frac1{4}\left(\sum \frac{\partial^2}{\partial z_i^2}\right) Q(z)^{s+1}.
$$

(10.5.2) Already
the case of cusps is nontrivial:
$$
(s+1)(s+\frac{5}{6})(s+\frac{7}{6})(x^2+y^3)^s=
\left(\frac1{27}\frac{\partial^3}{\partial y^3}+
\frac{y}{6}\frac{\partial^3}{\partial x^2\partial y}+
\frac{x}{8}\frac{\partial^3}{\partial x^3}
\right)
(x^2+y^3)^{s+1}.
$$

(10.5.3) Assume that $f$ defines an isolated singularity at the origin, and if
we set
$wt(z_i)=a_i$ then $f$ is weighted homogeneous of degree 1.  By
\cite{Yano78}, 
$$
\prod_i\frac{t^{a_i}-t}{1-t^{a_i}}=\sum_{\alpha\in \q} q_{\alpha}t^{\alpha}
\qtq{is  a finite sum, and}
\tilde b_f(s)=\prod_{\alpha:q_{\alpha}\neq 0}(s+\alpha).
$$

(10.5.4) Let $f=\sum z_i^{m_i}$.  One can easily compute using (10.5.3) that
$$
\text{largest root of $\tilde b_f$}=-\sum \frac1{m_i}.
$$
\enddemo

The following   observation relates the roots of 
Bernstein-Sato polynomials to the lc-threshold of $f$:

\proclaim{10.6 Theorem} Let $f=f(z_1,\dots,z_n)$ be
a polynomial or a convergent power series. Then 
$$
\text{largest root of $b_f(s)$}= -\text{(lc-threshold of $f$)}
$$
\endproclaim

\demop By definition, we have $b(s)f^{s}=Pf^{s+1}$ and conjugating it we
obtain
$\bar b(s)\bar f^{s}=\bar P\bar f^{s+1}$. Two differential operators do not
commute in general, but a holomorphic operator always commutes with an
antiholomorphic one, and, moreover
$$
Pf^{s+1}\cdot \bar P\bar f^{s+1}=(P\bar P)|f^2|^{s+1}.
$$
Let $\phi$ be any $C^{\infty}$ function supported in a small
neighborhood of the origin. Then we have the equality
$$
|b(s)|^2\int |f^2|^s\phi dm(z)=
\int |f^2|^{s+1}(P\bar P)\phi d m(z),
$$
where $dm(z)$ is the Lebesgue measure.

Let $c(f)$ be the lc-threshold of $f$.
As long as $s>-c(f)$, both sides are well defined and
finite.  If $\phi$ is positive and nonzero at the origin, then
$$
\lim_{s\to -c(f)^+}\int |f^2|^s\phi dm(z)=+\infty,
\qtq{and}
\int |f^2|^{-c(f)+1} (P\bar P)\phi d m(z)<\infty.
$$
This shows that $-c(f)$ is a root of $b_f$.

Assume that $t>-c(f)$ is a root of $b_f$. We obtain that
$$
\int |f^2|^{t+1} (P\bar P)\phi d m(z)=0\qtq{for every $\phi$. }
$$
This is a rather rare accident, but cannot be excluded
without knowing something about $P$.

The actual proof that $-c(f)$ is the largest root is unfortunately rather
complicated. It follows from the next result of \cite{Lichtin89, Thm. 5}, which
in turn is a modification of the arguments in \cite{Kashiwara76}:

\proclaim{10.7 Theorem} Set $D=(f=0)\subset \c^n$ and let $\pi:Y\to \c^n$ be a 
log resolution with exceptional divisors $D_i: i>0$. Set $D_0\deq
\pi^{-1}_*D$. Write
$$
K_Y=\sum_{i>0} a_iD_i,\qtq{and} \pi^*D=\sum_{i\geq 0} d_iD_i.
$$
Then every root of $b_f$ is of the form
$$
-\frac{a_i+e}{d_i} \qtq{for some $i\geq 0$ and $e\in \n$.}\qed
$$
\endproclaim

By (8.5) we know that $c(f)=\min_i \{(a_i+1)/d_i\}$, which shows that
$b_f$ does not have any root   bigger than $-c(f)$.\qed\enddemo

\demo{10.8 Remark} If $f=0$ defines  a rational singularity, then the largest
root of $b_f$ is the trivial root $-1$ and the lc-threshold is $1$, as it
should be.

In this case it is natural to try to connect the largest nontrivial root
of $b_f$ with some geometric data coming from the resolution. The following
may appear a rather natural candidate:
$$
-\text{(largest  root
of $\tilde b_f$)}\overset{?}\to{=}\inf_{i\: \pi(D_i)\subset \sing D}
\frac{a_i+1}{d_i}.
$$
Unfortunately, the right hand side depends on the resolution chosen.

More generally, it would be  of  interest to understand which exceptional
divisors give roots of  $b_f$ in (10.7). 
\enddemo

\head
11. Rational and Canonical Singularities
\endhead

The aim of this section to prove that canonical singularities are rational.
This result was proved by 
\cite{Elkik81; Flenner81}.
The essential part of these proofs was generalized in \cite{Fujita85} and 
treated   systematically in \cite{KaMaMa87, 1-3}.  The treatment given here  
uses duality theory only for CM schemes, and this simplification makes the
proofs even a little shorter. 

\proclaim{11.1 Theorem} Let $X$ be a normal variety over a field
of characteristic zero. 

(11.1.1) Assume that $\omega_X$  is locally free. 
  Then
$X$ has rational singularities iff $X$ has canonical singularities.

(11.1.2) Assume that $(X,D)$ is a klt pair.
  Then
$X$ has rational singularities.
\endproclaim

\demo{11.1.3 Remark} 
If $\omega_X$ is not locally free, then rational and klt are no longer
equivalent. 
Most rational singularities are not klt and not even log canonical. 
For instance, a normal surface singularity is klt iff it is a
quotient singularity (3.6), but there are many rational surface singularities
which are not quotient.
\enddemo

At the end (11.15) we present a result about deformations of rational
singularities. This is a generalization of \cite{Elkik78}. 
$$
*****
$$

 As the first step of the proof,  recall the
Leray spectral sequence for local cohomology and some of its immediate
consequences:

\proclaim{11.2 Theorem} Let $f:Y\to X$ be a proper morphism, $x\in X$ a closed
point, $F=f^{-1}(x)$ and $G$ a sheaf on $Y$. 

(11.2.1) There is a Leray spectral sequence
$E_2^{ij}=H^i_x(X,R^jf_*G)\Rightarrow H^{i+j}_F(Y,G)$.

(11.2.2) The spectral sequence gives  an injection $H^1_x(X,f_*G)\hookrightarrow
H^1_F(Y,G)$.

(11.2.3) If $R^if_*G=0$ for $i>0$ then $H^j_F(Y,G)=H^j_x(X,f_*G)$ for every $j$.

(11.2.4) If $\supp R^if_*G\subset \{x\}$ for $i\geq 0$ then
$H^j_F(Y,G)=R^jf_*G$ for every $j$.

(11.2.5) If $\supp R^if_*G\subset \{x\}$ for $1\leq i<k$ and
$H^i_F(Y,G)=0$ for $i\leq k$ then
$R^jf_*G=H^{j+1}_x(X,f_*G)$ for   $j=1,\dots,k-1$.
\endproclaim

\demop It is clear that $H^0_x(X,f_*G)= H^0_F(Y,G)$.
This gives a spectral sequence between the derived functors. The
construction is the same as for  the ordinary  Leray spectral sequence
(see e.g. \cite{Griffiths-Harris78, p.462}). 

Looking at the beginning of the spectral sequence gives (11.2.2).
Under the assumptions (11.2.3) or (11.2.4) the spectral sequence degenerates
at the $E_2$ term since all the nonzero $E_2^{ij}$ are in one row or
column.

Finally assume (11.2.5). Then the only nonzero $E_2^{ij}$ for $0\leq
i,j\leq k$ are those with $ij=0$. Thus for every $j<k$ there is only one possible
nonzero differential $d:R^jf_*G\to H^{j+1}_x(X,f_*G)$ which must be an
isomorphism since
$H^i_F(Y,G)=0$ for $i\leq k$.\qed\enddemo

\demo{11.3 Definition} Let $X$ be a scheme of pure dimension $n$ and $G$ a sheaf
on $X$. We say that $G$ is {\it CM}   (which is an abbreviation for {\it
Cohen--Macaulay})  if  it satisfies the following equivalent conditions
(cf. \cite{Hartshorne77, Exercise III.3.4}):

(11.3.1) for
every point
$x\in X$,
$\depth_xG=\codim (x,X)$, 

(11.3.2) $H^i_x(X,G)=0$ for
every 
$x\in X$ and 
$i<\codim (x,X)$.

We say that $X$ is CM if $\o_X$ is CM.
\enddemo

Basic properties of CM sheaves are recalled in the next  lemma.

\proclaim{11.4 Lemma}  (11.4.1)  Let $X$ be a regular
scheme and
$G$ a coherent sheaf. Then $G$ is CM iff it is locally free.

(11.4.2)  Let $f:X\to Y $ be a finite morphism of  schemes of pure
dimension
$n$ and
$G$ a coherent sheaf on $X$. Then $G$ is CM iff $f_*G$ is CM.
\endproclaim

\demop   The first part is proved in \cite{Matsumura86, 19.1}.
The second assertion follows from (11.2.3).\qed\enddemo

\proclaim{11.5 Proposition} Let $X$ be an $S_2$  scheme of pure dimension $n$ and
assume that $\omega_X$ exists. Then $\o_X$ is CM iff $\omega_X$ is CM. 

More generally, if $G$ is an $S_2$ sheaf then $G$ is CM iff $\hom
(G,\omega_X)$ is CM.
\endproclaim

\demop We may clearly suppose that $X$ is affine. Assume first that
there is a finite morphism $f:X\to Z$ onto a regular scheme $Z$ of dimension
$n$. 

Then $f_*\omega_X=\hom(f_*\o_X,\omega_Z)$ and $\omega_Z$ is a line bundle. 
Since $X$ is $S_2$, $f_*\o_X$ is reflexive hence 
$f_*\omega_X$ and $f_*\o_X$ are duals (up to a twist by a line bundle). 
Thus $\omega_X$ is CM iff $f_*\omega_X$ is locally free iff  $f_*\o_X$
is locally free iff $\o_X$ is CM.

The   general case is proved the same way since 
$f_*\hom (G,\omega_X)=\hom(f_*G,\omega_Z)$
(cf. \cite{Hartshorne77, Exercise III.6.10}).

By Noether normalization $f$ always exists  if $X$ is of
finite type. 
$f$ also exists if $X$ is the spectrum of a complete local ring.
The general case can be reduced to the latter by noting that a module over a
local ring is CM iff its completion is CM over the completion of the local
ring (this is a very special case of \cite{Matsumura86, 23.3}). \qed\enddemo

The next proposition is a collection of some duality statements. They are all
special cases of the general duality theorem, but they can also be derived from
ordinary duality easily.

\proclaim{11.6 Proposition}  Let $f:Y\to X$ be a proper morphism, $x\in X$ a closed
point, $F=f^{-1}(x)$ and $G$ a locally free sheaf on $Y$.  Assume  
that
$Y$ is CM of pure dimension $n$. 

(11.6.1) If $\supp R^if_*G\subset \{x\}$ then
$R^if_*G \dual H^{n-i}_F(\omega_Y\otimes G^{-1})$.

(11.6.2) If $\supp R^if_*G\subset \{x\}$ for every $i\geq 0$ then
$R^jf_*G \dual R^{n-j}f_*(\omega_Y\otimes G^{-1})$ for every $j\geq 0$.

(11.6.3) Assume that
 $R^if_*(\omega_Y\otimes G^{-1})=0$ for $i>0$. If
 $\supp R^{n-i}f_*G=\{x\}$ for some $i$ then
$$
H^i_x(X,f_*(\omega_Y\otimes G^{-1}))\dual R^{n-i}f_*G.
$$
\endproclaim

\demop By  duality (on $Y$) we obtain that
$\ext^{n-i}(\o_{mF},\omega_Y\otimes G^{-1})$ is dual to $H^i(\o_{mF}\otimes
G)$ for every
$m\geq 1$. The inverse limit of the $H^i(\o_{mF}\otimes G)$ is the
completion of
$R^if_*G$ at
$x$ which is   finite
dimensional by assumption.  
$H^{n-i}_F(Y,\omega_Y\otimes G^{-1})$ is the direct limit of the groups
$\ext^{n-i}(\o_{mF},\omega_Y\otimes G^{-1})$, and is therefore finite
dimensional by the above duality. 

Thus for $m\gg 1$ we obtain that 
$$
H^{n-i}_F(Y,\omega_Y\otimes G^{-1})=
\ext^{n-i}(\o_{mF},\omega_Y\otimes G^{-1})
\dual H^i(\o_{mF}\otimes G)
=R^if_*G.
$$

In order to show (11.6.2) assume for simplicity that $\omega_Y$ is locally
free. (This is the only case that we use later.)
Then $R^{n-j}f_*(\omega_Y\otimes G^{-1})$ is dual to 
$H^j_F(\omega_Y\otimes (\omega_Y\otimes G^{-1})^{-1})=
H^j_F( G)$ and  $H^j_F( G)=R^jf_*G$ by (11.2.4).

Finally consider (11.6.3). By (11.2.3) and (11.6.1) we see that  
$$
H^i_x(X,f_*(\omega_Y\otimes
G^{-1}))\cong H^i_F(Y,\omega_Y\otimes G^{-1})\dual R^{n-i}f_*G.\qed
  $$
\enddemo

The next two applications use these duality results to get information about the
depth of direct image sheaves.

\proclaim{11.7 Corollary} (cf. \cite{Fujita85}) Let $f:Y^n\to X^k$ be a proper
morphism between pure dimensional schemes, $Y$ CM. Assume that every irreducible
component of $Y$ dominates an irreducible component of $X$. Let $G$ be a locally
free sheaf on $Y$ such that  $R^if_*(\omega_Y\otimes G^{-1})=0$ for $i>n-k$.

Then $f_*G$ is $S_2$.
\endproclaim

\demop Pick $x\in X$ such that $j=\dim x\leq k-2$. We need to prove that
$\depth_xf_*G\geq 2$. By localization we are reduced to the case when 
$f:Y^{n-j}\to X^{k-j}$ is proper and $x\in X$ is closed.

By (11.2.2) there is an injection $H^1_x(X,f_*G)\hookrightarrow H^1_F(Y,G)$
and by (11.6.1) $H^1_F(Y,G)\dual R^{n-j-1}f_*(\omega_Y\otimes G^{-1})$.
Since  $n-j-1\geq n-(k-2)-1>n-k$, the latter group is zero by assumption.
\qed\enddemo

The following result is a refined version of (11.7). It is stated in the dual
form, since we use it mostly that way.

\proclaim{11.8 Corollary} Let $f:Y^n\to X^k$ be a proper morphism of pure
dimensional schemes, $Y$ CM. Assume that every irreducible
component of $Y$ dominates an irreducible component of $X$.  Let $G$ be a
locally free sheaf on $Y$ such that
 $R^if_*(\omega_Y\otimes G^{-1})=0$ for $i>0$. 
The following are equivalent:

(11.8.1) $R^if_*G=0$ for every $i>n-k$, and

(11.8.2) $f_*(\omega_Y\otimes G^{-1})$ is a CM sheaf.
\endproclaim

\demop      There is nothing to prove if $k=0$.
The assumptions are stable under localization in $X$, thus by
induction on $k$ we may assume that there is a closed point
$x\in X$ such that $\supp R^if_*G\subset \{x\}$ for every $i>n-k$
and $f_*(\omega_Y\otimes G^{-1})$ is a CM sheaf on $X-x$. 

Then $f_*(\omega_Y\otimes G^{-1})$ is a CM sheaf iff
$H^i_x(X,f_*(\omega_Y\otimes G^{-1}))=0$ for $i<k$ which, by (11.6.3), is
equivalent to $R^if_*G=0$ for every $i>n-k$.\qed\enddemo

\proclaim{11.9 Corollary} \cite{KKMS73, p.50} 
Let $X$ be a  normal scheme   and $f:Y\to X$ a resolution of singularities.
Assume that  $R^if_*\omega_Y=0$ for $i>0$.
The  following  conditions are equivalent:

(11.9.1) $R^if_*\o_Y=0$ for $i>0$. 

(11.9.2) $f_*\omega_Y=\omega_X$ and $\omega_X$ is a CM sheaf.

(11.9.3) $f_*\omega_Y=\omega_X$ and $X$ is  CM.
\endproclaim

\demop  (11.9.2) $\Leftrightarrow$ (11.9.3) was established in (11.5). 
$f_*\omega_Y$ is a subsheaf of $\omega_X$ and they are equal in codimension
one. Thus (11.9.2) is equivalent to the condition:  $f_*\omega_Y$ is a CM sheaf.
 (11.8) for
$G=\o_Y$ shows that the latter is equivalent to (11.9.1). 
\qed\enddemo

\demo{11.10 Definition} Let $X$ be an excellent scheme over a field of
characteristic zero. We say that
$X$ has {\it rational singularities} if 
 it
satisfies the     equivalent conditions of (11.9).
(By (2.17.6)   $R^if_*\omega_Y=0$ for every $i>0$ and for every  resolution of
singularities $f:Y\to X$ if $X$ is over a field of characteristic zero.)
\enddemo

\demo{11.11 Exercise} Let $X$ be a reduced and pure dimensional  scheme,
$f:Y\to X$ a resolution of singularities and $g:Z\to X$ a proper birational
morphism, $Z$ normal.  Show that there are natural inclusions
$$
f_*\omega_Y\subset f_*\omega_Z\subset \omega_X.
$$

In particular,  the conditions (11.9.2)  and (11.9.3) are independent of the 
choice of $f:Y\to X$.
\enddemo

The following result, due to \cite{Fujita85; KaMaMa87, 1-3}, is the main
technical result of the section.  We give a   simpler proof in a slightly more
general form:

\proclaim{11.12 Theorem}  Let $f:Y^n\to X^k$ be a proper morphism of pure
dimensional schemes, $Y$ CM. Assume that every irreducible
component of $Y$ dominates an irreducible component of $X$.  Let $L_1,L_2$ be line bundles on $Y$ and $E$
an effective Cartier divisor on $Y$. Assume that

(11.12.1)  $\codim (f(E),X)\geq 2$,

(11.12.2) $\omega_Y\cong L_1\otimes L_2\otimes \o_Y(E)$, and

(11.12.3) $R^if_*L_j(E)=0$ for $i>0$ and $j=1,2$.

\noindent Then $R^if_*L_j=0$ for $i>0$ and $j=1,2$.
\endproclaim

\demop The assumptions are stable under localization at a point of $X$. Thus by
induction on $\dim X$ we may assume that there is a closed point $x\in X$
such that $\supp R^if_*L_j\subset
\{x\}$ for $i>0$ and $j=1,2$.

The main part of the proof is to  establish two different dualities between the
sheaves 
$R^if_*L_j$. Note that $\omega_Y\otimes L_j^{-1}\cong L_{3-j}(E)$.

By (11.6.1) $R^if_*L_j$ is dual to $H^{n-i}_F(Y,L_{3-j}(E))$  and using
(11.2.3) we obtain that
$H^{n-i}_F(Y,L_{3-j}(E))=H^{n-i}_x(X,f_*L_{3-j}(E))$.

$f_*L_{3-j}$ is $S_2$ by (11.7), hence  $f_*L_{3-j}=f_*L_{3-j}(E)$ 
and so  $H^{n-i}_x(X,f_*L_{3-j}(E))=H^{n-i}_x(X,f_*L_{3-j})$ for every
$i,j$.  $H^{n-i}_x(X,f_*L_{3-j})=0$ for $i=n,n-1$ since $f_*L_{3-j}$ is $S_2$,
and this shows that
$$
R^nf_*L_j=R^{n-1}f_*L_j=0.\tag 11.12.4
$$
By (11.6.1) $H^{n-i}_F(Y,L_{3-j})$ is dual to $R^if_*L_j(E)=0$ for $i\geq 1$.
Therefore 
by (11.2.5)
$H^{n-i}_x(X,f_*L_{3-j})=R^{n-i-1}f_*L_{3-j}$ for $1\leq i\leq n-2$. 

Putting all these together we obtain that
$$
R^if_*L_j\dual R^{n-i-1}f_*L_{3-j}\qtq{for $1\leq i\leq n-2$.}\tag 11.12.5
$$

On the other hand, look at the exact sequence
$$
0\to L_j\to L_j(E)\to L_j(E)|E\to 0.
$$
By assumption $R^if_*L_j(E)=0$ for $i>0$ and we proved that 
$f_*L_j=f_*L_j(E)$. Thus $R^if_*( L_j(E)|E)=R^{i+1}f_* L_j$ for $i\geq 0$.
In particular, $\supp R^if_*( L_j(E)|E) \subset
\{x\}$ for $i\geq 0$.

By adjunction  $\omega_E\otimes ( L_1(E)|E)^{-1}\cong   L_2(E)|E$. Thus
by (11.6.2) $R^if_*( L_j(E)|E)\dual R^{n-1-i}f_*( L_{3-j}(E)|E)$. 
This gives that
$$
R^if_*L_j\dual R^{n-i+1}f_*L_{3-j}\qtq{for $1\leq i\leq n$.}\tag 11.12.6
$$

Put (11.12.5--6) together  to conclude that
$$
R^if_*L_j\cong R^{i-2}f_*L_j \qtq{for $3\leq i\leq n$.}\tag 11.12.7
$$
Starting with the vanishing (11.12.4) this completes the proof by descending
induction on $i$.
\qed\enddemo

The simplest application  of this vanishing is the first part of (11.1):

\proclaim{11.13 Corollary} Let $X$ be an excellent normal scheme over a field
of characteristic zero. Assume that $\omega_X$ exists and is locally free. 
  Then
$X$ has rational singularities iff $X$ has canonical singularities.
\endproclaim

\demop Let $f:Y\to X$ be a resolution of singularities. 
Assume that $X$ has rational singularities. Then $f_*\omega_Y=\omega_X$, hence
there is a natural map $f^*\omega_X\to \omega_Y$. This shows that
$\omega_Y=f^*\omega_X(E)$ for some effective divisor $E\subset Y$.

Conversely, assume that $X$ has canonical singularities, that is
$\omega_Y\cong f^*\omega_X\otimes \o_Y(E)$ for some effective divisor $E$.
Apply (11.12) with $L_1\cong \o_Y$ and $L_2\cong f^*\omega_X$. 

Then $L_2(E)\cong \omega_Y$, hence $R^if_*L_2(E)=0$ for $i>0$ by (2.17.6).
By the projection formula,
$R^if_*L_1(E)=\omega_X^{-1}\otimes R^if_*L_2(E)=0$.
Thus by (11.12) we see that $R^if_*\o_Y=R^if_*L_1=0$ for $i>0$.\qed\enddemo

If $\omega_X$ is not locally free, then   we show that every klt singularity is
rational. The sharpest technical result proved in \cite{Fujita85; KaMaMa87, 1-3}
asserts that if
$(X,D)$ is dlt then $X$ has rational singularities. We have not defined 
dlt (cf. \cite{Koll\'ar et al.92, 2.13}), but the proof requires only small
changes.

\proclaim{11.14 Corollary} Let $(X,D)$ be a klt pair over a
field of characteristic zero. 
  Then
$X$ has rational singularities.
\endproclaim

\demop Let $f:Y\to X$ be a log resolution.  Write
$$
K_Y\equiv f^*(K_X+D)+F, \qtq{where $\rup{F}$ is effective.}
$$
Set $\Delta=\rup{F}-F$; this is an effective  normal crossing divisor  such
that $\rdown{\Delta}=0$.  Set 
$$
L_1=\o_Y, \quad E=\rup{F}\qtq{and}L_2=\o(K_Y-E).
$$
  Then
$$
L_1(E)\equiv K_Y+\Delta-f^*(K_X+D),\qtq{and}
L_2(E)\equiv K_Y.
$$
$-f^*(K_X+D)$ is $f$-nef (even $f$-numerically trivial), thus (2.17.3)
applies and we get that 
$R^if_*L_j(E)=0$ for $i>0$ and $j=1,2$.
Thus by (11.12) we see that $R^if_*\o_Y=R^if_*L_1=0$ for $i>0$.\qed\enddemo

\cite{Elkik78} proved that a small deformation of  a rational singularity is
again rational. The following is a slight variation of her arguments,
which also says something about deformations of certain nonnormal schemes.

\proclaim{11.15 Theorem}  Let $X$ be an excellent  normal scheme over a field of
characteristic zero and $h:X\to T$   a flat morphism to the spectrum of a DVR
with closed point $0\in T$ and local parameter $t\in \o_T$. Set $X_0\deq
h^{-1}(0)$.  Let $\pi:\bar X_0\to X_0$ be the normalization. Assume that 
$\pi$ is an isomorphism in codimension one and that $\bar X_0$ has rational
singularities.

Then $\bar X_0\cong X_0$ and $X$ has rational singularities.
\endproclaim

\demop  Take a resolution of singularities $f:X'\to X$ such that
$X'_0+E\deq (f\circ h)^{-1}(0)$ is a divisor with normal crossings only
where $f:X'_0\to \bar X_0\to X_0$ is a resolution of singularities.

We denote $h^*t$ and $f^*h^*t$ again by $t$.  On $X'$ we have an exact sequence
$$
0\to \omega_{X'}@>t>>\omega_{X'}@>>>\omega_{X'_0+E}@>>> 0.
\tag 11.15.1
$$
On $X$ we get an exact sequence
$$
0\to \omega_{X}@>t>>\omega_{X}@>>>\tilde{\omega}_{X_0}@>>> 0,
\tag 11.15.2
$$
where $\tilde{\omega}_{X_0}\subset \omega_{X_0}$
(see e.g. \cite{Reid94, 2.13}).
For any pure dimensional scheme $Z$, $\omega_Z$ is $S_2$
(see e.g. \cite{Reid94, 2.12}), which  in particular  implies that
if $\pi:\bar Z\to Z$ is finite and an isomorphism in codimension one then
$\pi_*\omega_{\bar Z}\cong \omega_Z$.
Since $\bar X_0$ has rational singularities, this implies that 
$\omega_{X_0}=f_*\omega_{X'_0}$

We have a natural injection $f_*\omega_{X'}\hookrightarrow \omega_X$. This gives
the following commutative diagram
$$
\CD
0@>>> f_*\omega_{X'}@>t>> f_*\omega_{X'} @>>> f_*\omega_{X'_0+E}@>>> 
R^1f_*\omega_{X'}\\ 
@.\bigg\vert @.\bigg\vert @. @AAA\\
\vspace{-6pt}
(11.15.3) @. \bigg\vert @.\bigg\vert @. \omega_{X_0}=f_*\omega_{X'_0} @. \\
\vspace{-6pt}
@. @VVV @VVV @AAA \\
0@>>> \omega_{X} @>t>>\omega_{X} @>>>\tilde{\omega}_{X_0} @>>> 0
\endCD
%\tag 11.15.3
$$
$R^1f_*\omega_{X'}=0$ by (2.17.6), thus
$f_*\omega_{X'} @>>> f_*\omega_{X'_0+E}$ is surjective. 
The vertical maps
$$
\tilde{\omega}_{X_0}\to \omega_{X_0}=f_*\omega_{X'_0}\to f_*\omega_{X'_0+E}
$$
 are injections.
 By the commutativity of (11.15.3) we obtain that they are both
isomorphisms:
$$
\tilde{\omega}_{X_0}= f_*\omega_{X'_0}= f_*\omega_{X'_0+E}.\tag 11.15.4
$$
This implies that
$$
f_*\omega_{X'}\to \omega_X\to \omega_X/t\omega_X
$$
is surjective, hence by the  Nakayama lemma
$f_*\omega_{X'}\hookrightarrow \omega_X$ is an isomorphism.
Also by (11.15.4), $\omega_X/t\omega_X\cong f_*\omega_{X'_0}$. By (11.4.2)
we know that $f_*\omega_{X'_0}$ is a CM-sheaf, hence $\omega_X$ is also a
CM-sheaf. By (11.9) we conclude that $X$ has rational singularities.

 Therefore $X$ is CM and so is $X_0$.
This means that $X_0=\bar X_0$.\qed\enddemo

\end